\documentclass[a4paper,11 pt]{article}

\usepackage{amsmath}
\usepackage{amssymb}
\usepackage[english]{babel}
\usepackage{}
\newtheorem{thm}{Theorem}[section]

\newtheorem{remark}[thm]{Remark}
\newtheorem{remarks}[thm]{Remarks}

\newtheorem{defin}[thm]{Definition}
 \numberwithin{equation}{section} \textheight22cm\textwidth16.5cm
\newcommand\mi{\ensuremath{\ds{-\hskip-3.9mm\int}}}

\def\R{\mathbb R}

\def\epsilon{\varepsilon}
\def\o{\overline\Omega}

\def\te{\tilde{e}}
\def\ds{\displaystyle}

\begin{document}

\title{\textbf{Pulsating travelling fronts:\\ Asymptotics and
homogenization regimes}} \author{\begin{Large}Mohammad El
Smaily\end{Large}\footnote{E-mail:
mohammad.el-smaily@etu.univ-cezanne.fr}\\ Universit\'e Aix-Marseille
III, LATP, Facult\'e des Sciences et
Techniques,\\
 Avenue Escadrille Normandie-Niemen, F-13397 Marseille
Cedex 20, France. }

\date{28 March 2008}
\maketitle

 \textbf{$\mathcal{A}$bstract}. This paper is concerned with some nonlinear
propagation phenomena for reaction-advection-diffusion equations
with Kolmogrov-Petrovsky-Piskunov (KPP) type nonlinearities in
general periodic domains or in infinite cylinders with oscillating
boundaries. Having a variational formula for the minimal speed of
propagation involving eigenvalue problems ( proved in Berestycki,
Hamel and Nadirashvili \cite{BHN1}), we consider the minimal speed
of propagation as a function of diffusion factors, reaction factors
and periodicity parameters. There we study the limits, the
asymptotic behaviors and the variations of the considered functions
with respect to these parameters. Section \ref{homogenization} deals with
homogenization problem as an application of the results in the
previous sections in order to find the limit of the minimal speed
when the periodicity cell is very small.
\vskip 0.4 cm
\section{Introduction}
 This paper is a continuation in the
study of the propagation phenomena of pulsating travelling fronts in
a periodic framework corresponding to reaction-advection-diffusion
equations with heterogenous KPP (Kolmogrov, Petrovsky and Piskunov)
nonlinearities. We will precisely describe the
\textsl{heterogenous-periodic} setting, recall the extended notion
of \textsl{pulsating travelling fronts}, and then we move to announce
the main results. Let us first recall some of the basic
features of the \textsl{homogenous} KPP equations.

Consider the Fisher-KPP equation:
\begin{equation}\label{homeqin RN}
    u_t\,-\Delta u\,=\,f(u)\quad\hbox{in}\;\mathbb{R}^{N}.
\end{equation}
It was introduced in the celebrated papers of Fisher (1937) and
in \cite{KPP} originally motivated by models in biology. Here, the
main assumption is that $f$ is, say, a $C^1$ function satisfying
\begin{eqnarray}\label{KPPfunction}
   \left\{
     \begin{array}{ll}
       f(0)\,=\,f(1)=0, \;f'(1)<0,\; f'(0)>0, & \hbox{} \\
       \;f>0\;\hbox{in}\;(0,1),\;
   f<0\;\hbox{in}\,
   (1,+\infty), & \hbox{}
     \end{array}
   \right.
\end{eqnarray}
\begin{equation}
\, f(s)\leq f'(0)s \,,\forall s\in \,[0,1].
\end{equation}
  As examples of such nonlinearities, we have: $f(s)=s(1-s)$ and
$f(s)=s(1-s^2).$

 The important feature in (\ref{homeqin RN}) is that this equation has
a family of planar travelling fronts. These are solutions of the
form
\begin{eqnarray}\label{u=phi1}
\left\{
\begin{array}{l}
 \forall(t,x)\,\in\mathbb{R}\times\mathbb{R}^N,~u(t,x)\,=\,\phi(x\cdot e+ct),\\
    \phi(-\infty)=0\quad\hbox{and}\quad\phi(+\infty)=1,
    \end{array}\right.
\end{eqnarray}
where $e\in\mathbb{R}^N$ is a fixed vector of unit norm which is the
direction of propagation, and $c\,>\,0$ is the speed of the front.
The function $\phi:\,\mathbb{R}\mapsto\mathbb{R}$ satisfies
\begin{eqnarray}\label{u=phi2}
\left\{
\begin{array}{ll}-\phi^{''}+c\,\phi\,=f(\phi),\\
\phi(-\infty)=0\quad\hbox{and}\quad\phi(+\infty)=1.
\end{array}
\right.
\end{eqnarray}
In the original paper of Kolmogorov, Petrovsky and Piskunov, it was
proved that, under the above assumptions, there is a threshold value
$c^{*} = 2\sqrt{f^{'}(0)} > 0$ for the speed $c$. Namely, no fronts
exist for $c \,<\, c^{*},$ and, for each $c \,\geq\, c^{*},$ there
is a unique front of the type (\ref{u=phi1}-\ref{u=phi2}).
Uniqueness is up to shift in space or time variables.

Later, the homogenous setting was extended to a general
heterogenous periodic one. The heterogenous character appeared both
in the reaction-advection-diffusion equation and in the underlying
domain. The general form of these equations is
\begin{eqnarray}\label{heteq}
    \left\{
      \begin{array}{l}
 u_t =\nabla\cdot(A(z)\nabla u)\;+q(z)\cdot\nabla u+\,f(z,u),\; t\,\in\,\mathbb{R},\;z\,\in\,\Omega,
 \\
        \nu \cdot A\;\nabla u(t,z) =0,\;
        t\,\in\,\mathbb{R},\; z\,\in\,\partial\Omega,
      \end{array}
    \right.
\end{eqnarray}
where $\nu(z)$ is the unit outward normal on $\partial\Omega$ at the
point $z.$

 The propagation phenomena attached with equation
(\ref{heteq}) has been widely studied in many papers. Several properties
of pulsating fronts in periodic media and their speed of propagation
were given in several papers ( Berestycki, Hamel \cite{BH1},
Berestycki, Hamel, Nadirashvili \cite{BHN1}, and Berestycki, Hamel,
Roques \cite{BHR1, BHR2} and Xin \cite{JXin ejam}). In section \ref{periodic
framework}, we will recall the periodic framework and some known
results which motivate our study. The main results of this paper
are presented in sections \ref{as eps tends to zero} to \ref{variation of min speed w.r.t diffusion and period}.

\section{The periodic framework}\label{periodic framework}
\subsection{Pulsating travelling fronts in periodic domains}
In this section, we introduce the general setting with the precise
assumptions. Concerning the domain, let $N\geq\,1$ be the space
dimension, and let $d\;$ be an integer so that $1\leq\,d\leq\,N.$
For an element
$z=(x_1,x_2,\cdots,x_d,x_{d+1},\cdots,x_N)\in\,\mathbb{R}^{N},$ we
call $x=(x_1,x_2,\cdots,x_d)$ and $
y=(x_{d+1},\cdots,x_N)$ so that $z=(x,y).$ Let
$L_1,\cdots,L_d$ be $d$ positive real numbers, and let $\Omega$ be a
$C^3$ nonempty connected open subset of $\mathbb{R}^N$ satisfying
\begin{eqnarray}\label{comega}
    \left\{
      \begin{array}{l}
        \exists\,R\geq0\,;\forall\,(x,y)\,\in\,\Omega,\,|y|\,\leq\,R,  \\
        \forall\,(k_1,\cdots,k_d)\;\in\;L_1\mathbb{Z}\times\cdots\,\times L_d\mathbb{Z},
        \quad\displaystyle{\Omega\;=\;\Omega+\sum^{d}_{k=1}k_ie_i},
      \end{array}
    \right.
\end{eqnarray}
where $\;(e_i)_{1\leq i\leq N}\;$ is the canonical basis of
$\mathbb{R}^N.$
 In particular, since $d\geq1,$ the set $\Omega$ is \textbf{unbounded}.\\
In this periodic situation, we give the following definitions:
\begin{defin}[Periodicity cell] The set $C=\{\,(x,y)\in\Omega
;\; x_{1}\in(0,L_1),\cdots,x_{d}\in(0,L_d)\}$ is called the
periodicity cell of $\Omega.$
\end{defin}
\begin{defin}[$L$-periodic flows ] A field
$w:\Omega\rightarrow\,\mathbb{R}^N$ is said to be $L$-periodic with
respect to $x$
 if $w(x_1+k_1,\cdots,x_d+k_d\,,y)=w(x_1,\cdots,x_d,y)$ almost everywhere in
$ \Omega,$ and for all
$\displaystyle{k=(k_1,\cdots\,,k_d)\in\prod^{d}_{i=1}L_i\mathbb{Z}}.$
\end{defin}

Before going further on, we point out that this framework includes
several types of simpler geometrical configurations. The case of the
whole space $\mathbb{R}^{N}$ corresponds to $d = N,$ where
$L_1,\ldots,L_N$ are any positive numbers. The case of the whole
space $\mathbb{R}^{N}$ with a periodic array of holes can also be
considered. The case d = 1 corresponds to domains which have only
one unbounded dimension, namely infinite cylinders which may be
straight or have oscillating periodic boundaries, and which may or
may not have periodic holes. The case $2\leq d \leq N-1$ corresponds
to infinite slabs.

We are concerned with propagation phenomena for the
reaction-advection-diffusion equation (\ref{heteq}) set in the
periodic domain $\Omega.$ Such equations arise in combustion models
for flame propagation (see \cite{ronney}, \cite{Williams} and \cite{zeldovich}), as well as in models in biology and for
population dynamics of a species (see \cite{fife}, \cite{KKTS}, \cite{murray} and \cite{shigesada Kawsaki 1}). These equations are used in modeling
the propagation of a flame or of an epidemics in a periodic heterogenous medium.
The passive quantity $u$ typically stands for
the temperature or a concentration which diffuses in a periodic
excitable medium. However, in some sections we will ignore the
advection and deal only with reaction-diffusion equations.

 Let us now detail the assumptions concerning the coefficients in
 (\ref{heteq}). First, the diffusion matrix $A(x,y)=(A_{ij}(x,y))_{1\leq i,j\leq
 N}$ is a \textit{symmetric} $C^{2,\delta}(\,\overline{\Omega}\,)$ (with $\delta
\,>\,0$) matrix field satisfying
\begin{eqnarray}\label{cA}
    \left\{
      \begin{array}{l}
        A\; \hbox{is $L$-periodic with respect to}\;x, \vspace{3 pt}\\
        \exists\,0<\alpha_1\leq\alpha_2,\forall(x,y)\;\in\;\Omega,\forall\,\xi\,\in\,\mathbb{R}^N,\vspace{3 pt}\\
       \displaystyle{ \alpha_1|\xi|^2 \;\leq\;\sum_{1\leq i,j\leq N}\,A_{ij}(x,y)\xi_i\xi_j\,\;\leq\alpha_2|\xi|^2.}
      \end{array}
    \right.
\end{eqnarray}
The boundary condition $\nu\cdot A\nabla u(x,y) =0$ stands for
$\displaystyle{\sum_{1\,\leq\, i,j\leq
N}\nu_i(x,y)A_{ij}(x,y)\partial_{x_j}u(t,x,y)},$ and $\nu$ stands
for the unit outward normal on $\partial\Omega.$ We note that when
$A$ is the identity matrix, then this boundary condition reduces to
the usual Neumann condition $\partial_{\nu}u\,=\,0.$

The underlying advection $q(x,y)=(q_1(x,y),\cdots,q_N(x,y))$  is a
$C^{1,\delta}(\overline{\Omega})$ (with $\delta>0$) vector field
satisfying
\begin{equation}\label{cq}
    \left\{
      \begin{array}{ll}
        q\quad\hbox{is $L-$ periodic with respect to }\;x, & \hbox{} \\
        \nabla\cdot q\,=0\quad\hbox{in}\; \overline{\,\Omega\,}, \\
         q\cdot\nu\,=0\quad \hbox{on}\;\partial\Omega\,, \\
        \forall\,1\leq\,i\,\leq\,d, \quad\displaystyle{\int_{C}q_i\;dx\,dy \,=\,0}\hbox{.}
      \end{array}
    \right.
\end{equation}

Concerning the nonlinearity, let $f\,=\,f(x,y,u)$ be a nonnegative
function defined in $\overline{\Omega}\,\times[0,1],\;$ such that
\begin{eqnarray}\label{cf1}
    \left\{
      \begin{array}{ll}
        f\geq0, f\;\hbox{is $L$-periodic with respect to}\; x, \;\hbox{and of class}\;C^{1,\,\delta}(\overline{\Omega}\,\times\,[0,1]),\vspace{3 pt}  \\
        \forall\,(x,y)\,\in\,\overline{\Omega},\quad \displaystyle{f(x,y,0)=\,f(x,y,1)=\,0 } \hbox{,} \vspace{3 pt}\\
        \exists \,\rho\in\,(0,1),\;\forall(x,y)\,\in\overline{\Omega},\;\displaystyle{\forall\, 1-\rho\leq\,s\,\leq
s'\,\leq\,1,}\;
\displaystyle{f(x,y,s)\;\geq\,f(x,y,s') } \hbox{,} \vspace{3 pt}\\
        \forall\,s\in\,(0,1),\; \exists \,(x,y)\in\overline{\Omega}\;\;\hbox{such that}\;f(x,y,s)\,>\,0  \hbox{,} \\
        \forall\,(x,y)\,\in\,\overline{\Omega},\quad \displaystyle{f'_{u}(x,y,0)\,=\,\lim_{u\rightarrow\,0^+}\frac{f(x,y,u)}{u}\,>\,0}  \hbox{,}
      \end{array}
    \right.
\end{eqnarray}
 with the additional assumption
 \begin{equation}\label{cf2}
 \forall\, (x,y,s)\in\overline{\Omega}\times(0,1),\quad 0\,<\,f(x,y,s)\,\leq\,f'_u(x,y,0)\times\,s
 .
\end{equation}
 We denote by $\zeta(x,y):=\,f'_u(x,y,0),$ for each
$(x,y)\in\overline{\Omega}.$\\
The set of such nonlinearities contains two particular types of
functions: \begin{itemize}
\item The homogeneous  (KPP) type: $f(x,y,u)=g(u),$ where $g$
is a $C^{1,\delta}$ function that satisfies:
\begin{center}$g(0)=g(1)=\,0,\;g>0$ on $(0,1),\;g'(0)>0,\, g'(1)<0$ and
$0<g(s)\leq g'(0)s$ in $(0,1).$\end{center}
\item Another type of such nonlinearities consists
of functions $f(x,y,u)=h(x,y).\tilde{f}(u),$ such that $\tilde{f}$
is of the previous type, while $h$ lies in $
C^{1,\delta}(\overline{\Omega}), L$ -periodic with respect to $x,$
and positive in $\overline{\Omega}.$
\end{itemize}

Having this periodic framework, the notions of travelling fronts and
propagation were extended, in \cite{BH1}, \cite{BHN1}, \cite{KKTS},
\cite{Pap Xin} \cite{shigesada Kawsaki 1}, \cite{shigesada Kawsaki
2}, and \cite{Xin3} as follows:
\begin{defin}
Let $e\,=(e^1,\cdots,e^d)$ be an arbitrarily given vector in
$\mathbb{R}^d.$ A function $u=u(t,x,y)$ is called a pulsating
travelling front propagating in the direction of $e$ with an
effective speed $c\neq0,$ if $u$ is a classical solution of
\begin{eqnarray}\label{front}
    \left\{
      \begin{array}{ll}
u_t =\nabla\cdot(A(x,y)\nabla u)+q(x,y)\cdot\nabla u+\,f(x,y,u),\;
t\,\in\,\mathbb{R},\;(x,y)\,\in\,\Omega,
\vspace{3 pt}\\
        \nu \cdot A\;\nabla u(t,x,y) =0,\;
        t\,\in\,\mathbb{R},\;(x,y)\,\in\,\partial\Omega,\vspace{3 pt}\\
 \displaystyle{\forall\, k\in\prod^{d}_{i=1}L_i\mathbb{Z},\; \forall\,(t,x,y)\,\in\,\mathbb{R}\,\times\,\overline{\Omega}},
\quad\displaystyle{u(t-\frac{k\cdot e}{c},x,y)\,=\,u(t,x+k,y)}  \hbox{,}\vspace{3 pt} \\
         \displaystyle{\lim_{x\cdot e\,\rightarrow\,-\infty}u(t,x,y)\,=0,\;\hbox{and}\; \lim_{x\cdot e\rightarrow \,+\infty} u(t,x,y)\,=\,1}  \hbox{,}
         \vspace{3 pt}\\
          0\,\leq\,u\leq\,1,
      \end{array}
    \right.
\end{eqnarray}
 where the above limits hold locally in $t$ and uniformly in $y$ and in the
directions of $\mathbb{R}^d$ which are orthogonal to $e$ .
\end{defin}
\subsection{Some important known results concerning the propagation phenomena in a periodic
framework}\label{recall}
 Under the assumptions (\ref{comega}),
(\ref{cA}), (\ref{cq}), (\ref{cf1}) and (\ref{cf2}) set in the
previous subsection, Berestycki and Hamel \cite{BH1} proved that:
having a pre-fixed unit vector $e\in\mathbb{R}^{d},$ there exists
$c^{*}(e)>0$ such that pulsating travelling fronts propagating in the
direction $e$ (i.e satisfying (\ref{front})) with a speed of
propagation $c$ exist if and only if $c\geq c^{*}(e);$ moreover, the
pulsating fronts (within a speed $c\geq c^{*}(e)$) are increasing in
the time $t.$ The value
$\displaystyle{c^{*}(e)=c^{*}_{\Omega,A,q,f}(e)}$ is called the
\emph{minimal speed of propagation in the direction of} $e.$ Other
nonlinearities have been considered in the cases of the whole space
$\mathbb{R}^{N}$ or in the general periodic framework (see
\cite{BH1}, \cite{shigesada Kawsaki 1}, \cite{shigesada Kawsaki 2},
\cite{Xin1}, \cite{Xin2}, \cite{Xin3}, \cite{Xin4}).

Having the threshold value $\displaystyle{c^{*}_{\Omega,A,q,f}(e)},$
our paper aims to study the limits, the asymptotic behaviors, and
the variations of some parametric quantities. These parametric
quantities involve the parametric speeds of propagation of different
reaction-advection-diffusion problems within a diffusion factor
$\varepsilon>0,$ a reaction factor $B>0,$ or a periodicity parameter
$L.$ Thus, it is important to have a variational characterization
which shows the dependance of the minimal speed of propagation on
the coefficients $A,\;q$ and $f$ and on the geometry of the domain
$\Omega.$ In this context, Berestycki, Hamel, and Nadirashvili
\cite{BHN1} gave such a formulation for $c^{*}_{\Omega,A,q,f}(e)$
involving elliptic eigenvalue problems. We recall this variational
characterization in the following theorem:
\begin{thm}[Berestycki, Hamel, and Nadirashvili
\cite{BHN1}]\label{varthm}
 Let $e$ be a fixed unit vector in
$\mathbb{R}^d.$ Let $\tilde{e}=(e,0,\ldots,0)$ $\in\mathbb{R}^N.$
Assume that $\Omega,\,A$ and $f$ satisfy (\ref{comega}),(\ref{cA}),
(\ref{cf1}), and (\ref{cf2}). The minimal speed
$c^{*}(e)=c^{*}_{\Omega,A,q,f}(e)$ of pulsating fronts solving
(\ref{front}) and propagating in the direction of $e$ is given by
\begin{equation}\label{var}
    \displaystyle{c^*(e)=c^{*}_{\Omega,A,q,f}(e)\,=\,\min_{\lambda>0}\frac{k(\lambda)}{\lambda}},
\end{equation}
where $\displaystyle{k(\lambda)=k_{\Omega,e,A,q,\zeta}(\lambda)}$ is
the principal eigenvalue of the operator
$\displaystyle{L_{\Omega,e,A,q,\zeta,\lambda}}$ which is defined by
\begin{eqnarray}\label{Leq}
\begin{array}{ll}
\displaystyle{L_{\Omega,e,A,q,\zeta,\lambda}\psi\,:=}&\displaystyle{\,\nabla\cdot(A\nabla\psi)\,-2\lambda\tilde{e}\cdot
A\nabla\psi\,+q\cdot\nabla\psi\,}\vspace{3 pt}\\
&\displaystyle{+[\lambda^2\tilde{e}A\tilde{e}-\lambda\nabla\cdot(A\tilde{e})-\lambda
q\cdot\tilde{e}+\,\zeta]\psi}
\end{array}
\end{eqnarray}
acting on the set
$$
\begin{array}{ll}
E=&\left\{\right.\psi\in C^2(\overline{\Omega}), \psi\hbox{ is
$L$-periodic with respect to $x$ and}\;
\nu\cdot A\nabla\psi=\lambda(\nu A
\tilde{e}\psi)\;\hbox{on}\;\partial{\Omega}\left.\right\}.
\end{array}
$$
\end{thm}

The proof of formula (\ref{var}) is based on methods developed in
\cite{BH1}, \cite{BNir1} and \cite{BNV}. These are techniques of sub and
super-solutions, regularizing and approximations in bounded domains.

 We note that in formula (\ref{var}), the value of the minimal
speed $c^*(e)$ is given in terms of the direction $e,$ the domain
$\Omega,$ and the coefficients $A,q$ and $f^{\,'}_{u}(.,.,0).$
Moreover, it is important to notice that the dependence of $c^*(e)$
on the nonlinearity $f$ is only through the derivative of $f$ with
respect to $u$ at $u=0.$

 Before going further on, let us mention that formula
(\ref{var}) extends some earlier results about front propagation.
When $\Omega=\mathbb{R}^{N},\,A=Id$ and $f=f(u)$ (with
$f(u)\leq\,f^{\,'}(0)u$ in $[0,1]$), formula (\ref{var}) then
reduces to the well-known KPP formula
$c^*(e)\,=\,2\sqrt{f^{\,'}(0)}.$ That is the value of the minimal
speed of propagation of planar fronts for the \textsl{homogenous}
reaction-diffusion equation:
 $ u_t-\Delta u=f(u)\;\hbox{in}\;\mathbb{R}^{N}.$\footnote{In fact, the uniqueness, up to multiplication by a non-zero real number, of the first eigenvalue function of
  $L_{\mathbb{R}^{N},e,Id,f^{'}(0),\lambda}\psi=k(\lambda)\psi$ together with this particular situation, yield that the principal eigenfunction $\psi$ is constant and $k(\lambda)=\lambda^{2}+f^{\,'}(0)\; \hbox{for all}\;\lambda>0.$ Therefore
  by (\ref{var}), we have $\displaystyle{c^*(e)\,=\,\min_{\lambda>0}\,\left(\lambda+\frac{f^{\,'}(0)}{\lambda}\right)\,=\,2\sqrt{f^{\,'}(0)}}.$ }

 The above variational characterization
of the minimal speed of propagation of pulsating fronts in general
periodic excitable media will play the main role in studying the
dependence of the minimal speed $c^{*}(e)=c^{*}_{\Omega,A,q,f}(e)$
on the coefficients of reaction, diffusion, advection and on the
geometry of the domain. In this context, we have:

\begin{thm}[Berestycki,
Hamel, Nadirashvili \cite{BHN1}]\label{influence of f, A} Under the
assumptions (\ref{comega}), (\ref{cA}), and (\ref{cq}) on
$\Omega,\,A,$ and $q,$ let $f\,=\,f(x,y,u)$ $[$respectively
$g\,=\,g(x,y,u)]$ be a nonnegative nonlinearity satisfying
(\ref{cf1}) and (\ref{cf2}). Let $e$ be a fixed unit vector in
$\mathbb{R}^d,$ where $1\,\leq\,d\,\leq\,N,$

 a) If $f_{u}^{\,'}(x,y,0)\,\leq\,g_{u}^{\,'}(x,y,0)$ for all
$(x,y)\,\in\,\overline{\Omega},$ then
$$c_{\Omega,A,q,f}^*(e)\,\leq\,c_{\Omega,A,q,g}^*(e).$$ Moreover if
$\quad f_{u}^{\,'}(x,y,0)\,\leq\,,\,\not\equiv\,
\,g_{u}^{\,'}(x,y,0)\,$ in $\,\overline{\Omega},$ then
$c_{\Omega,A,q,f}^*(e)\,<\,c_{\Omega,A,q,g}^*(e).$

 b) The map  $B\mapsto c_{\Omega,A,q,Bf}^*(e)$ is increasing in
$B>0$ and
$$\displaystyle{\limsup_{B\,\rightarrow\,+\infty}\displaystyle{\frac{c_{\Omega,A,q,Bf}^*(e)}{\sqrt{B}}}\,<\,+\infty.}$$
Furthermore, if  $\;\Omega\,=\,\mathbb{R}^N\;$ or  if $\;\nu
A\tilde{e}\,\equiv\,0$ on $\,\partial\Omega,$ then
$\;\displaystyle{\liminf_{B\,\rightarrow\,+\infty}\displaystyle{\frac{c_{\Omega,A,q,Bf}^*(e)}{\sqrt{B}}}\,>\,0.}$

c) \begin{equation}\label{influence of A 0}
    \displaystyle{c_{\Omega,A,q,f}^*(e)\,\leq\,||(q.\tilde{e})^{-}||_{\infty}+2\sqrt{\max_{(x,y)\,\in\,\overline{\Omega}}\zeta(x,y)}\sqrt{\max_{(x,y)\,\in\,\overline{\Omega}}\;\tilde{e}A(x,y)\tilde{e}}},
\end{equation}
where
$\displaystyle{||(q.\tilde{e})^{-}||_{\infty}=\max_{(x,y)\in\overline{\Omega}}\left(q(x,y).\tilde{e}\right)^{-}}$
and $s^{-}=\max\,(-s,0)$ for each $s\in\mathbb{R}.$ Furthermore, the
equality holds in (\ref{influence of A 0}) if and only if
$\tilde{e}A\tilde{e}$ and $\zeta$ are constant,
$q.\tilde{e}\equiv\nabla\,.\,(A\tilde{e})\,\equiv\,0\;$ in
$\overline{\,\Omega}$ and $\nu.A\tilde{e}\,=\,0$ on $\partial\Omega$
$($in the case when $\partial\Omega\,\neq\,\emptyset)$.\vskip 0.3cm

  d) Assume furthermore that $f\,=\,f(u)$ and $q\equiv0$ in
$\overline{\Omega},$ then the map $\displaystyle{\beta\mapsto
c_{\Omega,\beta A,0,f}^*(e)}$ is increasing in $\beta>0.$
\end{thm}

 As a corollary of (\ref{influence of A 0}), we see that
$\displaystyle{\limsup_{M\rightarrow+\infty}\frac{c_{\Omega,M
A,q,f}^*(e)}{\sqrt{M}}\,\leq\,C}\;$ where $C$ is a positive
constant. Furthermore, part d) implies that a larger diffusion
speeds up the propagation in the absence of the advection field.

We mention that the existence of pulsating travelling fronts in space-time periodic media was proved in Nolen, Xin \cite{Nolen 1, Nolen 2},
Nolen, Rudd, Xin \cite{Nolen 3} and recently in Nadin \cite{N2,N3}. In \cite{N3}, Nadin  characterized the minimal speed of propagation and he studied the influence of the diffusion, the amplitude of the reaction term and the drift on the characterized speed.

After reviewing some results in the study of the KPP propagation
phenomena in a periodic framework, we pass now to announce new
results concerning the limiting behavior of the minimal speed of
propagation within a small (resp. large) diffusion and reaction
coefficients (in some particular situations of the general periodic
framework) and we will study the minimal speed as a function of the
period of the coefficients in the KPP reaction-diffusion-advection
(or reaction-diffusion) equation in the case where
$\Omega=\mathbb{R}^N$. The proofs will be shown in details in
section~\ref{proofs section}. The announced results will be applied
to find the homogenization limit of the minimal speeds of propagation. We believe
that this limit might help to find the homogenized equation in the
``KPP'' periodic framework (see section \ref{homogenization} for more details).

\section{The minimal speed within small diffusion factors or
within large period coefficients}\label{as eps tends to zero} In
this section, our problem is a reaction-diffusion equation with
absence of advection terms:
\begin{eqnarray}\label{beta A}
\left\{
  \begin{array}{ll}
    u_t =\,\beta\,\nabla\cdot(A(x,y)\nabla u)\;+\,f(x,y,u),\; t\,\in\,\mathbb{R},\;(x,y)\,\in\,\Omega,
    \vspace{3 pt}\\
        \nu \cdot A\;\nabla u(t,x,y) =0,\;
        t\,\in\,\mathbb{R},\;(x,y)\,\in\,\partial\Omega,
  \end{array}
\right.
\end{eqnarray}
where $\beta>0.$

We mention that (\ref{beta A}) is a reaction-diffusion problem
within a diffusion matrix $\beta A.$  Let $e$ be a unit direction in
$\mathbb{R}^d.$ Under the assumptions (\ref{comega}), (\ref{cA}),
(\ref{cf1}) and (\ref{cf2}), for each $\beta>0,$ there corresponds a
minimal speed of propagation $\displaystyle{c_{\Omega,\beta
A,0,f}^*(e)}$ so that a pulsating front with a speed $c$ and
satisfying (\ref{beta A}) exists
 if and only if $c\geq \displaystyle{c_{\Omega,\beta A,0,f}^*(e)}.$

  Referring to part c) of Theorem \ref{influence of f, A}, one gets
$\,0<\, \displaystyle{c_{\Omega,\beta A,0,f}^*(e)}\,\leq\,2\sqrt{\beta}\sqrt{M_0M},$
 for any $\beta>0,$ where
$\displaystyle{\,M_0\,=\,\max_{(x,y)\,\in\,\overline{\Omega}}\;\zeta(x,y)\,}\;$
and
$\displaystyle{\;M\,=\,\max_{(x,y)\,\in\,\overline{\Omega}}\;\tilde{e}A(x,y)\tilde{e}\,}.$\\
Consequently, there exists $C>0$ and independent of $\beta$ such
that
\begin{equation}\label{star}
\forall\beta>0,\,0<\, \displaystyle{\frac{c_{\Omega,\beta
A,0,f}^*(e)}{\sqrt{\beta}}\,\leq\,C}.
\end{equation}

The inequality (\ref{star}) leads us to investigate the limits of
$\displaystyle{\frac{c_{\Omega,\beta A,0,f}^*(e)}{\sqrt{\beta}}}$ as
$\beta\rightarrow0$ and as $\beta\rightarrow +\infty.$ The following
theorem gives the precise limit when the diffusion factor tends to
zero. However, it will not be  announced in the most general
periodic setting. We will describe the situation before the
statement of the theorem:

 The domain will be in the form
$\,\Omega=\mathbb{R}\times\omega\subseteq\,\mathbb{R}^N,$ where
$\,\omega\subseteq\mathbb{R}^{d}\times\mathbb{R}^{N-d-1}$
($d\geq0$). If $d=0,$ then $\omega$ is a $C^{3}$ connected, open
bounded subset of $\mathbb{R}^{N-1}.$ While, in the case where
$1\leq d\leq\, N-1,\;$ $\omega$ is a $(L_1,\ldots,L_d)$-periodic
open domain of $\,\mathbb{R}^{N-1}$ which satisfies (\ref{comega});
and hence, $\Omega$ is a $(l,L_1,\ldots,L_d)-$periodic subset of
$\mathbb{R}^{N}$ that satisfies (\ref{comega}) with $l>0$ and
arbitrary. An element of $\Omega=\mathbb{R}\times\omega$ will be
represented as $z=(x,y)$ where $x\in\mathbb{R}$ and
$y\in\omega\subseteq\mathbb{R}^{d}\times\mathbb{R}^{N-1-d}.$

The nonlinearity $f=f(x,y,u),$ in this section, is a KPP
nonlinearity defined on $\overline{\Omega}\times[0,1]$ that
satisfies
\begin{eqnarray}\label{nonlinearity 1}
     \left\{
      \begin{array}{ll}
        f\geq0, \;\hbox{and of class}\;C^{1,\,\delta}(\mathbb{R}\times\overline{\omega}\,\times\,[0,1]),  \\
        f \;\hbox{is $(l,L_1,\ldots,L_d)$-periodic with respect to}\;
        (x,y_1,\ldots,y_d),\;\hbox{when $d\geq1$},\vspace{3 pt}\\
f \;\hbox{is $l$-periodic in $x$},\;\hbox{when $d=0$},\vspace{3 pt}\\
        \forall\,(x,y)\,\in\,\overline{\Omega}=\mathbb{R}\times\overline{\omega},\; \displaystyle{f(x,y,0)=\,f(x,y,1)=\,0 } \hbox{,}\vspace{3 pt} \\
        \exists \,\rho\in\,(0,1),\;\forall (x,y)\,\in\overline{\Omega},\;\displaystyle{\forall\, 1-\rho\leq\,s\,\leq s'\,\leq\,1,\; f(x,y,s)\;\geq\,f(x,y,s') } \hbox{,} \vspace{3 pt}\\
        \forall\,s\in\,(0,1),\quad \exists \,(x,y)\,\in\,\overline{\Omega}\;\hbox{such that}\;f(x,y,s)\,>\,0  \hbox{,}
\end{array}
    \right.
\end{eqnarray}
together with the assumptions
\begin{eqnarray}\label{nonlinearity 2}
\left\{
  \begin{array}{ll}
     f'_{u}(x,y,0)\;\hbox{depends only on}\;y; \;\hbox{we
        denote
        by}\;\zeta(y)= f'_{u}(x,y,0),\;\forall(x,y)\in\overline{\Omega}.\vspace{3 pt}\\
        \forall\,(x,y)\,\in\,\overline{\Omega}=\mathbb{R}\times\overline{\omega},\quad \displaystyle{f'_{u}(x,y,0)\,=\,\zeta(y)\,>\,0}
        \hbox{,}\vspace{3 pt}\\
\forall\, (x,y,s)\,\in\,\overline{\Omega}\,\times\,(0,1),\quad
0\,<\,f(x,y,s)\,\leq\,\zeta(y)\,s \hbox{.}
  \end{array}
\right.
\end{eqnarray}
Notice that $f'_{u}(x,y,u)$ is assumed to depend only on $y,$ but
$f(x,y,u)$ may depend on $x.$

Lastly, concerning the diffusion matrix,
$A(x,y)=A(y)=(A_{ij}(y))_{1\leq i,j\leq
 N}$ is a $C^{2,\delta}(\,\overline{\Omega}\,)$ (with $\delta
\,>\,0$) symmetric matrix field whose entries are depending only on
$y,$ and satisfying
\begin{eqnarray}\label{cA Y}
    \left\{
      \begin{array}{l}
        A\; \hbox{is $(L_1,\ldots,L_d)$-periodic with respect to}\;(y_1,\ldots,y_d),\vspace{3 pt}\\
        \exists\,0<\alpha_1\leq\alpha_2,\;\forall\,y\in\omega,\forall\,\xi\,\in\,\mathbb{R}^N,\vspace{3 pt}\\
        \alpha_1|\xi|^2 \;\leq\;\sum \,A_{ij}(y)\xi_i\xi_j\,\;\leq\alpha_2|\xi|^2.
      \end{array}
    \right.
\end{eqnarray}

\begin{thm}\label{limit as eps}
 Let $e=(1,0,\ldots,0)\in\mathbb{R}^{N}$ and $\varepsilon>0.$ Let $\,\Omega=\mathbb{R}\times\omega\subseteq\,\mathbb{R}^N$ satisfy the form described in the previous page.
Under the assumptions (\ref{nonlinearity 1}), (\ref{nonlinearity
2}), and (\ref{cA Y}),
 consider the reaction-diffusion equation
\begin{eqnarray}\label{p eps}
    \left\{
      \begin{array}{rl}
        u_{t}(t,x,y)=&\varepsilon\, \nabla\cdot(A(y)\nabla
        u)(t,x,y)+\,f(x,y,u),\;\hbox{for}\;\;(t,x,y)\,\in\,\mathbb{R}\times\Omega\hbox{,}\vspace{3 pt}\\
   \nu\cdot A\nabla u=&0\quad\hbox{on}\quad\mathbb{R}\times\mathbb{R}\times
\partial\omega\hbox{.}
\end{array}
\right.
\end{eqnarray}
Assume, furthermore, that $A$ and $f$ satisfy one of the following
two alternatives:
\begin{eqnarray}\label{alt 1}
\left\{
  \begin{array}{ll}
  \exists\,\alpha>0,\;\forall y\in\omega,\;A(y)e=\alpha e,\vspace{3 pt}\\
  f^{\,'}_{u}(x,y,0)=\zeta(y),\,\hbox{for all} \;(x,y)\in\overline{\Omega,}
  \end{array}
\right.
\end{eqnarray}
or
\begin{eqnarray}\label{alt 2}
\left\{
  \begin{array}{ll}
  f^{\,'}_{u}(x,y,0)=\zeta\,\hbox{ is constant},\vspace{3 pt}\\
 \forall y\in\omega,\;A(y)e=\alpha(y)
 e,\;\hbox{where}\vspace{3 pt}\\
 \;y\mapsto\alpha(y)\;\hbox{is a positive, $(L_1,\ldots,L_d)-$periodic function
 over}\;\overline{\;\omega}.
  \end{array}
\right.
\end{eqnarray}
Then,
\begin{eqnarray}
\lim_{\varepsilon\rightarrow0^+}\,\frac{c_{\Omega,\varepsilon
A,0,f}^*(e)}{\sqrt{\varepsilon}}\,=\,2\sqrt{\max_{\overline{\omega}}\,\zeta}\sqrt{\max_{\overline{\omega}}\;eAe}.
\end{eqnarray}
\end{thm}

Before going further on, we mention that the family of domains for
which Theorem \ref{limit as eps} holds is wide. An infinite cylinder
$\displaystyle{\mathbb{R}\times B_{\mathbb{R}^{N-1}}(y_0,R)}$ (where
$R>0,$ and $\displaystyle{B_{\mathbb{R}^{N-1}}(y_0,R)}$ is the
Euclidian ball of center $y_0$ and radius $R$) is an archetype of
such domains. In these cylinders,
$\omega=\displaystyle{B_{\mathbb{R}^{N-1}}(y_0,R)},\;l\,$ is any
positive real number, and $d=0.$ The whole space $\mathbb{R}^{N}$ is
another archetype of the domain $\Omega$ where
$d=N-1,\;\omega=\mathbb{R}^{N-1},\;$ and $\;\{l,\,L_1,\ldots,L_d\}$
is any family of positive real numbers.
\begin{remark}
{\rm In Theorem \ref{limit as eps}, the domain
$\Omega=\mathbb{R}\times\omega$ is invariant in the direction of
$e=(1,0\ldots,0)$ which is parallel to $Ae$ ( in both cases
(\ref{alt 1}) and (\ref{alt 2})). Also, the assumption that the
entries of $A$ do not depend on $x,$ yields that
$\nabla.(Ae)\equiv0$ over $\Omega.$  On the other hand, it is easy
to find a diffusion matrix $A$ and a nonlinearity $f$ which satisfy,
together, the assumptions of Theorem \ref{limit as eps} while one of
$\;eAe(y)$ and $\zeta(y)$ is \textbf{not constant}. Referring to
part c) of Theorem \ref{influence of f, A}, one obtains:
$$\displaystyle{\forall\varepsilon>0,\quad0<\frac{c_{\Omega,\varepsilon
A,0,f}^*(e)}{\sqrt{\varepsilon}}\,\lneqq\,2\sqrt{\max_{y\,\in\,\overline{\omega}}\,\zeta(y)}\sqrt{\max_{y\,\in\,\overline{\omega}}\,eAe(y)}}.$$
However, Theorem \ref{limit as eps} implies that
$$\displaystyle{\lim_{\varepsilon\rightarrow0^+}\,\frac{c_{\Omega,\varepsilon
A,0,f}^*(e)}{\sqrt{\varepsilon}}\,=\,2\sqrt{\max_{y\,\in\,\overline{\omega}}\,\zeta(y)}\sqrt{\max_{y\,\in\,\overline{\omega}}\,eAe(y)}}.$$

On the other hand, if $\;\Omega=\mathbb{R}\times\omega$ as in
Theorem \ref{limit as eps}, $A=Id$ and $f=f(u),$ Theorem
\ref{influence of f, A} yields that $\displaystyle{\,
c_{_{\Omega,\varepsilon
Id,0,f}}^*(e)\,=\,2\sqrt{\varepsilon}\sqrt{f'(0)},} $ for all
$\varepsilon>0.$ $\hfill{\Box}$}\end{remark}

In the same context, one can also find the limit when the diffusion factor goes to zero, but in the presence of an advection field in the form of shear flows:
\begin{thm}\label{lim as eps but in presence of a shear flow}
 Assume that $e=(1,0,\cdots,0)\in\mathbb{R}^{N},$ the domain $\Omega=\mathbb{R}\times\omega$ has the same form as in Theorem \ref{limit as eps},
  and the coefficients $f$ and $A$ satisfy (\ref{nonlinearity 1}-\ref{nonlinearity 2}) and (\ref{cA Y}) respectively.
Assume, furthermore, that for all $y\in\overline{\omega},$ there exists $\alpha(y)$ positive so that $A(y)e=\alpha(y)e$ in $\overline{\omega}.$
Consider, in addition, an advective shear flow $q=(q_{_{1}}(y),0,\ldots,0)$
($y\in\overline{\omega}$) which is $(L_1,\cdots,L_d)-$periodic with respect to $y.$ Assume that
 $\varepsilon$ is a positive parameter and consider the parametric
reaction-advection-diffusion problem
\begin{eqnarray}\label{eps A with q}
\left\{
  \begin{array}{ll}
    u_t =\,\varepsilon\,\nabla\cdot(A(y)\nabla u)+\;q_{_{1}}(y)\,\partial_{x}u(t,x,y)\;+\,f(x,y,u),\quad t\,\in\,\mathbb{R},\;(x,y)\,\in\,\Omega,
    \vspace{3 pt}\\
        \nu \cdot A\;\nabla u(t,x,y) =0,\quad
        t\,\in\,\mathbb{R},\quad(x,y)\,\in\,\partial\Omega,
  \end{array}
\right.
\end{eqnarray}
where $q\not\equiv0$ over $\mathbb{R}\times\overline{\omega}$ and
$q$ has a zero average. Then,
\begin{equation}\label{lim as eps tend to 0 with q}
    \displaystyle{\lim_{\varepsilon\rightarrow0^{+}}c_{\Omega,\varepsilon
A,q,f}^*(e)=\max_{y\in\overline{\omega}}\,\left(-\,q_{_{1}}(y)\right)=\,\max_{\overline{\omega}}(-\,q.e)}.
\end{equation}
\end{thm}

The situation in this result is more general than that
considered in part b) of Corollary~4.5 in \cite{BHN2}. In details,
the coefficients $A$ and $f$ can be both non-constant. Meanwhile, in the result of
\cite{BHN2}, the coefficients considered were assumed to satisfy the alternative (\ref{alt 1}).

After having the exact value of
$\displaystyle{\lim_{\varepsilon\rightarrow0^+}\,\frac{c_{\Omega,\varepsilon
A,0,f}^*(e)}{\sqrt{\varepsilon}}},$ we move now to investigate the
limit of the minimal speed of propagation, considered as a function
of the period of the coefficients of the reaction-diffusion equation
set in the whole space $\mathbb{R}^{N},$
 when the periodicity parameter tends to $+\infty.$ By making some change in variables, we will find a link between this problem and Theorem \ref{limit as eps}:
\begin{thm}\label{lim as per L tend to infty}
Let $e=(1,0,\ldots,0)\in\mathbb{R}^{N}.$ An element
$z\in\mathbb{R}^N$ is represented as
$z=(x,y)\in\mathbb{R}\times\mathbb{R}^{N-1}.$ Assume that $f=f(x,y,u)$
and $A=A(y)$ satisfy (\ref{nonlinearity 1}), (\ref{nonlinearity 2})
and (\ref{cA Y}) with $\omega=\mathbb{R}^{N-1},\, d=N-1,\,$ and
$\;l=L_1=\ldots=L_{N-1}=1.$ (That is, the domain and the
coefficients of the equation are $(1,1,\ldots,1)$ periodic with
respect to $y$). Assume furthermore, that $A$ and $f$ satisfy either
(\ref{alt 1}) or (\ref{alt 2}). For each $L>0,\;$ and
$(x,y)\in\mathbb{R}^{N},$ let
 $\displaystyle{A_{_{L}}(y)=A(\frac{y}{L})}$ and
$\displaystyle{f_{_{L}}(x,y,u)=f(\frac{x}{L},\frac{y}{L},u)}.$
Consider the reaction-diffusion problem
\begin{eqnarray}\label{KPP L with u}
 \begin{array}{ll}
 u_t(t,x,y)&= \displaystyle{\nabla\cdot(A_{_{L}}\nabla u)(t,x,y)\,+\,f_{_{L}}(x,y,u),
    \,(t,x,y)\in\mathbb{R}\times\mathbb{R}^{N}}\vspace{4 pt}\\
 &=\displaystyle{\nabla\cdot(A(\frac{y}{L})\nabla
    u)(t,x,y)\,+\,f(\frac{x}{L},\frac{y}{L},u),
    \,(t,x,y)\in\mathbb{R}\times\mathbb{R}^{N}},
\end{array}
\end{eqnarray}
whose  coefficients are $(L,\ldots,L)$ periodic with respect to
$(x,y)\in\mathbb{R}^{N}.$ Then,
\begin{eqnarray}
\displaystyle{ \lim_{L\rightarrow\,+\infty}\,c_{\mathbb{R}^{N},
\displaystyle{A_{_{L}},0,f_{_{L}}}}^*(e)=2\sqrt{\max_{y\,\in\,\mathbb{R}^{N-1}}\,\zeta(y)}\sqrt{\max_{y\,\in\,\mathbb{R}^{N-1}}\,e.Ae(y)}}.
\end{eqnarray}
\end{thm}

The above theorem gives the limit of the minimal speed of
propagation in the direction of $e=(1,0,\cdots,0)$ as the
periodicity parameter $L\rightarrow+\infty.$ The domain is the whole
space $\mathbb{R}^{N}$ which is $(L,\cdots,L)-$periodic whatever the
positive number $L$. However, one can find $$\displaystyle{ \lim_{L\rightarrow\,+\infty}\,c_{\mathbb{R}^{N},
\displaystyle{A_{_{L}},Lq_{_L},f_{_{L}}}}^*(e)}$$ whenever $q$ is a shear flow advection. Namely, in the same manner that Theorem \ref{limit as eps} implies Theorem \ref{lim as per L tend to infty}, one can prove that Theorem \ref{lim as eps but in presence of a shear flow} implies
\begin{thm}\label{lim as L tends to infty thm in presence of q}
Let $e=(1,0,\ldots,0)\in\mathbb{R}^{N}.$ Assume that $f=f(x,y,u)$
and $A=A(y)$ satisfy (\ref{nonlinearity 1}), (\ref{nonlinearity 2})
and (\ref{cA Y}) with $\omega=\mathbb{R}^{N-1},\, d=N-1,\,$ and
$\;l=L_1=\ldots=L_{N-1}=1.$ (That is, the domain and the
coefficients of the equation are $(1,1,\ldots,1)$ periodic with
respect to $y$ in $\R^{N-1}$). Assume, furthermore, that for all $y\in\mathbb{R}^{N-1},$ there exists $\alpha(y)$ positive so that $A(y)e=\alpha(y)e$ in $\mathbb{R}^{N-1}.$ Let $q=(q_{_{1}}(y),0,\ldots,0)$
for all $y\in\mathbb{R}^{N-1}$ such that  $q_1\not\equiv0$ over $\mathbb{R}^{N-1},$ $q$ is $(1,\cdots,1)-$periodic with respect to $y$ and
$q_1$ has a zero average. Then,
\begin{equation}\label{lim as L tends to infty but with q}
\displaystyle{\lim_{L\rightarrow\,+\infty}\,c_{\mathbb{R}^{N},
\displaystyle{A_{_{L}},Lq_{_L},f_{_{L}}}}^*(e)=\max_{y\in\mathbb{R}^{N-1}}\,\left(-\,q_{_{1}}(y)\right)=\,\max_{\mathbb{R}^{N-1}}(-\,q.e)}.
\end{equation}
\end{thm}

 In the proof of Theorem \ref{lim as eps but in presence of a shear flow} (which implies Theorem \ref{lim as L tends to infty thm in presence of q}), the
 assumption that the advection $q$ is
in the form of shear flows plays an important role in reducing the elliptic equation involved by the variational formula
(\ref{var form with eps and q}) below. Namely, since $q=(q_1(y),0,\cdots,0)$ and since $e=(1,0,\cdots,0),$ then the terms $q(x,y)\cdot\nabla_{x,y}\psi$ and
$q(x,y)\cdot e$ (in the general elliptic equation)
become equal to $q_1(y)\partial_x\psi$ and $q_1(y)$ respectively. As a consequence, and due the uniqueness of the principal eigenfunction $\psi$ up to
multiplication by a constant, we are able to choose $\psi$ independent of $x,$ and hence, obtain a symmetric elliptic operator (without drift)
whose principal eigenvalue was given by the variational formula (\ref{rayleigh}) below (see section \ref{proofs section} for more details).
\begin{remark}\label{without shear}
After the above explanations, we find that the techniques used to prove Theorem \ref{lim as eps but in presence of a shear flow} which implies \ref{lim as L tends to infty thm in presence of q},
 will  no longer work in the presence a general periodic advection field satisfying (\ref{cq}).
\end{remark}

 Concerning the influence of advection, we mention that the limit of $\ds{\frac{c^{*}_{\Omega,A,\,Bq,f}(e)}{B}}$ as $B\rightarrow+\infty$ (in the general periodic setting) is not yet given explicitly as a function of the direction $e$ and the coefficients $A,$ $q$ and $f.$ For more details one can see Theorem 4.1 in \cite{BHN2}. However, the problem of front
 propagation in an infinite cylinder with an underlying shear flow
 was widely studied in Berestycki \cite{6}, Berestycki and Nirenberg
 \cite{BNir2}. In the case of
 strong advection, assume that $\Omega=\mathbb{R}\times\omega,$ where $\omega$ is a bounded smooth
 subset of $\mathbb{R}^{N-1},\;q=(q_{_{1}}(y),0,\cdots,0),\;y\in\omega,\;$ and $f=f(u)$ is a (KPP) nonlinearity. It was proved, in Heinze \cite{Heinze convection}, that
 \begin{equation}\label{Heinze limit of c(M)/M}
\displaystyle{\lim_{B\rightarrow+\infty}\frac{c_{\Omega,A,\,Bq,f}^*(e)}{B}=\gamma,}
 \end{equation}
 where $$\displaystyle{\gamma=\sup_{\psi\in\,D}\int_{\omega}\,q_{_{1}}(y)\,\psi^{2}\,dy},$$
 $$D=\left\{\psi\in H^{1}(w), \quad \int_{\omega}|\nabla\psi|^{2}\,dy\,\leq\,f^{'}(0),\;\hbox{and}\;\int_{\omega}\,\psi^{2}\,dy=1\right\}.$$

\section{The minimal speed
within large diffusion factors or within small period
coefficients}\label{as M goes to infty} After having the limit of
$\displaystyle{c_{\Omega,\varepsilon
A,0,f}^*(e)/{\sqrt{\varepsilon}}}$ as $\varepsilon\rightarrow0^+,$
and after knowing that this limit depends on
$\displaystyle{\max_{y\,\in\,\overline{w}}\,\zeta(y)}\;$ and
$\;\displaystyle{\max_{y\,\in\,\overline{w}}\,eAe(y)},$ we
investigate now the limit of
$\;\displaystyle{{c_{\Omega,MA,\,q,f}^*(e)}/{\sqrt{M}}}\,$ as
the diffusion factor $M$ tends to $+\infty,$ and we try to answer
this question in a situation which is more general than that we
considered in the previous section (in the case where the diffusion
factor was going to $0^{+}$). That is in the presence of an
advection field and in a domain $\Omega$ which satisfies
(\ref{comega}) and which may take more forms other than those of
section \ref{as eps tends to zero}. We will find that in the case of
large diffusion, the limit will depend on
$\displaystyle{\mi_{\!\!\!\!_C}\zeta(x,y)dx\,dy:=\frac{1}{|C|}\,\int_{C}\zeta(x,y)}dx\,dy$
$~~\hbox{and}\quad
\displaystyle{\mi_{\!\!\!\!_C}\tilde{e}A\tilde{e}(x,y)dx\,dy:=\frac{1}{|C|}\int_{C}\tilde{e}A\tilde{e}(x,y)dx\,dy},$
where $C$ denotes the periodicity cell of the domain $\Omega.$
\begin{thm}\label{lim as M}
Under the assumptions (\ref{comega}) for $\Omega,$ (\ref{cq}) for
the advection $q,\;$(\ref{cf1}) and (\ref{cf2}) for the nonlinearity
$f=f(x,y,u),$ let $e$ be any unit direction of $\,\mathbb{R}^{d}.$
Assume that the diffusion matrix $A=A(x,y)$ satisfies (\ref{cA})
together with $\nabla\cdot A\tilde{e}\equiv0\;$ over $\Omega,$ and
$\nu\cdot A\tilde{e}=0$ over $\partial\Omega.$ For each $M>0$ and
$\displaystyle{0\leq\gamma\leq\,{1}/{2}},$ consider the following
reaction-advection-diffusion equation
\begin{equation*}
    \left\{
      \begin{array}{l}
 \displaystyle{u_t =M\,\nabla\cdot(A(x,y)\nabla u)\;+\,\displaystyle{M^{\,\gamma}\,q(x,y)\cdot\nabla u}\,+\,f(x,y,u),\; t\,\in\,\mathbb{R},\;(x,y)\,\in\,\Omega,} \\
        \nu \cdot A\;\nabla u(t,x,y) =0,\;
        t\,\in\,\mathbb{R},\;(x,y)\,\in\,\partial\Omega.
      \end{array}
    \right.
\end{equation*}
Then
$$\displaystyle{\lim_{M\rightarrow+\infty}\frac{c_{\Omega,MA,\displaystyle{M^{\gamma}\,q},f}^*(e)}{\sqrt{M}}=2\sqrt{\mi_{\!\!\!\!_C}\tilde{e}A\tilde{e}(x,y)dx\,dy}\,\sqrt{\mi_{\!\!\!\!_C}\zeta(x,y)dx\,dy}},$$
where $C$ is the periodicity cell of $\Omega.$
\end{thm}
\begin{remarks}
\end{remarks}
\begin{itemize}
  \item The setting in Theorem \ref{lim as M} is more general than that in
Theorem \ref{limit as eps}, where:
$\Omega=\mathbb{R}\times\omega,\,\tilde{e}=(1,0,\ldots,0),$ and
$A\tilde{e}=\alpha(y)\tilde{e}$. Under the assumptions of Theorem
\ref{limit as eps}, the domain $\;\Omega\;$ is invariant in the
direction of $A\tilde{e},$ which is that of $\tilde{e}.$
Consequently, if $\nu$ denotes the outward normal on
$\partial\Omega=\mathbb{R}\times\partial\omega,$ one gets $\nu\cdot
A\tilde{e}=\alpha(y)\,\nu\cdot\tilde{e}=0$ over $\partial\Omega,$
while
$\displaystyle{\nabla\cdot(A\tilde{e})=\frac{\partial}{\partial
x}\alpha(y)=0}$ over $\Omega.$ Moreover, in Theorem \ref{limit as
eps}, we have only reaction and diffusion terms. That is $q\equiv0.$
Therefore, considering the setting of Theorem \ref{limit as eps},
and taking $\beta A$ as a parametric diffusion matrix,  one
consequently knows the limits of
$\;\displaystyle{\frac{c_{_{\Omega,\beta
A,0,f}}^*(e)}{\sqrt{\beta}}}\;$ as $\beta\rightarrow0^+$ (Theorem
\ref{limit as eps}) and as $\beta\rightarrow+\infty$ (Theorem
\ref{lim as M}).
  \item The other observation in Theorem \ref{lim as M} is
that the limit does not depend on the advection field $q.$ This
 may play an important role in drawing counter examples
to answer many different questions. For example, the variation
of the minimal speed of propagation with respect to the diffusion
factor and with respect to diffusion matrices which are symmetric
positive definite.
\item Another important feature, in Theorem \ref{lim as M}, is that the order of $M$ in the denominator of the ratio $\displaystyle{\frac{c_{\Omega,MA,\displaystyle{M^{\gamma}\,q},f}^*(e)}{\sqrt{M}}}$ is equal to  $\;\displaystyle{{1}/{2}}.\;$
It is independent of $\gamma.$ Consequently, \textbf{the case where
the advection is null} and there is only a reaction-diffusion
equation follows, in particular, from the previous theorem. That is
$$\displaystyle{\lim_{M\rightarrow+\infty}\frac{c_{_{\Omega,MA,0,f}}^*(e)}{\sqrt{M}}=2\sqrt{\mi_{\!\!\!\!_C}\tilde{e}A\tilde{e}(x,y)dx\,dy}\,\sqrt{\mi_{\!\!\!\!_C}\zeta(x,y)dx\,dy}}.$$
\item The previous point leads us to conclude that the presence of
an advection with a factor $M^{\gamma},$ where
$\displaystyle{0\leq\gamma\leq{1}/{2}},$ will have no more effect on
the ratio
$\displaystyle{\frac{c_{\Omega,MA,\displaystyle{M^{\gamma}\,q},f}^*(e)}{\sqrt{M}}}$
as soon as the diffusion factor $M$ gets very large.
\end{itemize}

As far as the limit of the minimal speed of propagation within small
periodic coefficients in the reaction-diffusion equation is
concerned, the following theorem, which mainly depends on
Theorem~\ref{lim as M}, treats this problem:
\begin{thm}\label{lim as per L tend to 0}
Let $\Omega=\mathbb{R}^{N}.$ Assume that $A=A(x,y),\;q=q(x,y)$ and
$f=f(x,y,u)$ are $(1,\ldots,1)-$periodic with respect to
$(x,y)\in\mathbb{R}^{N},$ and that they satisfy (\ref{cA}),
(\ref{cq}), (\ref{cf1}),and (\ref{cf2}) with $\;L_1=\ldots=L_N=1.$
Let $e$ be any unit direction of $\mathbb{R}^{N},$ such that
$\;\nabla\cdot A\tilde{e}\equiv0$ over $\,\mathbb{R}^{N}.$ For each
$L>0,\;$  let
$\displaystyle{A_{_{L}}(x,y)=A(\frac{x}{L},\frac{y}{L}),}$ $\ds{q_{_{L}}(x,y)=q(\frac{x}{L},\frac{y}{L}),}$
and $\displaystyle{f_{_{L}}(x,y,u)=f(\frac{x}{L},\frac{y}{L},u)},$
where $(x,y)\in\mathbb{R}^{N}.\;$ Consider the problem
\begin{eqnarray}\label{KPP L with u 2}
   \begin{array}{ll}
u_t(t,x,y)&= \displaystyle{\nabla\cdot(A_{_{L}}\nabla
    u)(t,x,y)\,+\,q_{_{L}}\cdot\nabla u(t,x,y)+f_{_{L}}(x,y,u),
    \,(t,x,y)\in\mathbb{R}\times\mathbb{R}^{N}},\\\\
 &= \displaystyle{\nabla\cdot(A(\frac{x}{L},\frac{y}{L})\nabla
    u)(t,x,y)\,+\,q(\frac{x}{L},\frac{y}{L})\cdot\nabla u(t,x,y)+f(\frac{x}{L},\frac{y}{L},u)},
\end{array}
\end{eqnarray}
whose  coefficients are $(L,\ldots,L)$ periodic with respect to
$(x,y)\in\mathbb{R}^{N}.$ Then,
\begin{eqnarray}
\nonumber\displaystyle{
\lim_{L\rightarrow\,0^+}\,c_{{\mathbb{R}^{N},\displaystyle{A_{_{L}},\,q_{_{L}},f_{_{L}}}}}^{*}(e)=2\sqrt{\mi_{\!\!\!\!_C}\tilde{e}A\tilde{e}(x,y)dx\,dy}\,\sqrt{\mi_{\!\!\!\!_C}\zeta(x,y)dx\,dy}},
\end{eqnarray}
where, in this setting, $\displaystyle{C=
[0,1]\times\cdots\times[0,1]\subset\mathbb{R}^{N}}.$
\end{thm}

The above result gives the limit in any space dimension. It depends
on the assumption $\nabla\cdot(A\tilde{e})\equiv0$ in
$\mathbb{R}^{N}.$ However, if one takes $N=1,$ and denotes the
diffusion coefficient by $a=a(x),\;x\in\mathbb{R},$ then the
previous result holds under the assumptions that $a$ satisfies
(\ref{cA}) and $\displaystyle{{da}/{dx}\equiv0}$ in $\mathbb{R}.$ In
other words, it holds when $a$ is a positive constant. Thus, it is
be interesting to mention that, in the one-dimensional case, the
above limit was given in \cite {EHR} and \cite{KKS} within a general diffusion
coefficient (which may be not constant over $\mathbb{R}$). In
details, assume that $f=f(x,u)= (\zeta(x)-u)u$ is a $1$-periodic
(KPP) nonlinearity satisfying (\ref{cf1}) with (\ref{cf2}), and
$\mathbb{R}\ni x\mapsto a(x)\;$ is a $1-$periodic function which
satisfies $0<\alpha_{1}\leq a(x)\leq\alpha_{2},$ for all
$x\in\mathbb{R},$ where $\alpha_1$ and $\alpha_{2}$ are two positive
constants. For each $L>0,\;$ consider the reaction-diffusion
equation
\begin{equation}\label{shigesada}
\partial_{t}\,u(t,x)=\frac{\partial}{\partial\,x}\,\left(a(\frac{x}{L})\frac{\partial\,u}{\partial
x}\right)(t,x)\;+\;\left[\zeta(\frac{x}{L})-u(t,x)\right]u(t,x)\quad\hbox{for}\quad(t,x)\in\mathbb{R}\times\mathbb{R}.
\end{equation}

It was derived in \cite{EHR} and, formally, in \cite{KKS} that
\begin{equation}\label{limit of shigesada}
\displaystyle{
\lim_{L\rightarrow\,0^+}\,c_{{\mathbb{R},\displaystyle{a_{_{L}},\,0,f_{_{L}}}}}^{*}(e)=2\,\sqrt{<a>_{_{H}}.\int_{0}^{1}\zeta(x)}},
\end{equation}

where $\displaystyle{<a>_{_{H}}}$ denotes the harmonic mean of the
map $x\mapsto a(x)$ over $[0,1].$

\section{The minimal speed within small or large reaction
coefficients}\label{reaction} In this section, the parameter of the
reaction-advection-diffusion problem is the coefficient $B$
multiplied by the nonlinearity $f.$ In fact, it follows from Theorem
1.6 in Berestycki, Hamel and Nadirashvili \cite{BHN1} (recalled via
Theorem \ref{influence of f, A} in the present paper) that the map
$\displaystyle{B\mapsto{c_{\Omega,A,q,Bf}^*(e)}/{\sqrt{B}}}$
remains, with the assumption $\nu.A\tilde{e}=0$ on $\partial\Omega,$
 bounded by two positive constants as $B$ gets very large. Therefore, it is interesting to
 find the limit of
 $\displaystyle{{c_{\Omega,A,q,Bf}^*(e)}/{\sqrt{B}}}$ as
 $B\rightarrow+\infty$ even in some particular situations. Moreover, it
 is important to find the limit of the same quantity as $B\rightarrow0^{+}.$
We start with the case where $B\rightarrow+\infty$ and then we move
to that where $B\rightarrow 0^{+}.$
\begin{thm}\label{lim as B tend to infty}
 Let $e=(1,0,\ldots,0)\in\mathbb{R}^{N}$ and $B>0.$ Assume that $\,\Omega=\mathbb{R}\times\omega\subseteq\,\mathbb{R}^N,\;A,\;\hbox{and}\;f$  satisfy the same assumptions of Theorem \ref{limit as eps}.
That is, $f$ and $A$ satisfy (\ref{nonlinearity 1}),
(\ref{nonlinearity 2}), and (\ref{cA Y}), and one of the two
alternatives (\ref{alt 1})-(\ref{alt 2}). Consider the
reaction-diffusion equation
\begin{eqnarray}\label{p B}
    \left\{
      \begin{array}{rl}
        u_{t}(t,x,y)=& \nabla\cdot(A(y)\nabla
        u)(t,x,y)+\,B\,f(x,y,u),\;\hbox{for}\;\;(t,x,y)\,\in\,\mathbb{R}\times\Omega\hbox{,}\\
   \nu\cdot A\nabla u=&0\quad\hbox{on}\quad\mathbb{R}\times\mathbb{R}\times
\partial\omega\hbox{.}
\end{array}
\right.
\end{eqnarray}

Then,
\begin{eqnarray}
\lim_{B\rightarrow+\infty}\,\frac{c_{\Omega, A,0,B
f}^*(e)}{\sqrt{B}}\,=\,2\sqrt{\max_{y\,\in\,\overline{\omega}}\,\zeta(y)}\sqrt{\max_{y\,\in\,\overline{\omega}}\;eAe(y)}.
\end{eqnarray}
\end{thm}

We mention that one can find the coefficients $A,$ and $f$ and the
domain $\Omega$ of the problem (\ref{p B}) satisfying all the
assumptions of Theorem \ref{lim as B tend to infty}, which are the
same of Theorem \ref{limit as eps}, including one of the
alternatives (\ref{alt 1})-(\ref{alt 2}) while one of $\zeta$ and
$eAe$ is \textbf{not constant}. Owing to Theorem 1.10 in
\cite{BHN1}, it follows that $$\displaystyle{\forall B>0,\quad
c_{\Omega, A,0,B
f}^*(e)\lneqq\,2\sqrt{B}\sqrt{\max_{y\in\overline{\omega}}\,\zeta(y)}\,\sqrt{\max_{y\in\overline{\omega}}\,eAe(y)}},$$
which is equivalent to saying that
$$\displaystyle{\frac{c_{\Omega,
A,0,Bf}^*(e)}{\sqrt{B}}\,\lneqq\,2\sqrt{\max_{y\,\in\,\overline{\omega}}\,\zeta(y)}\sqrt{\max_{y\,\in\,\overline{\omega}}\;eAe(y)}}.$$

Therefore, there are \textsl{heterogeneous} settings in which the
result found in Theorem \ref{lim as B tend to infty} does not follow
trivially.

We move now to study the limit when the reaction factor $B$ tends to
$0^+.$ However, the situation will be more general than that in
Theorem \ref{lim as B tend to infty} because it will consider
reaction-advection-diffusion equations rather than considering
reaction-diffusion equations only:
\begin{thm}\label{lim as B tend 0}
Under the assumptions (\ref{comega}) for $\Omega,$ (\ref{cq}) for
the advection $q,\;$(\ref{cf1}) and (\ref{cf2}) for the nonlinearity
$f=f(x,y,u),$ let $e$ be any unit direction of $\,\mathbb{R}^{d}.$
Assume that the diffusion matrix $A=A(x,y)$ satisfies (\ref{cA})
together with $\nabla\cdot A\tilde{e}\equiv0\;$ over $\Omega,$ and
$\nu\cdot A\tilde{e}=0$ over $\partial\Omega.$ For each $B>0$ and
$\displaystyle{\gamma\geq\,{1}/{2}},$ consider the following
reaction-advection-diffusion equation
\begin{equation*}
    \left\{
      \begin{array}{l}
 \displaystyle{u_t =\nabla\cdot(A(x,y)\nabla u)\;+\,\displaystyle{B^{\,\gamma}\,q(x,y)\cdot\nabla u}\,+\,B\,f(x,y,u),\quad t\,\in\,\mathbb{R},\;(x,y)\,\in\,\Omega,} \\
        \nu \cdot A\;\nabla u(t,x,y) =0,\quad
        t\,\in\,\mathbb{R},\quad(x,y)\,\in\,\partial\Omega.
      \end{array}
    \right.
\end{equation*}
Then
$$\displaystyle{\lim_{B\rightarrow0^{+}}\frac{c_{\Omega,A,\displaystyle{B^{\gamma}\,q},Bf}^*(e)}{\sqrt{B}}=2\sqrt{\mi_{\!\!\!\!_C}\tilde{e}A\tilde{e}(x,y)dx\,dy}\,\sqrt{\mi_{\!\!\!\!_C}\zeta(x,y)dx\,dy}},$$
where $C$ is the periodicity cell of $\Omega.$
\end{thm}

Having the above result one can mark a sample of notes:

The order of $B$ in the denominator of the ratio
$\;\displaystyle{c_{\Omega,A,\displaystyle{B^{\gamma}\,q},Bf}^*(e)/\sqrt{B}}\;$
is independent of $\gamma$ (it is equal to ${1}/{2}$). Thus,
whenever the advection is null, one gets
$$\displaystyle{\lim_{B\rightarrow0^{+}}\frac{c_{\Omega,A, 0,Bf}^*(e)}{\sqrt{B}}=2\sqrt{\mi_{\!\!\!\!_C}\tilde{e}A\tilde{e}(x,y)dx\,dy}\,\sqrt{\mi_{\!\!\!\!_C}\zeta(x,y)dx\,dy}}.$$

Therefore, one concludes that the presence of an advection with a
factor $B^{\gamma},$ where $\displaystyle{\gamma\geq\;{1}/{2}},$
will have no more effect on the limit of the ratio
$\displaystyle{c_{\Omega,A, \displaystyle{B^{\gamma}
q},Bf}^*(e)/\sqrt{B}}$ as the reaction factor $B$ gets very small.

On the other hand, it is easy to check that the assumptions in
Theorem \ref{lim as B tend 0} are more general than those in Theorem
\ref{lim as B tend to infty}. Consequently, once we are in the more
strict setting, which is that of Theorem \ref{lim as B tend to
infty}, we are able to know both limits of
$\displaystyle{{c_{\Omega,A,0,Bf}^*(e)}/{\sqrt{B}}}$ as
$B\rightarrow +\infty$ and as $B\rightarrow0^{+}.$

\section{Variations of the minimal speed with respect to diffusion
and reaction factors and with respect to periodicity
parameters}\label{variation of min speed w.r.t diffusion and period}
After having studied the limits and the asymptotic behaviors of the
of the functions
$\;\displaystyle{\varepsilon\mapsto\,c_{\Omega,\varepsilon
A,0,f}^*(e)/{\sqrt{\varepsilon}}},$ $\displaystyle{M\mapsto\,c_{\Omega,MA,\displaystyle{M^{\gamma}\,q},f}^*(e)/{\sqrt{M}}}$
(for very large $M$  and for
$0\leq\gamma\leq{1}/{2}),$
$\displaystyle{B\mapsto\,c_{\Omega,A,\displaystyle{B^{\gamma}\,q},Bf}^*(e)/{\sqrt{B}}}$ ($\gamma\geq\displaystyle{{1}/{2}}$) and
$L \mapsto c_{\mathbb{R}^{N},A_L,\,q_{_{L}},f_{_{L}}}^{*}(e),$ where
$L$ is a periodicity parameter, we move now to investigate the
variations of these functions with respect to the diffusion and
reaction factors and with respect the periodicity parameter $L.$ The
present section will be devoted to discuss and answer these
questions.

We sketch first the form of the domain.
 $\Omega\subseteq \mathbb{R}^{N}$ is assumed to be in the form $\mathbb{R}\times\omega$ which was taken
in section \ref{as eps tends to zero}. As a review,
$\,\Omega=\mathbb{R}\times\omega\subseteq\,\mathbb{R}^N,$ where
$\,\omega\subseteq\mathbb{R}^{d}\times\mathbb{R}^{N-d-1}$
($d\geq0$). If $d=0,$ the subset $\omega$ is a bounded open subset
of $\mathbb{R}^{N-1}.$ While, in the case where $1\leq d\leq\,
N-1,\;$ $\omega$ is a $(L_1,\ldots,L_d)$-periodic open domain of
$\,\mathbb{R}^{N-1}$ which satisfies (\ref{comega}); and hence,
$\Omega$ is a $(l,L_1,\ldots,L_d)-$ periodic subset of
$\mathbb{R}^{N}$ that satisfies (\ref{comega}) with $l>0.$ An
element of $\Omega=\mathbb{R}\times\omega$ will be represented as
$z=(x,y)$ where
$y\in\omega\subseteq\mathbb{R}^{d}\times\mathbb{R}^{N-1-d}.$ With a
domain of such form, we have:

\begin{thm}\label{variation of min speed with respect diffusion}
Let $\;e = (1,0,\ldots,0)\in\mathbb{R}^{N}.$ Assume that $\Omega$ has the form
$\mathbb{R}\times\omega$ which is described above, and that the
diffusion matrix $A = A(y)$ satisfies (\ref{cA Y}) together with the
assumption
\begin{equation}\label{Ae= alpha(y)e}
    A(x,y)e=A(y)e=\alpha(y)e,\;\hbox{for all}\; (x,y)\in\mathbb{R}\times\overline{\omega};
\end{equation}
where $y\mapsto\alpha(y)$ is a positive $(L_1,\ldots,L_d)-$ periodic
function defined over $\overline{\omega}.$ The nonlinearity $f$ is
assumed to satisfy (\ref{nonlinearity 1}) and (\ref{nonlinearity
2}). Moreover, one assumes that, at least, one of $\tilde {e}\cdot A\tilde {e}$ and $\zeta$ is not constant. Besides, the advection field $q$ (when it
exists) is in the form $q(x,y)=(q_{_{1}}(y),0,\ldots,0)$ where
$q_{_{1}}$ has a zero average over $C,$ the periodicity cell of
$\omega.$ For each $\beta>0$ consider the
reaction-advection-diffusion problem
\begin{equation*}
    \left\{
      \begin{array}{l}
 u_t =\beta\,\nabla\cdot(A(y)\nabla u)\;+\,\sqrt{\beta}\,q_{_{1}}(y)\,\partial_{x} u\,+\,f(x,y,u),\; t\,\in\,\mathbb{R},\;(x,y)\,\in\,\mathbb{R}\times\omega,\\
        \nu \cdot A\;\nabla u(t,x,y) =0,\quad
        t\,\in\,\mathbb{R},\;(x,y)\,\in\,\partial\Omega.
      \end{array}
    \right.
\end{equation*}

 Then the map $\displaystyle{\,\beta\mapsto\frac{c_{\Omega,\beta A,\sqrt{\beta}\,q,f}^*(e)}{\sqrt{\beta}}\,}$ is decreasing in
$\beta>0,$ and by Theorem \ref{lim as M}, one has
$$\displaystyle{\lim_{\beta\rightarrow+\infty}\frac{c_{\Omega,\beta A,\sqrt{\beta}\,q,f}^*(e)}{\sqrt{\beta}}=2\sqrt{\mi_{\!\!\!\!_C}\tilde{e}A\tilde{e}(y)dy}\,\sqrt{\mi_{\!\!\!\!_C}\zeta(y)dy}},$$
where $C$ is the periodicity cell of $\omega.$
\end{thm}

\begin{remark}
In the same setting of Theorem \ref{variation of min speed with
respect diffusion} but with no advection, that is $q_{_{1}}\equiv0,$
we still have $\displaystyle{\,\beta\mapsto\frac{c_{\Omega,\beta
A,0,f}^*(e)}{\sqrt{\beta}}\,}$ as a decreasing map in $\beta>0.$
Moreover, if one of the alternatives (\ref{alt 1})-(\ref{alt 2})
holds and there is no advection, Theorem \ref{limit as eps} yields
that
$$\displaystyle{\lim_{\beta\rightarrow0^+}\frac{c_{\Omega,\beta
A,0,f}^*(e)}{\sqrt{\beta}}=2\sqrt{\max_{y\in\overline{\omega}}\tilde{e}A\tilde{e}(y)}\,\sqrt{\max_{y\in\overline{\omega}}\zeta(y)}}.$$
\end{remark}

 The preceding result yields another one concerned in the variation of the minimal speeds with respect to the
periodicity parameter $L.$ In the following, the domain will be the
whole space $\mathbb{R}^{N}.$ We choose the diffusion matrix
$A(x,y)=A(y),$ the shear flow $q$ and reaction term $f$ to be
$(1,\ldots,1)$-periodic and to satisfy some restrictions. For each
$L
> 0,$ we assign the diffusion matrix $A_L(x,y) = A(
\displaystyle{\frac{x}{L},\frac{y}{L}}),$ the advection field
$q_{_{L}}(x,y)=q(\displaystyle{\frac{x}{L},\frac{y}{L}})$ and the
nonlinearity $\displaystyle{f_L = f(\frac{x}{L}, \frac{y}{L}, u)}$
and we are going to study the variation, with respect to the
periodicity parameter $L,$ of the minimal speed
$c_{\mathbb{R}^{N},A_{L},q_{_{L}},f_{L}}^*(e),$ which corresponds to
the reaction-advection-diffusion equation within the
$(L,\cdots,L)-$periodic coefficients
 $A_{_L},\,q_{_L}$ and $f_{_L}:$
\begin{thm}\label{variation of the min speed with respect to period}
Let $\;e = (1, 0,\ldots, 0) \in \mathbb{R}^N.$ An element $z \in
\mathbb{R}^N$ is represented as $z = (x, y)\in$~$
\mathbb{R}\times\mathbb{R}^{N-1}.$ Assume that $A(x,y) = A(y)$ (for
all $(x,y)\in\mathbb{R}^{N}$) and $f(x, y, u)$ satisfy (
(\ref{nonlinearity 1}), (\ref{nonlinearity 2}) and \ref{cA Y}) with
$\omega = \mathbb{R}^{N-1},\; d = N - 1,$ and $l = L_1 =\ldots =
L_{N-1} = 1.$ Assume furthermore, that for all $y \in
\mathbb{R}^{N-1}, \;A(x,y)e =A(y)e= \alpha(y)e,$  where
$y\mapsto\alpha(y)$ is a positive $(1, \ldots , 1)$-periodic
function defined over $\mathbb{R}^{N-1}$ and that, at least, one of $\tilde {e}\cdot A\tilde {e}$ and $\zeta$ is not constant. Let $q$ be an advection
field satisfying (\ref{cq}) and having the form
$\;\displaystyle{q(x,y)=(q_{_{1}}(y),0\ldots,0)}$ for each
$(x,y)\in\mathbb{R}^{N}.$ Consider the reaction-advection-diffusion
problem,
\begin{eqnarray}\label{KPP L }
\begin{array}{c}
\forall\, (t,x,y)\in\mathbb{R}\times\mathbb{R}^{N},\vspace{3pt}\\
    \displaystyle{u_t(t,x,y)=\nabla\cdot(A_{_{L}}(y)\nabla
    u)(t,x,y)+(q_{_{1}})_{_{L}}(y)\partial_{x}u(t,x,y)+f_{_{L}}(x,y,u)},
\end{array}
\end{eqnarray}
 whose coefficients are $(L, \ldots , L)-$periodic with respect to
$(x, y) \in \mathbb{R}^N.$

 Then, the map $\displaystyle{L
\mapsto c_{\mathbb{R}^{N},A_{_L},\,q_{_{L}},f_{_L}}^{*}(e)}\;$ is
increasing in $L>0.$
\end{thm}
\begin{remark}
The assumptions of Theorem \ref{variation of the min speed with respect to period} can not be fulfilled whenever $N=1.$ However, assuming that $N=1$ and that the function $$\displaystyle{\frac{\zeta}{<\zeta>_A}+\frac{<a>_H}{a}}$$ is not identically equal to $2$ (where $a(x)$ is the diffusion factor, $<a>_H$ and $<\zeta>_A$ are, respectively, the harmonic mean of $x\mapsto a(x)$ and arithmetic mean of $x\mapsto\zeta(x)$ over $[0,1]$), it was proved, in {\rm \cite{EHR}}, that $\displaystyle{L
\mapsto c_{\mathbb{R}^{N},a_{_L},\,q_{_{L}},f_{_L}}^{*}(e)}$ is increasing in $L$ when $L$ is close to $0.$ In particular, if $a$ is constant and $\zeta$ is not constant, or if $\mu$ is constant and $a$ is not constant, then $L\mapsto c_{\mathbb{R}^{N},a_{_L},\,q_{_{L}},f_{_L}}^{*}(e)$ is increasing when $L$ is close to $0.$
\end{remark}
Concerning now the variation with respect to the reaction factor
$B,$ we have the following:
\begin{thm}\label{variation of c*(B)over rad B}
Assume that $\Omega=\mathbb{R}\times\omega$ and the coefficients
$A,\,q$ and $f$ satisfy the same assumptions of Theorem
\ref{variation of min speed with respect diffusion}. Let
$e=(1,0\ldots,0)$ and for each $B>0,$ consider the
reaction-advection-diffusion problem
\begin{equation*}
    \left\{
      \begin{array}{l}
 u_t =\nabla\cdot(A(y)\nabla u)\;+\,\sqrt{B}\,q_{_{1}}(y)\,\partial_{x} u\,+\,B f(x,y,u),\quad t\,\in\,\mathbb{R},\;(x,y)\,\in\,\mathbb{R}\times\omega,
 \\
        \nu \cdot A\;\nabla u(t,x,y) =0,\quad
        t\,\in\,\mathbb{R},\;(x,y)\,\in\,\partial\Omega.
      \end{array}
    \right.
\end{equation*}

Then, the map $\displaystyle{B\mapsto \frac{c_{\Omega,
A,\sqrt{B}\,q,Bf}^*(e)}{\sqrt{B}}}$ is increasing in $B>0.$
\end{thm}

As a first note, we mention that Theorem \ref{variation of c*(B)over
rad B} holds also in the case where there is no advection. On the other hand, Berestycki, Hamel and Nadirashvili \cite{BHN1}
proved that the map $\displaystyle{B\mapsto {c_{\Omega,
A,\,q,Bf}^*(e)}}$ is increasing in $B>0$ under the assumptions
(\ref{comega}), (\ref{cA}), (\ref{cq}), (\ref{cf1}), and (\ref{cf2})
which are less strict than the assumptions considered in our present
theorem. However, the present theorem is concerned in the variation
of the map $\displaystyle{B\mapsto \frac{c_{\Omega,
A,\sqrt{B}\,q,Bf}^*(e)}{\sqrt{B}}}$ rather than that of
$\displaystyle{B\mapsto c_{\Omega, A,\,q,Bf}^*(e)}.$

\begin{remark} Owing to the same justifications given after Theorem \ref{lim as L tends to infty thm in presence of q},
 one concludes the importance of taking, in section \ref{variation of min speed w.r.t diffusion and period}, an advection in the form of shear flows.
To study the variations of the minimal speeds as in Theorems \ref{variation of min speed with respect diffusion}, \ref{variation of the min speed with respect to period}
and \ref{variation of c*(B)over rad B}, but in a more general framework (general advection fields, general diffusion, etc...), formula \ref{var} remains
an important tool. However, we will no longer have variational formulations as
 (\ref{-mu(lambda',beta)=min}) below. These problems remains open in the general periodic framework.
\end{remark}

\section{Proofs of the announced results}\label{proofs section}
In this section, we are going to demonstrate the Theorems announced
in sections \ref{as eps tends to zero}, \ref{as M goes to infty},
\ref{reaction}, and \ref{variation of min speed w.r.t diffusion and
period}. We will proceed in 4 subsections,
each devoted to proving the results announced in a
corresponding section.
\subsection{Proofs of Theorems \ref{limit
as eps}, {\ref{lim as eps but in presence of a shear flow}} and \ref{lim as per L tend to infty}}
\paragraph{Proof of Theorem \ref{limit as eps}.}
Under the assumptions of Theorem \ref{limit as eps}, we can apply
the variational formula (\ref{var}) of the minimal speed.
Consequently,
\begin{equation}\label{var form with eps}
  c_{\Omega,\varepsilon
A,0,f}^*(e)\,=\,\displaystyle{\min_{\lambda\,>\,0}\frac{k_{\Omega,e,\,\displaystyle{\varepsilon
A},\,0,\zeta}(\lambda)}{\lambda}},
\end{equation}

where $k_{\Omega,e,\varepsilon A,0,\zeta}(\lambda)$ is the first
eigenvalue (for each $\lambda,\,\varepsilon\,>\,0$) of the
eigenvalue problem
\begin{eqnarray}\label{eigeneq1}
\left\{
    \begin{array}{ll}
     \displaystyle{L_{\Omega,e,\,\displaystyle{\varepsilon A},\,0,\zeta,\lambda}\;\psi}&=\displaystyle{k_{\Omega,e,\,\displaystyle{\varepsilon A},\,0,\zeta}(\lambda)\;\psi(x,y)\;\hbox{over}\;\mathbb{R}\times\omega;}
    \vspace{3pt} \\
  \nu\cdot A\nabla\psi&=0 \quad
  \hbox{on}\;\mathbb{R}\,\times\,\partial\omega,
    \end{array}
  \right.
\end{eqnarray}
 and
 \begin{eqnarray*}
     L_{\Omega,e,\,\displaystyle{\varepsilon A},\,0,\zeta,\lambda}\psi(x,y)\,&=&\,\varepsilon\nabla\cdot\left(A(y)\nabla\psi(x,y)\right)\,-2\,\varepsilon\,\lambda
Ae\cdot\nabla\psi(x,y)\,+\\
&&\,\left[\varepsilon\,\lambda^{2}eA(y)e\,-\lambda\,\varepsilon\nabla\cdot(A(y)e)+\,\zeta(y)\right]\psi(x,y),
 \end{eqnarray*}
for all  $(x,y)\in\mathbb{R}\times\omega.$

Initially, the boundary condition in (\ref{eigeneq1}) is $ \nu\cdot
A\nabla\psi\,=\,\lambda\,\nu\cdot Ae$ on
$\partial\Omega=\mathbb{R}\times\partial\omega;$ where $\nu(x,y)$ is
the unit outward normal at $(x,y)\in\partial\Omega.$ However,
$\Omega=\mathbb{R}\times\omega$ is invariant in the direction of $e$
which is that of $Ae$ in both alternatives (\ref{alt 1}) and
(\ref{alt 2}). Consequently, $\nu\cdot Ae\equiv0$ on
$\partial\Omega.$

We recall that for all $\lambda>0,$ and for all $\varepsilon>0,$ we
have $k_{\Omega,e,\,\displaystyle{\varepsilon
A},\,0,\zeta}(\lambda)\,>\,0.$ Also, the first eigenfunction of
(\ref{eigeneq1}) is positive over
$\overline{\Omega}=\,\mathbb{R}\times\overline{\omega},$ and it is
unique up to multiplication by a non zero constant.

In our present setting, whether in (\ref{alt 1}) or (\ref{alt 2})
and due to the assumption (\ref{nonlinearity 2}), one concludes that
the coefficients
 in $ \displaystyle{L_{\Omega,e,\,\displaystyle{\varepsilon A},\,0,\zeta,\lambda}}$
are independent of $x.$ Moreover, in both alternatives (\ref{alt 1})
and (\ref{alt 2}), the direction of $Ae$ is the same of
$e=(1,0,\cdots,0).$ On the other hand, since
$\Omega=\mathbb{R}\times\omega,$ then for each
$(x,y)\in\partial\Omega,$ we have $\nu(x,y)=(0;\nu_{\omega}(y)),$
where $\nu_{\omega}(y)$ is the outward unit normal on
$\partial\omega$ at $y.$ Consequently, the first eigenfunction of
(\ref{eigeneq1}) is independent of $x$ and \textbf{the eigenvalue
problem (\ref{eigeneq1}) is reduced to}
\begin{eqnarray}\label{eigeneq2}
\left\{
    \begin{array}{rl}
   \displaystyle{ L_{\Omega,e,\,\displaystyle{\varepsilon
   A},0,\zeta,\lambda}\phi}:&=\,\varepsilon\,\nabla\cdot\left(A(y)\nabla\phi(y)\right)\,+\,\left[\varepsilon\,\lambda^{2}eA(y)e\,+\,\zeta(y)\right]\phi(y)\,\\
   &=\,\displaystyle{k_{\Omega,e,\,\displaystyle{\varepsilon A},\,0,\zeta}(\lambda)}\,\phi\quad\hbox{over}\;\omega; \\
  \nu(x,y)\cdot A(y)\nabla\phi(y)&=\left(0;\nu_{\omega}(y)\right)\cdot A(y)\nabla\phi(y)=0 \quad \hbox{on}\;\mathbb{R}\times\partial\omega,
    \end{array}
  \right.
\end{eqnarray}
where $\phi\,=\,\phi(y)$ is positive over $\overline{\omega},$
$L-$periodic (since the domain $\omega$ and the coefficients of
$L_{\Omega,e,\,\displaystyle{\varepsilon A},0,\zeta,\lambda}$ are
$L-$periodic), unique up to multiplication by a constant, and
belongs to $C\,^{2}(\overline{\omega}).$

In the case where $d\geq1,$ let $C\subseteq\mathbb{R}^{N-1}$ denote
the periodicity cell of $\omega.$ Otherwise, $d=0$ and one takes
$C=\omega.$ In both cases, $C$ is bounded. Multiplying the first
line of (\ref{eigeneq2}) by $\phi,$ and integrating by parts over
$C,$ one gets
\begin{equation}\label{_k1}
    -\,\displaystyle{k_{\Omega,e,\,\displaystyle{\varepsilon A},\,0,\zeta}(\lambda)}\,=\,\displaystyle{\frac{\varepsilon\displaystyle{\int_{C}\,\nabla\phi\cdot A(y)\nabla\phi\,dy}-\displaystyle{\int_{C}\left[\varepsilon\lambda^{2}eA(y)e\,+\,\zeta(y)\right]\,\phi^{2}(y)\,dy}}{\displaystyle{\int_{C}{\phi^2}(y)\,dy}}.}
\end{equation}

One also notes that, in this present setting, the operator
$L_{\Omega,e,\,\varepsilon A,0,\zeta,\lambda}$ is self-adjoint and
its coefficients are $(L_1,\ldots,L_d)-$periodic with respect
$(y_1,\ldots,y_d).$ Consequently,
$-\,\displaystyle{k_{\Omega,e,\,\displaystyle{\varepsilon
A},\,0,\zeta}(\lambda)}$ has the following variational
characterization:
\begin{equation}\label{k2}
    -\,\displaystyle{k_{\Omega,e,\,\displaystyle{\varepsilon A},\,0,\zeta}(\lambda)}=\displaystyle{\min_{\varphi\in
    H^{1}(C)\setminus\{0\}}\frac{\varepsilon\displaystyle{\int_{C}\,\nabla\varphi\cdot A(y)\nabla\varphi\,dy}-\displaystyle{\int_{C}\left[\varepsilon\lambda^{2}eA(y)e\,+\,\zeta(y)\right]\,\varphi^{2}(y)\,dy}}{\displaystyle{\int_{C}{\varphi^2}(y)\,dy}}.}
\end{equation}

In what follows, we will assume that (\ref{alt 1}) is the
alternative that holds. That is, $eAe=\alpha$ is constant. The proof
can be imitated easily whenever we assume that (\ref{alt 2}) holds.

The function $y\mapsto\zeta(y)$ is continuous and
$(L_1,\ldots,L_d)-$periodic over $\overline{\omega},$ whose
periodicity cell $C$ is a bounded subset of $\mathbb{R}^{N-1}$
(whether $d=0$ or $d\geq1$). Let
$y_0\in\overline{C}\subseteq\,\overline{\omega}$ such that
$\displaystyle{\max_{y\in\overline{w}}\,\zeta(y)=\zeta(y_0)}$
(trivially, this also holds when $\zeta$ is constant). Consequently,
we have
$$\forall\;\varphi\in \,H^{1}(C)\setminus\{0\},\;\displaystyle{\frac{\varepsilon\displaystyle{\int_{C}\,\nabla\varphi\cdot\,A\nabla\varphi}-\displaystyle{\int_{C}(\varepsilon\alpha\lambda^{2}\,+\,\zeta(y))\varphi^{2}}}{\displaystyle{\int_{C}{\varphi^2}(y)\,dy}}\,\geq\,-\left[\varepsilon\alpha\lambda^{2}\,+\,\zeta(y_0)\right].}$$
This yields that
\begin{equation}\label{k<}
\forall\,\varepsilon>0,\,\forall\,\lambda>0,\,-\,\displaystyle{k_{\Omega,e,\,\displaystyle{\varepsilon
A},\,0,\zeta}(\lambda)}\,\geq\,-\left[\varepsilon\alpha\lambda^{2}\,+\,\zeta(y_0)\right].
\end{equation}
Consequently,
\begin{equation}\label{m}
    \displaystyle{\forall\,\varepsilon>0,\,\forall\,\lambda>0,\,\frac{\displaystyle{k_{\Omega,e,\,\displaystyle{\varepsilon A},\,0,\zeta}(\lambda)}}{\lambda}\,\leq\,\lambda\,\alpha\varepsilon\,+\frac{\zeta(y_0)}{\lambda}}.
\end{equation}

However, the function
$\displaystyle{{\lambda}\,\mapsto\,\lambda\alpha\varepsilon\,+\frac{\zeta(y_0)}{\lambda}}$
attains its minimum, over $\mathbb{R}^{+},$ at
$\lambda(\varepsilon)\,=\,\displaystyle{\sqrt{\frac{\zeta(y_0)}{\alpha\varepsilon}}}$.
This minimum is equal to
$2\sqrt{\zeta(y_0)}\times\sqrt{\alpha\,\varepsilon}.$ From (\ref{m}), we conclude that
$$\displaystyle{\frac{\displaystyle{k_{\Omega,e,\,\displaystyle{\varepsilon
A},\,0,\zeta}(\lambda(\varepsilon))}}{\lambda(\varepsilon)}\,\leq\,2\sqrt{\alpha\varepsilon}\sqrt{\zeta(y_0)}}.$$

Finally, (\ref{var}) implies that
$\displaystyle{c_{\Omega,\varepsilon
A,0,f}^*(e)\,=\,\displaystyle{\min_{\lambda\,>\,0}\frac{k_{\Omega,e,\,\displaystyle{\varepsilon
A},\,0,\zeta}(\lambda)}{\lambda}}\,\leq\,2\sqrt{\alpha\,\varepsilon}\sqrt{\zeta(y_0),}}$
or equivalently
\begin{equation}\label{c(eps) over eps <}
   \forall\varepsilon\,>0,\, \displaystyle{\frac{c_{\Omega,\varepsilon
A,0,f}^*(e)}{\sqrt{\varepsilon}}\,\leq\,2\sqrt{\alpha}\sqrt{\zeta(y_0).}}
\end{equation}

We pass now to prove the other sense of the inequality for
$\displaystyle{\liminf_{\varepsilon\rightarrow0^+}\displaystyle{\frac{c_{\Omega,\varepsilon
A,0,f}^*(e)}{\sqrt{\varepsilon}}}}.$ We will consider formula
(\ref{k2}), and then organize a suitable function $\psi$ which leads
us to a lower bound of
$\displaystyle{\liminf_{\varepsilon\rightarrow0^+}\frac{c_{\Omega,\varepsilon
A,0,f}^*(e)}{\sqrt{\varepsilon}}}.$

We have $\zeta(y_0)\,>\,0.$ Let $\delta$ be such that
$0<\delta<\zeta(y_{0}).$ Thus
$\displaystyle{0<\zeta(y_{0})-\delta<\max_{\overline{\omega}}\,\zeta(y)}.$
The continuity of $\zeta,\;$ over
$\overline{C}\subseteq\overline{\omega},$ yields that there exists
an open and bounded set $U\subset\overline{C}$ such that
\begin{equation}\label{setU}
    \forall\,y\in\,\overline{U},~\zeta(y_{0})-\delta\leq\zeta(y).
\end{equation}

Designate by $\psi,$ a function in $\mathcal{D}(C)$ (a
$C^{\infty}(C)$ function whose support is compact), with
supp$\,\psi\,\subseteq\,\overline{U},$ and $\displaystyle{\int_{U}\psi^{2}=1.}$ One
will have,
\begin{eqnarray*}
\begin{array}{ll}
\forall\lambda>0,\,\forall\,\varepsilon\,>0,&\vspace{3 pt}\\
   \nonumber  -\, k_{\Omega,e,\,\displaystyle{\varepsilon
A},\,0,\zeta}(\lambda)&\leq\,\displaystyle{\varepsilon\displaystyle{\int_{U}\,\nabla\psi\cdot A(y)\nabla\psi\,dy}-\displaystyle{\int_{U}\left[\varepsilon\lambda^{2}eA(y)e\,+\,\zeta(y)\right]\,\psi^{2}(y)\,dy}}\vspace{3 pt}\\
                            \nonumber  &\leq\,\displaystyle{\varepsilon\displaystyle{\int_{U}\,\nabla\psi\cdot A(y)\nabla\psi\,dy}-\left[\varepsilon\lambda^{2}\alpha\,+\,\zeta(y_0)-\delta\right]\displaystyle{\int_{U}\,\psi^{2}(y)\,dy}}\vspace{3 pt}\\
\nonumber&\leq\displaystyle{\varepsilon\displaystyle{\int_{U}\,\alpha_2|\nabla\psi|^2}\,-\left[\varepsilon\lambda^{2}\alpha\,+\,\zeta(y_0)-\delta\right]},\;\hbox{by}\;(\ref{cA
Y}),
\end{array}
 \end{eqnarray*}
or equivalently
\begin{equation}\label{k() over lambda >}
\displaystyle{\frac{k_{\Omega,e,\,\displaystyle{\varepsilon
A},\,0,\zeta}(\lambda)}{\lambda}\,\geq\,\lambda\alpha\varepsilon\,+\,\displaystyle{\frac{1}{\lambda}}\;\beta(\varepsilon),}
\end{equation}
where
$\beta(\varepsilon)=\zeta(y_0)-\delta-\varepsilon\,\displaystyle{\displaystyle{\int_{U}\alpha_2|\nabla\psi|^{2}}}.$
Choosing
$0<\varepsilon<\,\displaystyle{\frac{\zeta(y_0)-\delta}{\displaystyle{\alpha_2\int_{U}|\nabla\psi|^{2}}}}\;$
(this is possible), we get  $\beta(\varepsilon)\,>0.$

 The map
$\lambda\mapsto\lambda\alpha\varepsilon\,+\,\displaystyle{\frac{1}{\lambda}}\;\beta(\varepsilon)$
attains its minimum, over $\mathbb{R}^+,$ at
$\lambda(\varepsilon)\,=\,\displaystyle{\sqrt{\frac{\beta(\varepsilon)}{\varepsilon\alpha}}}.$
This minimum is equal to
$2\displaystyle{\sqrt{\varepsilon\,\alpha}\sqrt{\beta(\varepsilon)}}.$

Now, referring to formula (\ref{k() over lambda >}), one gets
 $$\hbox{For $\varepsilon\;$ small enough,}\quad\displaystyle{\frac{k_{\Omega,e,\,\displaystyle{\varepsilon
A},\,0,\zeta}(\lambda)\,}{\lambda}\,\geq\,2\displaystyle{\sqrt{\varepsilon\alpha}\,\sqrt{\beta(\varepsilon)}}\quad
\hbox{for all $\lambda>0.$}}$$ Together with (\ref{var}), we
conclude that
\begin{equation}
\, \hbox{for $\varepsilon$ small enough,}\quad\displaystyle{\frac{
c_{\Omega,\varepsilon A,0,f}^*(e)}{\sqrt{\varepsilon}}} \;\geq\;
2\sqrt{\beta(\varepsilon)}\sqrt{\alpha}.
\end{equation}
 Consequently, \begin{eqnarray}
                \nonumber \displaystyle{\liminf_{\varepsilon\rightarrow0^+}\frac{ c_{\Omega,\varepsilon
A,0,f}^*(e)\,}{\sqrt{\varepsilon}}}\,&\geq&\,\displaystyle{\liminf_{\varepsilon\rightarrow0^+}2\sqrt{\beta(\varepsilon)}\sqrt{\alpha}\,} \vspace{3 pt}\\
                \nonumber &=&2\sqrt{\zeta(y_0)-\delta}\sqrt{\alpha}\quad\hbox{(since $\psi$ is independent of
$\varepsilon$)},
               \end{eqnarray}
and this holds for all $0<\delta<\zeta(y_{0}).$ Therefore, one can
conclude that
\begin{equation}\label{liminf c(epsilon) over rad eps}
    \displaystyle{\liminf_{\varepsilon\rightarrow0^+}\frac{ c_{\Omega,\varepsilon
A,0,f}^*(e)\,}{\sqrt{\varepsilon}}}\,\geq\,2\sqrt{\alpha}\sqrt{\zeta(y_0).}
\end{equation}

Finally, the inequalities (\ref{c(eps) over eps <}) and (\ref{liminf
c(epsilon) over rad eps}) imply that
$\displaystyle{\lim_{\varepsilon\rightarrow0^+}\frac{
c_{\Omega,\varepsilon A,0,f}^*(e)\,}{\sqrt{\varepsilon}}}$ exists,
and it is equal to
$$\displaystyle{2\sqrt{\alpha}\sqrt{\zeta(y_0)}=2\sqrt{\max_{\overline{\omega}}eA(y)e}\sqrt{\max_{\overline{\omega}}\zeta(y)}}.$$

We note that the same ideas of this proof can be easily applied in
the case where the assumption (\ref{alt 2}) holds. In (\ref{alt 2}),
we have $\zeta$ is constant; however, $eAe$ is not in general.
Meanwhile the converse is true in the case (\ref{alt 1}).  The
little difference is that, in the case of (\ref{alt 2}), we choose
the subset $U$ (of the proof done above) around the point $y_0$
where $eAe$ attains its maximum and then we continue by the same way
used above. \hfill$\Box$\\

\textbf{Proof of Theorem \ref{lim as eps but in presence of a shear flow}.} We have
\begin{equation}\label{var form with eps and q}
  c_{\Omega,\varepsilon
A,0,f}^*(e)\,=\,\displaystyle{\min_{\lambda\,>\,0}\frac{k_{\Omega,e,\,\displaystyle{\varepsilon
A},\,q,\zeta}(\lambda)}{\lambda}},
\end{equation}
where (due to the facts that $q$ is a shear flow, $e=(1,0,\cdots,0)$ and $e$ is an eigenvector of the matrix $A(y)$ for all $y\in\overline{\omega})$ $k_{\Omega,e,\,\displaystyle{\varepsilon
A},\,q,\zeta}(\lambda)$ is the principal eigenvalue of the problem
 \begin{eqnarray*}
    \left\{
    \begin{array}{ll}
      \displaystyle{L_{\Omega,e,\,\varepsilon A,\,q,\zeta,\lambda}}\psi(x,y)&=\;\displaystyle{k_{\Omega,e,\,\varepsilon A,\,q,\zeta}(\lambda)}\,\psi(x,y)\quad\hbox{over}\;\mathbb{R}\times\omega;\vspace{3 pt}\\
  \nu\cdot A\nabla\psi&=0 \quad
  \hbox{on}\;\mathbb{R}\,\times\,\partial\omega,
    \end{array}
  \right.
  \end{eqnarray*}
with \begin{equation}\label{eigen equation eps A with q}
\begin{array}{lll}
\displaystyle{L_{\Omega,e,\varepsilon
A,\,q,\zeta,\lambda}}\,\psi&=&\varepsilon\nabla\cdot\left(A(y)\nabla\psi\right)-2\varepsilon\lambda\,\alpha(y)\,\partial_{x}\psi+\,q_{_{1}}(y)\partial_{x}\psi\vspace{4pt}\\
&&+\left[\varepsilon\,\lambda^{2}eA(y)e\,-\lambda q_{_{1}}(y)+\zeta(y)\right]\psi\;\hbox {over}\;\mathbb{R}\times\omega.
\end{array}
\end{equation}

The uniqueness of the principal eigenfunction $\psi$ up to multiplication by a constant, yields that one can choose $\psi$ independent of $x.$ Hence, the elliptic operator
$\displaystyle{L_{\Omega,e,\varepsilon
A,\,q,\zeta,\lambda}}$ can be reduced to the symmetric operator
$$\displaystyle{L_{\Omega,e,\varepsilon
A,\,q,\zeta,\lambda}}\,\psi=\varepsilon\nabla\cdot\left(A(y)\nabla\psi\right)+\left[\varepsilon\,\lambda^{2}eA(y)e\,-\lambda q_{_{1}}(y)+\zeta(y)\right]\psi.$$
 Consequently,
\begin{equation}\label{rayleigh}
    \begin{array}{ll}
\forall \lambda>0,\,\forall \varepsilon>0,\quad-\displaystyle{k_{\Omega,e,\,\varepsilon A,\,q,\zeta}(\lambda)}=
\vspace{3 pt}\\
\displaystyle{\min_{\varphi\in
    H^{1}(C)\setminus\{0\}}\frac{\displaystyle{\varepsilon\int_{C}\nabla\varphi\cdot A(y)\nabla\varphi dy}+\lambda\int_{C}q_{_{1}}(y)\varphi^{2}-\displaystyle{{\int_{C}\left[{\lambda}^{2}\varepsilon eA(y)e\,+\,\zeta(y)\right]\,\varphi^{2}(y)\,dy}}}{\displaystyle{\int_{C}\varphi^{2}(y)\,dy}}.}\\
\end{array}
\end{equation}
Formula (\ref{rayleigh}) yields that
$$\forall \lambda>0,\,\forall \varepsilon>0,~-\displaystyle{k_{\Omega,e,\,\varepsilon A,\,q,\zeta}(\lambda)}\geq -\lambda\ds\max_{y\in\overline{\omega}}\left(-q_1(y)\right)-{\lambda}^{2}\varepsilon \ds\max_{y\in\overline{\omega}}eA(y)e-\max_{y\in\overline{\omega}}\zeta(y),$$
or equivalently
$$\ds{\forall \lambda>0,\,\forall \varepsilon>0,~\displaystyle{\frac{k_{\Omega,e,\,\varepsilon A,\,q,\zeta}(\lambda)}{\lambda}}\leq \ds\max_{y\in\overline{\omega}}\left(-q_1(y)\right)+{\lambda}\varepsilon \ds\max_{y\in\overline{\omega}}eA(y)e+\ds\frac{\ds\max_{y\in\overline{\omega}}\zeta(y)}{\lambda}.}$$
Putting $\lambda=\lambda(\varepsilon)=\ds{\sqrt{\frac{\max_{y\in\overline{\omega}}\zeta(y)}{\varepsilon\max_{y\in\overline{\omega}e\cdot A(y)e}}}>0}$
into the last inequality yields that
$$\displaystyle{\min_{\lambda\,>\,0}\frac{k_{\Omega,e,\,\displaystyle{\varepsilon
A},\,q,\zeta}(\lambda)}{\lambda}}\leq \max_{y\in\overline{\omega}}\left(-q_1(y)\right)+ 2\sqrt{\epsilon}\sqrt{\max_{y\in\overline{\omega}}e\cdot A(y)e}\sqrt{\max_{y\in\overline{\omega}}\zeta(y)},$$
and hence,
\begin{equation}\label{lim sup as eps with q}
\displaystyle{\limsup_{\varepsilon\rightarrow0^+}
c_{\Omega,\varepsilon A,q,f}^*(e)\leq\max_{y\in\overline{\omega}}\left(-q_1(y)\right).}
\end{equation}

Now, we take $y_0\in C$ ($C$ is the periodicity cell of $\omega$) such that $\max_{y\in\overline{\omega}}\left(-q_1(y)\right)=-q_1(y_0)>0$ (since $q$ is
periodic with respect to $y,$ $q_1\not \equiv 0$ and $q_1$ has a zero average)
and we take $\delta>0$ such $-q_1(y_0)-\delta>0.$ It follows, from the continuity of $q_1,$ that there exists an open subset $U\subset C$ such that
$y_0\in U$ and $$\forall y\in \overline{U},~-q_1(y)\geq \max_{y\in\overline{\omega}}\left(-q_1(y)\right)-\delta.$$

Let $\psi$ be a function in $\mathcal{D}(C)$ with
supp$\,\psi\,\subseteq\,\overline{U},$ and $\displaystyle{\int_{U}\psi^{2}=1.}$ Referring to (\ref{rayleigh}), it follows that

\begin{equation}\label{k() over lambda > with q}
\forall\lambda>0,\forall \epsilon>0, ~~\displaystyle{\frac{k_{\Omega,e,\,\displaystyle{\varepsilon
A},\,q,\zeta}(\lambda)}{\lambda}\,\geq\,-q_1(y_0)-\delta+\lambda\varepsilon\ds\min_{y\in\overline{\omega}}e\cdot Ae\,+\,\displaystyle{\frac{1}{\lambda}}\;\beta(\varepsilon),}
\end{equation}
where
$\beta(\varepsilon)=\ds\min_{y\in\overline{\omega}}\zeta(y)-\varepsilon\,\displaystyle{\displaystyle{\int_{U}\alpha_2|\nabla\psi|^{2}}}>0$ for a small enough
$\epsilon>0$ ($\alpha_2>0$ is the constant appearing in (\ref{cA Y})).

It follows from (\ref{k() over lambda > with q}) that
$$\forall\lambda>0,\forall \epsilon>0, ~~\displaystyle{\frac{k_{\Omega,e,\,\displaystyle{\varepsilon
A},\,q,\zeta}(\lambda)}{\lambda}\,\geq\,-q_1(y_0)-\delta+2\sqrt{\epsilon}\sqrt{\ds\min_{y\in\overline{\omega}}e\cdot Ae}\sqrt{\beta(\epsilon)}}.$$
Together with (\ref{var form with eps and q}), and since $\delta>0$ is arbitrary, one gets
\begin{equation}\label{liminf c(epsilon) with q}
    \displaystyle{\liminf_{\varepsilon\rightarrow0^+} c_{\Omega,\varepsilon
A,q,f}^*(e)\,\,\geq\,-q_1(y_0)=\max_{y\in \overline{\omega}}(-q_{1}(y)).}
\end{equation}
Finally, (\ref{lim sup as eps with q}) and (\ref{liminf c(epsilon) with q}) complete the proof of Theorem \ref{lim as eps but in presence of a shear flow}.\hfill$\Box$\\
\vskip 0.2cm
\textbf{ Proof of Theorem \ref{lim as per L tend to infty}.}
Consider the change of variables
$$v(t,x,y)=u(t,Lx,Ly),\quad(t,x,y)\in\mathbb{R}\times\mathbb{R}\times\mathbb{R}^{N-1}.$$
The function $u$ satisfies (\ref{KPP L with u}) if and only if $v$
satisfies
\begin{eqnarray}\label{KPP L with v}
v_t(t,x,y)=\displaystyle{\frac{1}{L^{^{2}}}\;\nabla\cdot(A(y)\nabla
    v)(t,x,y)+f(x,y,v)\;\hbox{over}\;\mathbb{R}\times\mathbb{R}\times\mathbb{R}^{N-1}}.
\end{eqnarray}
Consequently,
\begin{equation}\label{beta * in terms of c*}
\displaystyle{  \forall\,L>0,\;c_{\mathbb{R}^{N},
\displaystyle{A_{_{L}},0,f_{_{L}}}}^*(e)\,=\,L\,c_{\mathbb{R}^{N},\frac{1}{L^{^{2}}}
A,0,f}^*(e)}
\end{equation}

Taking
$\varepsilon=\displaystyle{{1}/{L^{^{2}}}},$ and applying Theorem \ref{limit as eps} to problem (\ref{KPP L with v}), one
then has
\begin{equation}\label{lim beta*}
\displaystyle{\lim_{L\rightarrow+\infty}\frac{c_{\mathbb{R}^{N},\frac{1}{L^{^{2}}} A,0,f}^*(e)}{\displaystyle{\sqrt{\frac{1}{L^{^{2}}}}}}\,=\,\lim_{\varepsilon\rightarrow0^+}\,\frac{c_{\mathbb{R}^{N},
\varepsilon A,0,f}^*(e)}{\sqrt{\varepsilon}}\,=\,2\,\sqrt{\max_{y\,\in\,\mathbb{R}^{N-1}}\zeta(y)}\sqrt{\max_{y\,\in\,\mathbb{R}^{N-1}}eA(y)e}}.
\end{equation}

Finally, (\ref{beta * in terms of c*}) together with (\ref{lim
beta*}) complete the proof of Theorem \ref{lim as per L tend to
infty}.\hfill$\Box$\vskip0.45cm

\subsection{Proofs of Theorems \ref{lim as M} and \ref{lim as per L tend to 0}}
 \paragraph{Proof of Theorem \ref{lim as M}.} The proof will be
 divided into three steps:

 \underline{Step 1.} According to Theorem \ref{varthm}, and since $\nu\cdot A\tilde{e}=0$ on $\partial\Omega,$ the minimal
 speeds $c_{\Omega,MA,\displaystyle{M^{\gamma}q},f}^*(e)$ are given by:
 \begin{equation*}
   \displaystyle{ \forall
   M>0,\;c_{\Omega,MA,\displaystyle{M^{\gamma}\,q},f}^*(e)=\min_{\lambda>0}\frac{k_{\Omega,e,\,\displaystyle{M
A},\,\displaystyle{M^{\gamma}\,q},\zeta}(\lambda)}{\lambda}},
 \end{equation*}
where $k_{\Omega,e,\,\displaystyle{M
A},\,\displaystyle{M^{\gamma}\,q},\zeta}(\lambda)$ and
$\psi^{\lambda,M}$ denote the unique eigenvalue and the positive
$L$-periodic eigenfunction of the problem
\begin{eqnarray*}\label{L lambda M}
\begin{array}{cc}
\displaystyle{M\nabla\cdot(A\nabla\psi^{\lambda,M})-2M\lambda\tilde{e}\cdot
A\nabla\psi^{\lambda,M}+\displaystyle{M^{\gamma}q\cdot\nabla\psi^{\lambda,M}}+[\lambda^{2}M\,\tilde{e}A\tilde{e}-\lambda
M^{\gamma}q\cdot\tilde{e}+\zeta]\psi^{\lambda,M}}\vspace{4 pt}\\
=k_{\Omega,e,\,\displaystyle{M
A},\,\displaystyle{M^{\gamma}\,q},\zeta}(\lambda)\psi^{\lambda,M}\;\;
\hbox{in}\; \Omega,
\end{array}
\end{eqnarray*}
with $\nu\cdot
A\nabla\psi^{\lambda,M}=0$ on $\partial\Omega.$

For each $\lambda>0$ and $M>0,$ let $\lambda^{\,'}=\lambda\sqrt{M},$
and let $\;k_{\Omega,e,\,\displaystyle{M
A},\,\displaystyle{M^{\gamma}\,q},\zeta}(\lambda)=\mu(\lambda^{'},M).\;$
Consequently,
\begin{equation}\label{var with mu}
 \displaystyle{ \forall
   M>0,\;\frac{c_{\Omega,MA,\displaystyle{M^{\gamma}\,q},f}^*(e)}{\sqrt{M}}=\min_{\lambda^{\,'}>0}\frac{\mu(\lambda^{\,'},M)}{\lambda^{\,'}}},
\end{equation}
where $\mu(\lambda^{\,'},M)$ and $\psi^{\lambda^{\,'},M}$ are the
first eigenvalue and the unique, positive $L-$periodic (with respect
to $x$ ) eigenfunction of
\begin{eqnarray}\label{L lambda ' M}
\begin{array}{cc}
\displaystyle{M\nabla\cdot(A\nabla\psi^{\lambda^{'},M})-2\lambda{'}\sqrt{M}\tilde{e}\cdot
A\nabla\psi^{\lambda^{'},M}+\displaystyle{M^{\gamma}q\cdot
\nabla\psi^{\lambda^{'},M}}}\vspace{4 pt}\\
+\left[{\lambda^{'}}^{2}\tilde{e}A\tilde{e}-\displaystyle{\frac{\lambda^{'}}{
M^{^{\frac{1}{2}-\gamma}}}}\,q\cdot\tilde{e}+\zeta\right]\psi^{\lambda^{'},M}\;=\;\mu(\lambda^{'},M)\psi^{\lambda^{'},M}\hbox{ in } \Omega,
\end{array}
\end{eqnarray}
with $\nu\cdot A\nabla\psi^{\lambda^{\,'},M}=0$ on
$\partial\Omega.$

Owing to the uniqueness, up to multiplication by positive constants,
of the first eigenfunction of (\ref{L lambda ' M}), one may assume
that:
 \begin{equation}\label{norm L2 of psi s}
    \forall \,\lambda^{'}>0,\;\forall
    \,M>0,\;||\psi^{\lambda^{'},M}||_{L^{2}(C)}\,=1.
 \end{equation}

Moreover, for each
$\displaystyle{M>0,\;\min_{\lambda^{\,'}>0}\frac{\mu(\lambda^{\,'},M)}{\lambda^{\,'}}}$
is attained at $\lambda^{\,'}_{M}>0.$ Thus,
\begin{equation}\label{var with mu and lambda'}
 \displaystyle{ \forall
   M>0,\;\frac{c_{\Omega,MA,\displaystyle{M^{\gamma}\,q},f}^*(e)}{\sqrt{M}}=\min_{\lambda^{\,'}>0}\frac{\mu(\lambda^{\,'},M)}{\lambda^{\,'}}=\frac{\mu(\lambda^{\,'}_{M},M)}{\lambda^{\,'}_{M}}}.
\end{equation}

The above characterization of
$\displaystyle{{c_{\Omega,MA,\displaystyle{M^{\gamma}\,q},f}^*(e)}/{\sqrt{M}}}$
will be used in the next steps in order to prove that
$\displaystyle{\liminf_{M\rightarrow+\infty}{c_{\Omega,MA,\displaystyle{M^{\gamma}\,q},f}^*(e)}/{\sqrt{M}}}$
(resp.
$\displaystyle{\limsup_{M\rightarrow+\infty}{c_{\Omega,MA,\displaystyle{M^{\gamma}\,q},f}^*(e)}/{\sqrt{M}}}$
) is greater than (resp. less than)
$\displaystyle{2\sqrt{\mi_{\!\!\!\!_C}\tilde{e}A\tilde{e}(x,y)dx\,dy}\,\sqrt{\mi_{\!\!\!\!_C}\zeta(x,y)dx\,dy}};$
and hence, complete the proof.\vskip0.35cm

\underline{Step 2.} Fix $\lambda^{\,'}>0$ and $M>0.$ We divide
(\ref{L lambda ' M}) by $\psi^{\lambda^{'},M}$ then, using the facts
$\nabla.A\tilde{e}\equiv0$ in $\Omega$ and $\nu\cdot A\tilde{e}=0$
on $\partial\Omega,$ we integrate by parts over the periodicity cell
$C.$ It follows from (\ref{cq}) and the $L-$periodicity of
$A,\,\zeta\;$ and $\psi^{\lambda^{'},M}\;$ that

\begin{equation}\label{divide by psi and integrate}
    \displaystyle{\int_{C}\frac{\nabla\psi^{\lambda^{'},M}\cdot A\nabla\psi^{\lambda^{'},M}}{{\left(\psi^{\lambda^{'},M}\right)}^{2}}+{\lambda^{\,'}}^{\,2}\int_{C}\tilde{e}A\tilde{e}+\int_{C}\zeta=\mu(\lambda^{\,'},M)|C|},
\end{equation}
where $|C|$ denotes the Lebesgue measure of $C.$ Let $$
\displaystyle{m_0=\mi_{\!\!\!\!_C}\tilde{e}A\tilde{e}=\frac{1}{|C|}\int_{C}\tilde{e}A(x,y)\tilde{e}\,dx\,dy}\quad\hbox{ and }\quad\displaystyle{m=\mi_{\!\!\!\!_C}\zeta(x,y)\,dx\,dy}.$$ One concludes that\begin{eqnarray*}
                       \forall\lambda^{\,'}>0,\,\forall
M>0,\quad\mu(\lambda^{\,'},M)\;\geq\;{\lambda^{\,'}}^{\,2}\mi_{\!\!\!\!_C}\tilde{e}A\tilde{e}+\mi_{\!\!\!\!_C}\zeta={\lambda^{\,'}}^{\,2}m_{0}+m,
                      \end{eqnarray*}
whence\begin{eqnarray}\label{mu(lambda ') over lambda'}
\displaystyle{\forall\lambda^{\,'}>0,\;\forall
M>0,\quad\frac{\mu(\lambda^{\,'},M)}{\lambda^{\,'}}\geq\,\lambda^{\,'}m_{0}+\frac{m}{\lambda^{\,'}}}.
      \end{eqnarray}

The right side of (\ref{mu(lambda ') over lambda'}) attains its
minimum over $\mathbb{R}^{+}$ at
$\displaystyle{\lambda^{\,'}_{0}=\sqrt{\frac{m}{m_0}}}.$ This
minimum is equal to $\displaystyle{2\sqrt{m_{0}m}}.$

Consequently, for any $M>0,$ $\displaystyle{\frac{c_{\Omega,MA,\displaystyle{M^{\gamma}\,q},f}^*(e)}{\sqrt{M}}=\min_{\lambda^{\,'}>0}\frac{\mu(\lambda^{\,'},M)}{\lambda^{\,'}}\geq
2\sqrt{m_{0}m}}.$ This yields that
\begin{equation}\label{liminf as M tends }
\displaystyle{\liminf_{M\rightarrow+\infty}\;\frac{c_{\Omega,MA,\displaystyle{M^{\gamma}\,q},f}^*(e)}{\sqrt{M}}\;\geq\,2\sqrt{\mi_{\!\!\!\!_C}\tilde{e}A\tilde{e}(x,y)dx\,dy}\,\sqrt{\mi_{\!\!\!\!_C}\zeta(x,y)dx\,dy}}.
\end{equation}

\underline{Step 3.} Fix $\lambda^{\,'}>0$ and $M>0.$ Multiply
(\ref{L lambda ' M}) by $\psi^{\lambda^{'},M}$ and integrate by
parts over $C.$ Owing to the $L-$periodicity of $\Omega,\,A,\,\zeta$
and $\psi^{\lambda^{'},M},$ and due to the facts that
$\displaystyle{\int_{C}{\left(\psi^{\lambda^{'},M}\right)}^2=1,\,\nabla\cdot
A\tilde{e}\equiv 0}$ in $\Omega,$ and that $\nu\cdot A\tilde{e}=0$
on $\partial\Omega,$ together with (\ref{cq}), one gets
\begin{eqnarray}\label{multiply by psi lambda ' and integrate}
 \begin{array}{c}
\displaystyle{-M\int_{C}\nabla\psi^{\lambda^{'},M}\cdot A\nabla\psi^{\lambda^{'},M}+{\lambda^{\,'}}^{\,2}\int_{C}\tilde{e}A\tilde{e}\,{\left(\psi^{\lambda^{'},M}\right)}^{2}+\int_{C}\zeta\,{\left(\psi^{\lambda^{'},M}\right)}^{2}}\vspace{4 pt}\\
-\,\displaystyle{\frac{\lambda^{'}}{M^{^{\frac{1}{2}-\gamma}}}\int_{C}q\cdot\tilde{e}\,{\left(\psi^{\lambda^{'},M}\right)}^{2}}\;=\;\mu(\lambda^{\,'},M),
 \end{array}
\end{eqnarray}
whence
\begin{equation*}\label{}
    \forall \lambda^{\,'}>0,\,\forall M>0,\;
    0<\mu(\lambda^{\,'},M)\leq \,{\lambda^{\,'}}^{\,2}\alpha+\beta+\,\displaystyle{\frac{\lambda^{'}}{M^{^{\frac{1}{2}-\gamma}}}||\left(q\cdot\tilde{e}\right)^{-}||_{\infty}},
\end{equation*}
where
$\displaystyle{\alpha=\max_{(x,y)\in\overline{\Omega}}\tilde{e}A\tilde{e}(x,y)}\;$
and
$\;\displaystyle{\beta=\max_{(x,y)\in\overline{\Omega}}\zeta(x,y)}.$
Together with (\ref{mu(lambda ') over lambda'}), one gets
\begin{equation}\label{mu (lambda',M) is bounded below and above }
 \forall \lambda^{\,'}>0,\,\forall M>0,\;
    0<{\lambda^{\,'}}^{2}\,m_{0}+m\leq\mu(\lambda^{\,'},M)\leq
{\lambda^{\,'}}^{\,2}\alpha+\beta+\,\displaystyle{\frac{\lambda^{'}}{M^{^{\frac{1}{2}-\gamma}}}||\left(q\cdot\tilde{e}\right)^{-}||_{\infty}}.
\end{equation}

If $\gamma=\displaystyle{\frac{1}{2}},$ then
$\displaystyle{\frac{\lambda^{'}}{\displaystyle{M^{^{\frac{1}{2}-\gamma}}}}||\left(q\cdot\tilde{e}\right)^{-}||_{\infty}=\lambda^{'}||\left(q\cdot\tilde{e}\right)^{-}||_{\infty}.}$
On the other hand, if $0\leq\gamma<\displaystyle{\frac{1}{2}},$ then
$$\displaystyle{\frac{\lambda^{'}}{M^{^{\frac{1}{2}-\gamma}}}||\left(q\cdot\tilde{e}\right)^{-}||_{\infty}\rightarrow\,0\quad\hbox{as}\quad M\rightarrow+\infty.}$$

Consequently, the right side of (\ref{mu (lambda',M) is bounded
below and above }) is bounded above by a positive constant $B$ which
does not depend on $M$ and $\gamma.$ This yields that
$$\forall \lambda^{'}>0,\quad\displaystyle{
0<\limsup_{M\rightarrow+\infty}\mu(\lambda^{'},M)<+\infty.}$$

On the other hand, it follows from (\ref{cA}) and (\ref{multiply by
psi lambda ' and integrate}) that $\forall \lambda^{\,'}>0,\forall
M>0,$
\begin{eqnarray*}
  \begin{array}{ll}
 0&\leq\displaystyle{\alpha_1\int_{C}|\nabla\psi^{\lambda^{'},M} |^{2}}\leq\displaystyle{\int_{C}\nabla\psi^{\lambda^{'},M}\cdot A\nabla\psi^{\lambda^{'},M}}\vspace{4 pt}\\
 &\leq\displaystyle{\frac{1}{M}\,\left[-\mu(\lambda^{'},M)+{\lambda^{'}}^{2}\displaystyle{\int_{C}\tilde{e}A\tilde{e}\,{\left(\psi^{\lambda^{'},M}\right)}^{2}}+\displaystyle{\int_{C}\zeta\,{\left(\psi^{\lambda^{'},M}\right)}^{2}}-\,\displaystyle{\frac{\lambda^{'}}{M^{^{\frac{1}{2}-\gamma}}}\int_{C}q\cdot\tilde{e}\,{\left(\psi^{\lambda^{'},M}\right)}^{2}}\right]}\vspace{4 pt}\\
&<\,\displaystyle{\frac{B}{M}}.
  \end{array}
\end{eqnarray*}
Meanwhile, $\displaystyle{\lim_{M\rightarrow+\infty}\frac{B}{M}}=0,$
one then gets
\begin{equation}\label{boudedness of psi lambda' in H1}
    \left\{
      \begin{array}{ll}
     \forall\lambda^{\,'}>0,\quad
\displaystyle{\lim_{M\rightarrow+\infty}\int_{C}|\nabla\psi^{\lambda^{'},M}}
    |^{2}=0,\\
\forall\lambda^{\,'}>0,\,\forall
M>0,\quad\displaystyle{\int_{C}\left(\psi^{\lambda^{'},M}\right)^{2}}=1.
      \end{array}
    \right.
\end{equation}

Fix $\lambda^{\,'}>0,$ and let $(M_{n})_{n}$ be a sequence
converging to $+\infty\;$ as $\;n\rightarrow+\infty$ and such that
$\displaystyle{\mu(\lambda^{'},M_n)\rightarrow\,\displaystyle{l^{\,\lambda^{'},(M_n)}}}$
as $\;n\rightarrow+\infty.$ It follows, from (\ref{boudedness of psi
lambda' in H1}), that
$\displaystyle{||\psi^{\lambda^{'},M_n}||_{H^{1}(C)}\rightarrow1}$
as $n\rightarrow+\infty.$ Thus, the sequence
$\displaystyle{(\psi^{\lambda^{'},M_n})_n}$ is bounded in
$H^{1}(C).$ Therefore, there exists a function
$\psi^{\lambda^{'},\infty}\in H^{1}(C)$ such that,  up to extraction
of some subsequence, the functions
$\displaystyle{(\psi^{\lambda^{'},M_n})_n}$ converge in $L^{2}(C)$
strong, $H^{1}(C)$ weak and almost everywhere in $C,$ to the
function $\psi^{\lambda^{'},\infty}.$ Consequently, and owing to
(\ref{boudedness of psi lambda' in H1}), $\psi^{\lambda^{'},\infty}$
satisfies

    \begin{equation}\label{*}
       \displaystyle{\int_{C}{\left(\psi^{\lambda^{'},\infty}\right)}^{2}=1},\;\hbox{and}
    \end{equation}

\begin{equation}\label{**}
  \displaystyle{\left(\int_{C}|\nabla\psi^{\lambda^{'},\infty}|^2\right)^{\frac{1}{2}}\leq\liminf_{M_n\rightarrow+\infty}\left(\int_{C}|\nabla\psi^{\lambda^{'},M_n}|^{2}\right)^{\frac{1}{2}}=0}.
\end{equation}

From (\ref{**}), it follows that for all $\lambda^{'}>0,$ the
function $\psi^{\lambda^{'},\infty}$ is almost everywhere constant
over $C.$ On the other hand, the elliptic regularity applied on
equation (\ref{L lambda ' M}) for $M=M_n,$ implies that
$\forall\lambda^{'}>0,$ the function $\psi^{\lambda^{'},\infty}$ is
continuous over $\overline{C}.$ Consequently, referring to
(\ref{*}), one gets
\begin{equation}\label{***}
\displaystyle{\forall\lambda^{'}>0,\quad\psi^{\lambda^{'},\infty}=\frac{1}{\sqrt{|C|}}}\;\;\hbox{over}\;\;\overline{C}.
\end{equation}

Consider now equation (\ref{L lambda ' M}). Fix $\lambda^{'},$ take
$M=M_n,$ and integrate by parts over $C.$ It follows, from
(\ref{cA}), (\ref{cq}) and the assumptions
$\nabla.A\tilde{e}\equiv0$ over $\Omega$ with $\nu.A\tilde{e}=0$ on
$\partial\Omega,$ that
$\displaystyle{\int_{C}M_n\,\nabla\cdot(A\nabla\psi^{\lambda^{'},M_n})}=0,$ $\displaystyle{\int_{C}\,-2\lambda{'}\sqrt{M_n}\tilde{e}\cdot
A\nabla\psi^{\lambda^{'},M_n}=0},$ and $\ds{\int_{C}q\cdot\nabla\psi^{\lambda^{'},M_n}=0.}$ Hence,
\begin{equation}\label{integrate L lamda ' by parts}
-\,\displaystyle{\frac{\lambda^{'}}{M_{n}^{^{\frac{1}{2}-\gamma}}}\int_{C}q\cdot\tilde{e}\,\psi^{\lambda^{'},M_{n}}}+\displaystyle{{\lambda^{'}}^{2}\int_{C}\tilde{e}\cdot
A\tilde{e}\,\psi^{\lambda^{'},M_n}+\int_{C}\zeta\,\psi^{\lambda^{'},M_n}=\mu(\lambda^{\,'},M_n)\int_{C}\psi^{\lambda^{'},M_n}}.
\end{equation}
Meanwhile, the functions $\psi^{\lambda^{'},M_n}$ converge to the
constant function $\psi^{\lambda^{'},\infty}$ in $L^{2}(C)$ strong;
and hence, in $L^{1}(C)$ strong ( $C$ is bounded, so $L^{2}(C)$ is
embedded in $L^{1}(C)$). Let $M_n\rightarrow +\infty$ in
(\ref{integrate L lamda ' by parts}):

In case $\gamma=\displaystyle{{1}/{2}},$ one has
$$\displaystyle{\frac{\lambda^{'}}{M_{n}^{^{\frac{1}{2}-\gamma}}}\int_{C}q\cdot\tilde{e}\,\psi^{\lambda^{'},M_{n}}=\lambda^{'}\int_{C}q\cdot\tilde{e}\,\psi^{\lambda^{'},M_{n}}\rightarrow\lambda^{'}\psi^{\lambda^{'},\infty}\displaystyle{\int_{C}q\cdot\tilde{e}=0}},$$ as
$\displaystyle{n\rightarrow+\infty\;}$ (from (\ref{cq})). Also, in
the case $0\leq\gamma<{1}/{2},\;$ one trivially has
$$\displaystyle{\frac{\lambda^{'}}{M_{n}^{^{\frac{1}{2}-\gamma}}}\int_{C}q\cdot\tilde{e}\,\psi^{\lambda^{'},M_{n}}\rightarrow0}\;\hbox{ as }\; n\rightarrow+\infty.$$
Moreover, $\tilde{e} A\tilde{e}$ and
$\zeta$ are in $L^{\infty}(C).$ Thus, as $M_n \rightarrow+\infty$ in
(\ref{integrate L lamda ' by parts}), we get
$$\displaystyle{{\lambda^{'}}^{2}\psi^{\lambda^{'},\infty}\int_{C}\tilde{e}A\tilde{e}\,+\psi^{\lambda^{'},\infty}\int_{C}\zeta\,=\displaystyle{l^{\,\lambda^{'},(M_n)}}\;\psi^{\lambda^{'},\infty}|C|}.$$

One concludes that
\begin{equation}\label{mu(lambda', inftty)}
   \displaystyle{
   \forall\lambda^{'}>0,\quad\frac{\displaystyle{l^{\,\lambda^{'},(M_n)}}}{\lambda^{'}}=\lambda^{\,'}\mi_{\!\!\!\!_C}\tilde{e} A\tilde{e}\,+\,\frac{\displaystyle{\mi_{\!\!\!\!_C}\zeta}}{\lambda^{\,'}}=\lambda^{'}m_0+\frac{m}{\lambda^{'}}}.
\end{equation}
Whence for
$\displaystyle{\lambda^{'}=\lambda^{'}_{0}=\sqrt{\frac{m}{m_0}}},$
one gets
$\displaystyle{\frac{\displaystyle{l^{\,\lambda_{0}^{'},(M_n)}}}{\lambda^{\,'}_{0}}=2\sqrt{m_{0}m}}.$

On the other hand, for all $M_n,$
\begin{equation}
    \displaystyle{\quad\frac{c_{\Omega,M_{n}A,\displaystyle{M_{n}^{\gamma}\,q},f}^*(e)}{\sqrt{M_n}}}=\displaystyle{\inf_{\lambda^{'}>0}\frac{\mu(\lambda^{\,'},M_n)}{\lambda^{\,'}}}
   \leq
   \displaystyle{\frac{\mu(\lambda^{\,'}_{0},M_n)}{\lambda^{\,'}_{0}}}.
\end{equation}
Passing $M_n \rightarrow+\infty,$ one gets
$\displaystyle{\limsup_{M_n\rightarrow+\infty}\;\frac{c_{\Omega,M_{n}A,\displaystyle{M_{n}^{\gamma}\,q},f}^*(e)}{\sqrt{M_n}}\leq\frac{\displaystyle{l^{\,\lambda_{0}^{'},(M_n)}}}{\lambda^{\,'}_{0}}=2\sqrt{m_0m}},$
and this holds for all sequences $\{M_n\}_{n}$ converging to
$+\infty.$ Thus,
\begin{equation}\label{limsup as M tends to infty}
\displaystyle{\limsup_{M\rightarrow+\infty}\;\frac{c_{\Omega,MA,\displaystyle{M^{\gamma}\,q},f}^*(e)}{\sqrt{M}}\leq\;2\sqrt{\mi_{\!\!\!\!_C}\tilde{e}A\tilde{e}(x,y)\,dxdy}\;\sqrt{\mi_{\!\!\!\!_C}\zeta(x,y)\,dxdy}}.
\end{equation}

Having (\ref{liminf as M tends }) together with (\ref{limsup as M tends to
infty}), the proof of Theorem \ref{lim as M} is complete.\hfill$\Box$

\paragraph{Proof of Theorem \ref{lim as per L tend to 0}.} We will consider the change of
variables similar to that made in the proof of Theorem \ref{lim as
per L tend to infty}:
$$v(t,x,y)=u(t,Lx,Ly),\quad(t,x,y)\in\mathbb{R}\times\mathbb{R}^{N}.$$

 After the same calculations done there, one gets
 that  $u$ satisfies (\ref{KPP L with u 2}) if and only if $v$
satisfies
\begin{eqnarray}\label{KPP L with v 2}
v_t(t,x,y)=\displaystyle{\frac{1}{L^{^{2}}}\nabla\cdot(A(x,y)\nabla
    v)(t,x,y)+ \frac{1}{L}\,q\cdot\nabla v(t,x,y) +f(x,y,v)\;\hbox{over}\;\mathbb{R}\times\mathbb{R}^{N}}.
\end{eqnarray}

Consequently,
\begin{equation}\label{beta * in terms of c* 2}
\displaystyle{  \forall\,L>0,\;c_{\mathbb{R}^{N},
\displaystyle{A_{_{L}},q_{_{L}},f_{_{L}}}}^*(e)\,=\,L\,c_{\mathbb{R}^{N},\frac{1}{L^{^{2}}}
A,\frac{1}{L}q,f}^*(e)}.
\end{equation}

On the other hand, the coefficients and the domain of problem
(\ref{KPP L with v 2}) satisfy all the assumptions of Theorem
\ref{lim as M}. Taking $\displaystyle{M={1}/{L^{^{2}}}}$  and
$\displaystyle{\gamma={1}/{2}},$ then (\ref{KPP L with v 2}) can be
rewritten as
\begin{equation*}
v_t(t,x,y)=\displaystyle{M\,\nabla\cdot(A(x,y)\nabla
    v)(t,x,y)+ M^{\frac{1}{2}}\,q\cdot\nabla v(t,x,y) +f(x,y,v)\;\hbox{over}\;\mathbb{R}\times\mathbb{R}^{N}}.
\end{equation*}
In this situation, the periodicity cell of the whole space
$\mathbb{R}^{N}$ is $\displaystyle{C=[0,1]\times\cdots\times[0,1]}.$

It follows, from Theorem \ref{lim as M}, that
\begin{eqnarray}\label{lim beta 2*}
\begin{array}{ll}
\displaystyle{\lim_{L\rightarrow0^+}\frac{c_{\mathbb{R}^{N},\frac{1}{L^{^{2}}}
A,\frac{1}{L}q,f}^*(e)}{\displaystyle{\sqrt{\frac{1}{L^{^{2}}}}}}}\,
&=\displaystyle{\lim_{M\rightarrow+\infty}\frac{c_{\mathbb{R}^{N},M\,
A,M^{\frac{1}{2}}\,q,f}^*(e)}{\displaystyle{\sqrt{M}}}}\\
&=\displaystyle{\,2\,\sqrt{\mi_{\!\!\!\!_C}\tilde{e}A\tilde{e}(x,y)dx\,dy}\,\sqrt{\mi_{\!\!\!\!_C}\zeta(x,y)dx\,dy}.}
\end{array}
\end{eqnarray}

Having (\ref{beta * in terms of c* 2}) together with (\ref{lim
beta 2*}), the proof of Theorem \ref{lim as per L tend to 0} is complete.
\hfill$\Box$

\subsection{Proofs of Theorems \ref{lim as B tend to infty} and \ref{lim as B tend 0} }

\textbf{Proof of Theorem \ref{lim as B tend to infty}.} The main
ideas of this proof are similar to those in the demonstration of
Theorem \ref{limit as eps}. Applying the variational formula
(\ref{var}) of the minimal speed, one gets
\begin{equation}\label{var form with B}
  c_{\Omega,
A,0,\,Bf}^*(e)\,=\,\displaystyle{\min_{\lambda\,>\,0}\frac{k_{\Omega,e,\,
A,\,0,\,B\zeta}(\lambda)}{\lambda}},
\end{equation}

where $k_{\Omega,e,A,0,B\zeta}(\lambda)$ is the first eigenvalue
(for each $\lambda,\,B\,>\,0$) of the eigenvalue problem:
\begin{eqnarray}\label{eigeneq11}
\left\{
    \begin{array}{ll}
     \displaystyle{L_{\Omega,e,\, A,\,0,B\zeta,\lambda}\;\psi(x,y)}&=\displaystyle{k_{_{\Omega,e,\, A,\,0,B\zeta}}(\lambda)\;\psi(x,y)\;\hbox{over}\;\mathbb{R}\times\omega;}
     \vspace{4 pt}\\
  \nu\cdot A\nabla\psi&=0 \quad
  \hbox{on}\;\mathbb{R}\,\times\,\partial\omega,
    \end{array}
  \right.
\end{eqnarray}
 and \begin{eqnarray*}
 L_{\Omega,e,\, A,\,0,B\,\zeta,\lambda}\psi(x,y)\,&=&\,\nabla\cdot \left(A(y)\nabla\psi(x,y)\right)\,-2\,\lambda
Ae\cdot\nabla\psi(x,y)\,+\vspace{4 pt}\\
&&\left[\lambda^{2}eA(y)e\,-\lambda\,\nabla\cdot(A(y)e)+\,B\zeta(y)\right]\psi(x,y),
 \end{eqnarray*}
for each $(x,y)\in\mathbb{R}\times\omega.$

We recall that for all $\lambda>0,$ and for all $B>0,$ we have
$k_{\Omega,e,\, A,\,0,\displaystyle{B\zeta}}(\lambda)\,>\,0.$ Also,
the first eigenfunction of (\ref{eigeneq11}) is positive over
$\overline{\Omega}=\,\mathbb{R}\times\overline{\omega},$ and it is
unique up to multiplication by a non zero constant.

Moreover, whether in (\ref{alt 1}) or (\ref{alt 2})
and due to (\ref{nonlinearity 2}), one concludes that
the coefficients
 in $ \displaystyle{L_{\Omega,e,\,A,\,0,\displaystyle{B\zeta},\lambda}}$
are independent of $x.$ Hence, the first eigenfunction of
(\ref{eigeneq11}) is independent of $x$ and the eigenvalue
problem (\ref{eigeneq11}) is reduced to
\begin{eqnarray}\label{eigeneq21}
\left\{
    \begin{array}{rl}
   \displaystyle{ L_{\Omega,e,\,
   A,0,\displaystyle{B\zeta},\lambda}\phi}&:=\,\nabla\cdot\left(A(y)\nabla\phi(y)\right)\,+\,\left[\lambda^{2}eA(y)e\,+\,B\zeta(y)\right]\phi(y)\vspace{3 pt}\\
&=\,\displaystyle{k_{\Omega,e,\, A,\,0,
B\zeta}(\lambda)}\,\phi\quad\hbox{over}\;\omega;
   \vspace{3pt}\\
  \nu(x,y)\cdot A(y)\nabla\phi(y)&=\left(0;\nu_{\omega}(y)\right)\cdot A(y)\nabla\phi(y)=0 \quad \hbox{on}\;\mathbb{R}\times\partial\omega,
    \end{array}
  \right.
\end{eqnarray}
where $\phi\,=\,\phi(y)$ is positive over $\overline{\omega},$
$L-$periodic (since the domain $\omega$ and the coefficients of
$L_{\Omega,e,\,A,0,\displaystyle{B\zeta},\lambda}$ are
$L-$periodic), unique up to multiplication by a constant, and
belongs to $C\,^{2}(\overline{\omega}).$

In the case where $d\geq1,$ let $C\subseteq\mathbb{R}^{N-1}$ denote
the periodicity cell of $\omega.$ Otherwise, $d=0$ and one takes
$C=\omega.$ In both cases, $C$ is bounded. Multiplying the first
line of (\ref{eigeneq21}) by $\phi,$ and integrating by parts over
$C,$ one gets
\begin{equation}\label{_k11}
    -\,\displaystyle{k_{\Omega,e,\, A,\,0,\displaystyle{B\zeta}}(\lambda)}\,=\,\displaystyle{\frac{\displaystyle{\int_{C}\,\nabla\phi\cdot A(y)\nabla\phi\,dy}-\displaystyle{\int_{C}\left[\lambda^{2}eA(y)e\,+\,B\,\zeta(y)\right]\,\phi^{2}(y)\,dy}}{\displaystyle{\int_{C}{\phi^2}(y)\,dy}}.}
\end{equation}

One also notes that, in this present setting, the operator
$L_{\Omega,e,\, A,0,\displaystyle{B\zeta},\lambda}$ is self-adjoint
and its coefficients are $(L_1,\ldots,L_d)-$periodic with respect
$(y_1,\ldots,y_d).$ Consequently, $-\,\displaystyle{k_{\Omega,e,\,
A,\,0,\displaystyle{B\zeta}}(\lambda)}$ has the following
variational characterization:
\begin{equation}\label{k21}
    -\,\displaystyle{k_{\Omega,e,\, A,\,0,\displaystyle{B\zeta}}(\lambda)}=\displaystyle{\min_{\varphi\in
    H^{1}(C)\setminus\{0\}}\frac{\displaystyle{\int_{C}\,\nabla\varphi\cdot A(y)\nabla\varphi\,dy}-\displaystyle{\int_{C}\left[\lambda^{2}eA(y)e\,+\,B\,\zeta(y)\right]\,\varphi^{2}(y)\,dy}}{\displaystyle{\int_{C}{\varphi^2}(y)\,dy}}.}
\end{equation}

In what follows, we will assume that (\ref{alt 1}) is the
alternative that holds. That is, $eAe=\alpha$ is constant. The proof
can be imitated easily whenever we assume that (\ref{alt 2}) holds.

The function $y\mapsto\zeta(y)$ is continuous and
$(L_1,\ldots,L_d)-$periodic over $\overline{\omega},$ whose
periodicity cell $C$ is a bounded subset of $\mathbb{R}^{N-1}$
(whether $d=0$ or $d\geq1$). Let
$y_0\in\overline{C}\subseteq\,\overline{\omega}$ such that
$\displaystyle{\max_{y\in\overline{w}}\,\zeta(y)=\zeta(y_0)}$
(trivially, this also holds when $\zeta$ is constant). Consequently,
we have
$$\forall\;\varphi\in \,H^{1}(C)\setminus\{0\},\;\displaystyle{\frac{\displaystyle{\int_{C}\,\nabla\varphi\cdot \,A\nabla\varphi}-\displaystyle{\int_{C}(\alpha\lambda^{2}\,+\,B\,\zeta(y))\varphi^{2}}}{\displaystyle{\int_{C}{\varphi^2}(y)\,dy}}\,\geq\,-\left[\alpha\lambda^{2}\,+\,B\,\zeta(y_0)\right].}$$
This yields that
\begin{equation}\label{k< 1}
\forall\,B>0,\,\forall\,\lambda>0,\,-\,\displaystyle{k_{\Omega,e,\,
A,\,0,\displaystyle{B\,\zeta}}(\lambda)}\,\geq\,-\left[\alpha\lambda^{2}\,+\,B\,\zeta(y_0)\right].
\end{equation}
Consequently,
\begin{equation}\label{m1}
    \displaystyle{\forall\,B>0,\,\forall\,\lambda>0,\,\frac{\displaystyle{k_{\Omega,e,\, A,\,0,\displaystyle{B\zeta}}(\lambda)}}{\lambda}\,\leq\,\lambda\,\alpha\,+\frac{B\,\zeta(y_0)}{\lambda}}.
\end{equation}

However, the function
$\displaystyle{{\lambda}\,\mapsto\,\lambda\alpha\,+\left({B\,\zeta(y_0)}/{\lambda}\right)}$
attains its minimum, over $\mathbb{R}^{+},$ at
$\lambda(B)\,=\,\displaystyle{\sqrt{\frac{B\,\zeta(y_0)}{\alpha}}}$.
This minimum is equal to $2\sqrt{B\zeta(y_0)}\times\sqrt{\alpha}.$

 From (\ref{m1}), we conclude that:
$\displaystyle{\frac{\displaystyle{k_{\Omega,e,\,
A,\,0,\displaystyle{B\zeta}}(\lambda(B))}}{\lambda(B)}\,\leq\,2\sqrt{B\,\alpha}\sqrt{\zeta(y_0)}}.$

Finally, (\ref{var}) implies that
$$\displaystyle{c_{\Omega,
A,0,\displaystyle{Bf}}^*(e)\,=\,\displaystyle{\min_{\lambda\,>\,0}\frac{k_{\Omega,e,\,
A,\,0,\displaystyle{B\zeta}}(\lambda)}{\lambda}}\,\leq\,2\sqrt{B\,\alpha}\sqrt{\zeta(y_0),}}$$
or equivalently
\begin{equation}\label{c(B) over B < }
   \forall\,B\,>0,\, \displaystyle{\frac{c_{\Omega,
A,0,Bf}^*(e)}{\sqrt{B}}\,\leq\,2\sqrt{\alpha}\sqrt{\zeta(y_0).}}
\end{equation}

We pass now to prove the other sense of the inequality for
$\displaystyle{\liminf_{B\rightarrow+\infty}\displaystyle{\frac{c_{\Omega,
A,0,Bf}^*(e)}{\sqrt{B}}}}.$ We will consider formula (\ref{k2}), and
then organize a suitable function $\psi$ which leads us to a lower
bound of
$\displaystyle{\liminf_{B\rightarrow+\infty}\displaystyle{\frac{c_{\Omega,
A,0,Bf}^*(e)}{\sqrt{B}}}}.$

We have $\zeta(y_0)\,>\,0.$ Let $\delta$ be such that
$0<\delta<\zeta(y_{0}).$ Thus
$\displaystyle{0<\zeta(y_{0})-\delta<\max_{\overline{\omega}}\,\zeta(y)}.$
The continuity of $\zeta,\;$ over
$\overline{C}\subseteq\overline{\omega},$ yields that there exists
an open and bounded set $U\subset\overline{C}$ such that
\begin{equation}\label{setU1}
    \zeta(y_{0})-\delta\leq\zeta(y),\;\forall\,y\in\,\overline{U}.
\end{equation}

Designate by $\psi,$ a function in $\mathcal{D}(C)$ (a
$C^{\infty}(C)$ function whose support is compact), with
supp$\,\psi\,\subseteq\,\overline{U},$ and $\displaystyle{\int_{U}\psi^{2}=1.}$ One
will have,
\begin{eqnarray*}
\begin{array}{ll}
   \forall\lambda>0,\,\forall\,B\,>0,&\vspace{3 pt}\\
    -\, k_{\Omega,e,\,
A,\,0,\displaystyle{B\,\zeta}}(\lambda)&\,\leq\,\displaystyle{\int_{U}\,\nabla\psi\cdot A(y)\nabla\psi\,dy-\displaystyle{\int_{U}\left[\lambda^{2}eA(y)e\,+\,B\,\zeta(y)\right]\,\psi^{2}(y)\,dy}}\vspace{3 pt}\\
                            \nonumber  &\leq\,\displaystyle{\int_{U}\,\nabla\psi\cdot A(y)\nabla\psi\,dy-\left[\lambda^{2}\alpha\,+\,B\,(\zeta(y_0)-\delta)\right]}\;\hbox{(by (\ref{setU1}))}\vspace{3 pt}\\
&\leq\displaystyle{\int_{U}\,\alpha_2|\nabla\psi|^2\,-\left[\lambda^{2}\alpha\,+\,B\,\left(\zeta(y_0)-\delta\right)\right]}\;\hbox{by}\;(\ref{cA
Y}),
\end{array}
 \end{eqnarray*}
or equivalently
\begin{equation}\label{k() over lambda > 1}
\displaystyle{\frac{k_{\Omega,e,\,
A},\,0,\displaystyle{B\,\zeta}(\lambda)}{\lambda}}\,\geq\,\lambda\alpha\,+\,\displaystyle{\frac{B}{\lambda}}\;\rho(B),
\end{equation}
where
$\rho(B)=\zeta(y_0)-\delta-\displaystyle{\frac{1}{B}}\,\displaystyle{\int_{U}\alpha_2|\nabla\psi|^{2}}.$
Choosing $B$ large enough, we get $\rho(B)\,>0$ (this is possible
since $\zeta(y_0)-\delta>0\;$ and also
$\;\displaystyle{\int_{U}\alpha_2|\nabla\psi|^{2}}>0$).
 The map
$\lambda\mapsto\lambda\alpha\,+\,\displaystyle{\frac{B}{\lambda}}\;\rho(B)$
attains its minimum, over $\mathbb{R}^+,$ at
$\lambda(\varepsilon)=\displaystyle{\sqrt{\frac{B\,\rho(B)}{\alpha}}}.$
This minimum is equal to
$2\displaystyle{\sqrt{B\,\alpha}\,\sqrt{\rho(B)}}.$

Now, referring to formula (\ref{k() over lambda > 1}), one gets:
 $$\;\hbox{for $B$ large enough},\;\displaystyle{\frac{k_{\Omega,e,\,
A,\,0,\displaystyle{B\zeta}(\lambda)}\;}{\lambda}\,\geq\,2\displaystyle{\sqrt{B\,\alpha}\,\sqrt{\rho(B)}}\quad\hbox{for
all $\lambda>0$}.}$$ Together with (\ref{var}), we conclude that
\begin{equation}
\, \hbox{for $B$ large enough,}\quad\displaystyle{\frac{ c_{\Omega,
A,0,\displaystyle{Bf}}^*(e)}{\sqrt{B}}} \;\geq\;
2\sqrt{\rho(B)}\sqrt{\alpha}.
\end{equation}
 Consequently, \begin{eqnarray}
                \nonumber \displaystyle{\liminf_{B\rightarrow+\infty}\frac{ c_{\Omega,
A,0,\displaystyle{B f}}^*(e)\,}{\sqrt{B}}}\,&\geq&\,\displaystyle{\liminf_{B\rightarrow+\infty}2\sqrt{\rho(B)}\sqrt{\alpha}\,} \\
                \nonumber &=&2\sqrt{\zeta(y_0)-\delta}\sqrt{\alpha}\quad\hbox{(since $\psi$ is independent of
$B$)},
               \end{eqnarray}
and this holds for all $0<\delta<\zeta(y_{0}).$ Therefore, one can
conclude that
\begin{equation}\label{liminf c(B) over rad B}
    \displaystyle{\liminf_{B\rightarrow+\infty}\frac{ c_{\Omega,
A,0,Bf}^*(e)\,}{\sqrt{B}}}\,\geq\,2\sqrt{\alpha}\sqrt{\zeta(y_0).}
\end{equation}

Finally, the inequalities (\ref{c(B) over B < }) and (\ref{liminf
c(B) over rad B}) imply that
$\;\displaystyle{\lim_{B\rightarrow+\infty}\frac{ c_{\Omega,
A,0,Bf}^*(e)\,}{\sqrt{B}}}\;$ exists, and it is equal to
$\displaystyle{2\sqrt{\alpha}\sqrt{\zeta(y_0)}\,=\,2\sqrt{\max_{\overline{\omega}}eA(y)e}\sqrt{\max_{\overline{\omega}}\zeta(y)}}.$

The above proof was done while assuming that the alternative
(\ref{alt 1}) holds. The same ideas of this proof can be easily
applied in the case where alternative (\ref{alt 2}) holds. In
(\ref{alt 2}), we have $\zeta$ is constant; however, $eAe$ is not in
general. Meanwhile the converse is true in the case (\ref{alt 1}).
The little difference is that, in the case of (\ref{alt 2}), we
chsose the subset $U$ (of the proof done above) around the point
$y_0$ where $eAe$ attains its maximum and then we continue by the
same way
used above. \hfill$\Box$\\

\textbf{Proof of Theorem \ref{lim as B tend 0}.}
 According to Theorem \ref{varthm}, and since
$\nu\cdot A\tilde{e}=0$ on $\partial\Omega,$ the minimal
 speeds $c_{\Omega,A,\displaystyle{B^{\gamma}\,q},\, Bf}^*(e)$ are given by:
 \begin{equation*}
   \displaystyle{ \forall
   B>0,\;c_{\Omega,A,\displaystyle{B^{\gamma}\,q},\, Bf}^*(e)=\min_{\lambda>0}\frac{k_{\Omega,e,\,
A,\,\displaystyle{B^{\gamma}\,q},\,B\zeta}(\lambda)}{\lambda}},
 \end{equation*}
where $k_{\Omega,e,\,
A,\,\displaystyle{B^{\gamma}\,q},\,B\zeta}(\lambda)$ and
$\psi^{\lambda,B}$ denote the unique eigenvalue and the positive
$L$-periodic eigenfunction of the problem
\begin{eqnarray*}\label{L lambda B}
\displaystyle{\nabla\cdot(A\nabla\psi^{\lambda,B})-2\lambda\tilde{e}\cdot
A\nabla\psi^{\lambda,B}+\displaystyle{B^{\gamma}q\cdot
\nabla\psi^{\lambda,B}}+\left[\lambda^{2}\,\tilde{e}A\tilde{e}-\lambda
B^{\gamma}q\cdot\tilde{e}+\,B\,\zeta\right]\psi^{\lambda,B}}\\
=k_{\Omega,e,\,
A,\,\displaystyle{B^{\gamma}\,q},\,B\zeta}(\lambda)\;\psi^{\lambda,B}\;\;\hbox{in
$\Omega,$ with $\nu\cdot A\nabla\psi=\nu\cdot
A\nabla\psi^{\lambda,B}=0$ on $\,\partial\Omega.$ }
\end{eqnarray*}

For each $\lambda>0$ and $B>0,$ let
$\lambda^{\,'}=\displaystyle{{\lambda}/{\sqrt{B}}},$ and let
$\;k_{\Omega,e,\,
A,\,\displaystyle{B^{\gamma}\,q},\,B\zeta}(\lambda)=\mu(\lambda^{'},B).\;$
Consequently,
\begin{equation}\label{var with mu 2}
 \displaystyle{ \forall
   B>0,\;\frac{c_{\Omega,A,\displaystyle{B^{\gamma}\,q},\,Bf}^*(e)}{\sqrt{B}}=\min_{\lambda^{\,'}>0}\frac{\mu(\lambda^{\,'},B)}{\lambda^{\,'}\,B}},
\end{equation}
where $\mu(\lambda^{\,'},B)$ and $\psi^{\lambda^{\,'},B}$ are the
first eigenvalue and the unique, positive $L-$periodic (with respect
to $x$) eigenfunction of
\begin{eqnarray}\label{L lambda ' B}
\begin{array}{c}
\displaystyle{\nabla\cdot(A\nabla\psi^{\lambda^{'},B})-2\lambda{'}\sqrt{B}\tilde{e}\cdot
A\nabla\psi^{\lambda^{'},B}+\displaystyle{B^{\gamma}q\cdot\nabla\psi^{\lambda^{'},B}}}\vspace{4 pt}\\
+\left[{\lambda^{'}}^{2}B\,\tilde{e}A\tilde{e}-\displaystyle{\lambda^{'}
B^{^{\gamma+\frac{1}{2}}}}\,q\cdot \tilde{e}+\,B\zeta\right]\psi^{\lambda^{'},B}\;=\;\mu(\lambda^{'},B)\psi^{\lambda^{'},B}\hbox{ in } \Omega,
\end{array}
\end{eqnarray}
with $\nu\cdot A\nabla\psi^{\lambda^{\,'},B}=0$ on
$\partial\Omega.$

Owing to the uniqueness, up to multiplication by positive constants,
of the first eigenfunction of (\ref{L lambda ' B}), one may assume
that:
 \begin{equation}\label{norm L2 of psi}
    \forall \,\lambda^{'}>0,\;\forall
    \,B>0,\;||\psi^{\lambda^{'},B}||_{L^{2}(C)}\,=1.
 \end{equation}

Moreover, for each
$\displaystyle{B>0,\;\min_{\lambda^{\,'}>0}\frac{\mu(\lambda^{\,'},B)}{\lambda^{\,'}\,B}}$
is attained at $\lambda^{\,'}_{B}>0.$ Thus,
\begin{equation}\label{var with mu and lambda' 2}
 \displaystyle{ \forall
   B>0,\;\frac{c_{\Omega,A,\displaystyle{B^{\gamma}\,q},\,Bf}^*(e)}{\sqrt{B}}=\min_{\lambda^{\,'}>0}\frac{\mu(\lambda^{\,'},B)}{\lambda^{\,'}\,B}=\frac{\mu(\lambda^{\,'}_{B},B)}{B\,\lambda^{\,'}_{B}}}.
\end{equation}

Having the above characterization, one can now imitate the steps 2 and 3 in the proof of Theorem \ref{lim as M} to prove that
$$\displaystyle{\liminf_{B\rightarrow
0^+}{c_{\Omega,A,\displaystyle{B^{\gamma}\,q},\,Bf}^*(e)}/{\sqrt{B}}}$$
(resp.
$\displaystyle{\limsup_{B\rightarrow0^+}{c_{\Omega,A,\displaystyle{B^{\gamma}\,q},\,Bf}^*(e)}/{\sqrt{B}}}$
) is greater than (resp. less than)
$$\displaystyle{2\sqrt{\mi_{\!\!\!\!_C}\tilde{e}A\tilde{e}(x,y)dx\,dy}\,\sqrt{\mi_{\!\!\!\!_C}\zeta(x,y)dx\,dy}};$$
and hence, complete the proof of Theorem \ref{lim as B tend 0}.\hfill$\Box$\vskip 0.3cm

\subsection{Proofs of Theorems  \ref{variation of min speed with respect diffusion}, \ref{variation of the min speed with respect to period}, and \ref{variation of c*(B)over rad
B}}

\vskip0.3cm

\textbf{Proof of Theorem \ref{variation of min speed with respect
diffusion}.} Referring to Theorem \ref{varthm}, it follows that for
each $\beta>0,$ we have:
$$\displaystyle{\frac{c_{\Omega,\beta A,\sqrt{\beta}\,q,f}^*(e)}{\sqrt{\beta}}=\min_{\lambda>0}\frac{k_{\Omega,e,\,\beta A,\,\sqrt{\beta}\,q,\zeta}(\lambda)}{\lambda\sqrt{\beta}}},$$
where $k_{\Omega,e,\,\beta A,\,\sqrt{\beta}\,q,\zeta}(\lambda)$ is
the first eigenvalue of the problem
\begin{eqnarray}
    \left\{
    \begin{array}{rl}
      \displaystyle{L_{\Omega,e,\,\beta A,\,\sqrt{\beta}\,q,\zeta,\lambda}}\psi(x,y)&=\;\displaystyle{k_{\Omega,e,\,\beta A,\,\sqrt{\beta}\,q,\zeta}(\lambda)}\,\psi(x,y)\hbox{ over }\mathbb{R}\times\omega;\vspace{3 pt}\\
  \nu.A\nabla\psi&=0
  \hbox{ on }\mathbb{R}\,\times\,\partial\omega,
    \end{array}
  \right.
  \end{eqnarray}
 where \begin{eqnarray*}
\displaystyle{L_{\Omega,e,\beta
A,\sqrt{\beta}\,q,\zeta,\lambda}}\,\psi&=&\beta\nabla\cdot\left(A(y)\nabla\psi\right)-2\beta\lambda\,\alpha(y)\,\partial_{x}\psi+\sqrt{\beta}\,q_{_{1}}(y)\partial_{x}\psi\\
&&+\left[\beta\,\lambda^{2}eA(y)e\,-\lambda
\sqrt{\beta}\,q_{_{1}}(y)+\zeta(y)\right]\psi\;\hbox {over}\;
\mathbb{R}\times\omega.
\end{eqnarray*}

The boundary condition follows so from the facts that
$\Omega=\mathbb{R}\times \omega,\;e=(1,0,\ldots,0)$ and that
$A(y)e=\alpha(y)e$ over $\omega.$ These yield that $\nu\cdot Ae=0$
over $\partial\Omega$ and $\nabla\cdot Ae=0.$ Moreover, for each
$(x,y)\in\partial\Omega,$ we have $\nu(x,y)=(0;\nu_{\omega}(y)),$
where $\nu_{\omega}(y)$ is the outward unit normal on
$\partial\omega$ at $y.$

On the other hand, the function $\psi$ is positive ,
$(L_1,\ldots,L_d)-$periodic with respect to $y,$ and unique up to
multiplication by non-zero constants. Meanwhile, the coefficients
$A,\;q$ and  $\zeta$ are independent of $x.$ Thus the eigenfunction
$\psi$ will be independent of $x$ and our eigenvalue problem is
reduced to
\begin{eqnarray}\label{reduced eigen value problem 1}
    \left\{
    \begin{array}{c}
    \beta\nabla\cdot\left(A(y)\nabla\psi(y)\right)+\left[\beta\,\lambda^{2}eA(y)e\,-\lambda \sqrt{\beta}q_{_{1}}(y)+\zeta(y)\right]\psi(y)\vspace{3 pt}\\
    =\displaystyle{k_{\Omega,e,\,\beta A,\,\sqrt{\beta}q,\zeta}(\lambda)}\psi(y)\hbox{ for all } y\in\omega,
    \vspace{3 pt}\\
  \nu(x,y)\cdot A(y)\nabla\psi(y)\;=\;\left(0;\nu_{\omega}(y)\right)\cdot A(y)\nabla\psi(y)\;=\;0 \; \hbox{ on }\;\mathbb{R}\times\partial\omega.
    \end{array}
  \right.
\end{eqnarray}

For each $\lambda>0$ and $\beta>0,$ let
$\lambda^{'}=\lambda\sqrt{\beta},$ and let
$\displaystyle{k_{\Omega,e,\,\beta
A,\,\sqrt{\beta}\,q,\zeta}(\lambda)}=\mu(\lambda^{'},\beta).$
Since for each $\beta>0,
\;\displaystyle{\min_{\lambda>0}\frac{k_{\Omega,e,\,\beta
A,\,\sqrt{\beta}\,q,\zeta}(\lambda)}{\lambda}}$ is attained at
$\lambda(\beta),$ it follows that
\begin{equation}\label{var form with lambda'}
   \forall\,\beta>0,\quad \displaystyle{\frac{c_{\Omega,\beta A,\sqrt{\beta}\,q,f}^*(e)}{\sqrt{\beta}}=\min_{{\lambda}^{'}>0}\frac{\mu(\lambda^{'},\beta)}{\lambda^{'}}},
\end{equation}
where $\mu(\lambda^{'},\beta)$ is the first eigenvalue of the
problem:
\begin{eqnarray}\label{ reduced eigenvalue problem 2}
\left\{
  \begin{array}{rl}
   L_{\lambda^{'}}^{\beta}\psi&=\beta\nabla\cdot\left(A(y)\nabla\psi\right)+\left[{\lambda^{'}}^{2}eA(y)e-\lambda^{'}q_{_{1}}(y)+\zeta(y)\right]\psi=\mu(\lambda^{'},\beta)\psi\hbox{ in }\omega,\vspace{5 pt} \\
  \nu\cdot A\nabla\psi&=0 \hbox{ on }\partial\omega.
  \end{array}
\right.
\end{eqnarray}

The elliptic operator $L_{\lambda^{'}}^{\beta}$ in (\ref{ reduced
eigenvalue problem 2}) is self-adjoint. Consequently, the first
eigenvalue $\mu(\lambda^{'},\beta)$ has the following
characterization:\footnote{To have an idea, multiply (\ref{
reduced eigenvalue problem 2}) by the positive,
$(L_1,\ldots,L_d)-$periodic function $\psi$ and integrate by parts
over the periodicity cell $C$ of the the domain $\omega.$}
\begin{eqnarray}\label{-mu(lambda',beta)=min}
    \begin{array}{ll}
\forall \lambda^{'}>0,\,\forall \beta>0,\quad-\mu(\lambda^{'},\beta)=
\vspace{4 pt}\\
\displaystyle{\min_{\varphi\in
    H^{1}(C)\setminus\{0\}}\frac{\displaystyle{\beta\int_{C}\nabla\varphi\cdot A(y)\nabla\varphi dy}+\lambda^{'}\int_{C}q_{_{1}}(y)\varphi^{2}-\displaystyle{{\int_{C}\left[{\lambda^{'}}^{2}eA(y)e+\zeta(y)\right]\varphi^{2}(y)dy}}}{\displaystyle{\int_{C}\varphi^{2}(y)dy}}}\vspace{3 pt}\\
\quad\quad\quad\quad=\displaystyle{\min_{\varphi\in
    H^{1}(C)\setminus\{0\}}\,R(\lambda^{'},\beta,\varphi).}
\end{array}
\end{eqnarray}

For each $\lambda^{'}$ and $\beta>0,\;\varphi\mapsto
R(\lambda^{'},\beta,\varphi)$ attains its minimum over
$H^{1}(C)\setminus\{0\}$ at $\psi^{\lambda^{'},\beta},$ the
eigenfunction of the problem (\ref{ reduced eigenvalue problem 2}).
On the other hand, $\beta\mapsto R(\lambda^{'},\beta,\varphi)$ is
increasing as an affine function in $\beta.$ Consequently, fixing
$\lambda^{'}>0\;$ and taking $\;\beta>\beta^{'}>0:$
\begin{eqnarray}
\begin{array}{ll}
-\mu(\lambda^{'},\beta)&=R(\lambda^{'},\beta,\psi^{\lambda^{'},\beta})>R(\lambda^{'},\beta^{'},\psi^{\lambda^{'},\beta})\\
&\geq \displaystyle{\min_{\varphi\in
    H^{1}(C)\setminus\{0\}}\,R(\lambda^{'},\beta^{'},\varphi)}=-\mu(\lambda^{'},\beta^{'}).
\end{array}
\end{eqnarray}

In other words, for all $\lambda^{'}>0,$ the function $\beta\mapsto
\mu(\lambda^{'},\beta)$ is decreasing. Concerning now the function
$\beta\mapsto \displaystyle{{c_{\Omega,\beta
A,\sqrt{\beta}\,q,f}^*(e)}/{\sqrt{\beta}}},$ one takes randomly $\beta>\beta^{'}>0,$ hence
\begin{eqnarray}
\begin{array}{rl}
\nonumber\displaystyle{\frac{c_{\Omega,\beta^{\,'}
A,\sqrt{\beta^{\,'}}\,q,f}^*(e)}{\sqrt{\beta^{'}}}}&=\displaystyle{\frac{\mu(\lambda^{'}(\beta^{'}),\beta^{'})}{\lambda^{'}(\beta^{'})}}
>\displaystyle{\frac{\mu(\lambda^{'}(\beta^{'}),\beta)}{\lambda^{'}(\beta^{'})}}\\
\nonumber&\geq\displaystyle{\min_{{\lambda}^{'}>0}\frac{\mu(\lambda^{'},\beta)}{\lambda^{'}}}=
\displaystyle{\frac{c_{\Omega,\beta
A,\sqrt{\beta}\,q,f}^*(e)}{\sqrt{\beta}},}
\end{array}
\end{eqnarray}
which means that the function $\beta\mapsto
\displaystyle{{c_{\Omega,\beta
A,\sqrt{\beta}\,q,f}^*(e)}/{\sqrt{\beta}}}$ is decreasing.

Finally, when $\beta\rightarrow+\infty,$ one can easily check that
the hypothesis of Theorem \ref{lim as M} are satisfied; hence, one
has the limit at $+\infty,$ and that completes the
proof of Theorem \ref{variation of min speed with respect diffusion}.\hfill$\Box$
\vskip0.3cm

\textbf{Proof of Theorem \ref{variation of the min speed with
respect to period}.} Consider the change of variables
$v(t,x,y)=u(t,Lx,Ly),$ for any $\quad(t,x,y)\in\mathbb{R}\times\mathbb{R}^{N}.$
One consequently has,
\begin{equation}\label{beta * in terms of c* 3}
\displaystyle{  \forall\,L>0,\;c_{\mathbb{R}^{N},
\displaystyle{A_{_{L}},q_{_{L}},f_{_{L}}}}^*(e)\,=\,L\,c_{\mathbb{R}^{N},\frac{1}{L^{^{2}}}
A,\frac{1}{L}q,f}^*(e)}.
\end{equation}

 Taking $\displaystyle{\beta={1}/{L^{^{2}}}},$ then
\begin{equation*}
v_t(t,x,y)=\displaystyle{\beta\,\nabla\cdot(A(y)\nabla
    v)(t,x,y)+ \sqrt{\beta}\,q_{_{1}}(y)\,\partial_{x}\,v(t,x,y) +f(x,y,v)\;\hbox{over}\;\mathbb{R}\times\mathbb{R}^{N}}.
\end{equation*}
 Owing to Theorem \ref{variation of min
speed with respect diffusion}, the function
$\beta\mapsto\displaystyle{{c_{\mathbb{R}^{N},\beta
A,\displaystyle{\sqrt{\beta}\,q},f}^*(e)}/{\sqrt{\beta}}}\;$ is
decreasing in $\beta>0.$ Besides, $L\mapsto1/L^{2}$ is decreasing in $L>0.$ Together with (\ref{beta * in terms of c* 3}), one obtains that
the function $\displaystyle{L\mapsto c_{\mathbb{R}^{N},
\displaystyle{A_{_{L}},q_{_{L}},f_{_{L}}}}^*(e)}\;$ is increasing in
$L>0$ which completes the proof of Theorem \ref{variation of the
min speed with respect to period}.\hfill$\Box$
\vskip0.3cm

\textbf{Proof of Theorem \ref{variation of c*(B)over rad B}.}
Referring to Theorem \ref{varthm}, it follows that for each $B>0,$
we have:
$$\displaystyle{\frac{c_{\Omega, A,\sqrt{B}\,q,Bf}^*(e)}{\sqrt{B}}=\min_{\lambda>0}\frac{k_{\Omega,e,\, A,\,\sqrt{B}\,q,B\zeta}(\lambda)}{\lambda\sqrt{B}}}.$$

Owing to the same justifications explained in the proof of Theorem
\ref{variation of min speed with respect diffusion}, $k_{\Omega,e, A,\sqrt{B}q,B\zeta}(\lambda)$ is the first
eigenvalue of the problem
\begin{eqnarray}\label{reduced eigen value problem 21}
    \left\{
    \begin{array}{c}
    \nabla\cdot\left(A(y)\nabla\psi(y)\right)+\left[\lambda^{2}e\cdot Ae-\lambda \sqrt{B}q_{_{1}}(y)+B\zeta(y)\right]\psi(y)=\displaystyle{k_{\Omega,e,A,\sqrt{B}q,B\zeta}(\lambda)}\psi\hbox{ in }\omega,
    \vspace{4 pt}\\
  \nu(x,y)\cdot A(y)\nabla\psi(y)=\left(0;\nu_{\omega}(y)\right)\cdot A(y)\nabla\psi(y)=0 \;\hbox{ on }\;\mathbb{R}\times\partial\omega.
    \end{array}
  \right.
\end{eqnarray}

For each $\lambda>0$ and $B>0,$ let $\lambda^{'}=\displaystyle{{\lambda}/{\sqrt{B}}}$ and
$\displaystyle{k_{\Omega,e,\,
A,\,\sqrt{B}\,q,B\zeta}(\lambda)}=\mu(\lambda^{'},B).$ The first
eigenvalue $\mu(\lambda^{'},B)$ has the following characterization:
\begin{eqnarray}\label{-mu(lambda',B) over lambda' B=min}
    \begin{array}{ll}
\forall \lambda^{'}>0,\,\forall
B>0,\;-\,\displaystyle{\frac{\mu(\lambda^{'},B)}{\lambda^{'}B}}=\\
\displaystyle{\min_{\underset{\displaystyle{||\varphi||_{L^{^{2}}(C)}=1}}{\displaystyle{\varphi\in
    H^{1}(C)\setminus\{0\};}}}\frac{\displaystyle{\int_{C}\nabla\varphi\cdot A(y)\nabla\varphi\,
    dy}}{\lambda^{'}B}+\displaystyle{\int_{C}q_{_{1}}\varphi^{2}}-\lambda^{'}\displaystyle{\int_{C}eAe\varphi^{2}}-\displaystyle{\frac{\displaystyle{\int_{C}\zeta(y)\varphi^{2}(y)\,dy}}{\lambda^{'}}}}\\
   \quad\quad\quad\quad \quad= \displaystyle{\min_{\underset{\displaystyle{||\varphi||_{L^{^{2}}(C)}=1}}{\displaystyle{\varphi\in
    H^{1}(C)\setminus\{0\}}}}R(\lambda^{'},B,\varphi)}.
\end{array}
\end{eqnarray}

On the other hand, $B\mapsto R(\lambda^{'},B,\varphi)$ is decreasing
in $B>0.$ Consequently, fixing $\lambda^{'}>0\;$ and taking
$\;0<B<B^{'},$
$$\begin{array}{ll}
-\,\displaystyle{\frac{\mu(\lambda^{'},B)}{\lambda^{'}B}}=R(\lambda^{'},B,\psi^{\lambda^{'},B})>R(\lambda^{'},B^{'},\psi^{\lambda^{'},B})&\geq
\displaystyle{\,\min_{\underset{\displaystyle{||\varphi||_{L^{^{2}}(C)}=1}}{\displaystyle{\varphi\in
    H^{1}(C)\setminus\{0\};}}}R(\lambda^{'},B^{'},\varphi)}\vspace{3 pt}\\
&=-\,\displaystyle{\frac{\mu(\lambda^{'},B^{'})}{\lambda^{'}B^{'}}}.
    \end{array}$$

In other words, for all $\lambda^{'}>0,$ the function
$\displaystyle{B\mapsto {\mu(\lambda^{'},B)}/{\lambda^{'}B}}$ is
increasing in $B>0.$  Now, we take randomly $0<B<B^{'}.$ Thus,
\begin{eqnarray}
\begin{array}{rl}
\nonumber\displaystyle{\frac{c_{\Omega,
A,\sqrt{B^{\,'}}\,q,B^{\,'}f}^*(e)}{\sqrt{B^{'}}}}&=\displaystyle{\min_{{\lambda}^{'}>0}\frac{\mu(\lambda^{'},B^{'})}{\lambda^{'}B^{'}}}=\displaystyle{\frac{\displaystyle{\mu(\lambda^{'}_{B^{'}},B^{'})}}{\lambda^{'}_{B^{'}}\times B^{'}}}\\
\nonumber&>\displaystyle{\frac{\mu(\lambda^{'}_{B^{'}},B)}{\lambda^{'}_{B^{'}}\times
B}}\geq\displaystyle{\min_{{\lambda}^{'}>0}\frac{\mu(\lambda^{'},B)}{\lambda^{'}B}}=
\displaystyle{\frac{c_{\Omega, A,\sqrt{B}\,q,B f}^*(e)}{\sqrt{B}},}
\end{array}
\end{eqnarray}
which means that $B\mapsto \displaystyle{{c_{\Omega,
A,\sqrt{B}\,q,Bf}^*(e)}/{\sqrt{B}}}$ is increasing in $B>0.$\hfill$\Box$

\section{Applications to homogenization problems}\label{homogenization}

The reaction-advection-diffusion problem set in a heterogenous
periodic domain $\Omega$ satisfying (\ref{comega}) generates a
homogenization problem:

Let $e\in\mathbb{R}^{d}$ be a vector of unit norm. Assume that
$\Omega,\,A,\,q,\,$ and $f$ are $\,(L_1,\ldots,L_d)-$ periodic and
that they satisfy (\ref{comega}), (\ref{cA}), (\ref{cq}),
(\ref{cf1}) and (\ref{cf2}).

For each $\varepsilon>0,\;$ let
$\displaystyle{\Omega^{\varepsilon}=\varepsilon\,\Omega}\;$ and
consider the following re-scales:
$$\forall (x,y)\in\Omega^{\varepsilon},\quad \displaystyle{A_{\varepsilon}(x,y)}=
\displaystyle{
A\left(\frac{x}{\varepsilon},\frac{y}{\varepsilon}\right),}\;\displaystyle{q_{_{\displaystyle{\varepsilon}}}(x,y)=q\left(\frac{x}{\varepsilon},\frac{y}{\varepsilon}\right),}~~
\hbox{and}~~
\displaystyle{f_{_{\displaystyle{\varepsilon}}}(x,y)=f\left(\frac{x}{\varepsilon},\frac{y}{\varepsilon}\right)}.$$

The coefficients
$\displaystyle{A_{\varepsilon},\,q_{_{\displaystyle{{\varepsilon}}}}},\;$
and $\displaystyle{f_{_{\displaystyle{{\varepsilon}}}}}$ together
with the domain $\displaystyle{\Omega^{\varepsilon}}$ are
$\displaystyle{(\varepsilon\,L_{1},\ldots,\varepsilon\,L_{d})}-$periodic,
and they satisfy similar properties to those of $A,\,q,\,f$ and
$\Omega.$

Consider the parametric reaction-advection-diffusion problem
\begin{eqnarray*}
                                                              \left(P_{\varepsilon}\right)\left\{
                                                                \begin{array}{ll}
 \displaystyle{u^{\varepsilon}_t(t,x,y) =\displaystyle{\nabla\cdot(A_{\varepsilon}\nabla u^{\varepsilon})(t,x,y)\;+\,q_{\varepsilon}\cdot\,\nabla u^{\varepsilon}\,+\,f_{\varepsilon}(x,y,u^{\varepsilon}),\; t\in\mathbb{R},\;(x,y)\in\Omega^{\varepsilon},}}
 \\\\
       \displaystyle{ \displaystyle{\nu^{\varepsilon}}\cdot A_{\varepsilon}\;\nabla u^{\varepsilon}(t,x,y) =0,\quad
        t\,\in\,\mathbb{R},\;(x,y)\,\in\,\partial\Omega^{\varepsilon},}
                                                                \end{array}
                                                              \right.
                                                            \end{eqnarray*}
where $\displaystyle{\nu^{\varepsilon}(x,y)}$ denotes the outward
unit normal on $\displaystyle{\partial\Omega^{\varepsilon}}$ at the
point $(x,y).$

Owing to the results found by Berestycki and Hamel in section $6$ of
\cite{BH1}, and since the coefficients
$A_{\varepsilon},\,\displaystyle{f_{_{\displaystyle{\varepsilon}}}}\;$
and $\;\displaystyle{q_{_{\varepsilon}}}\;$ together with the domain
$\Omega^{\varepsilon}$ satisfy all the necessary assumptions, it
follows that the problem $(P_{\varepsilon})$ admits a minimal speed
of propagation $c_{\Omega^{\varepsilon},\, A_{\varepsilon},\,
q_{\varepsilon},f_{\varepsilon}}^{*}(e)>0\;$ such that
$(P_{\varepsilon})$ has a solution $u^{\varepsilon}$ in the form of
a pulsating front within a speed $c$ if and only if $c\geq
c_{\Omega^{\varepsilon},\, A_{\varepsilon},\,
q_{\varepsilon},f_{\varepsilon}}^{*}(e)>0.$

In this section, we investigate the limit of the parametric minimal
speeds $c_{\Omega^{\varepsilon},\, A_{\varepsilon},\,
q_{\varepsilon},f_{\varepsilon}}^{*}(e)$  (whose parameter is
$\varepsilon$) of the problems
$\displaystyle{(P_{\varepsilon})_{\varepsilon>0}}$ as
$\varepsilon\rightarrow0^{+}.$ In other words, we search the limit
of these minimal speeds as the periodicity cell
$\;C^{\varepsilon}=\varepsilon\,C\;$ becomes a very small size.
On the other hand, we study although not the most general setting, the
variation of the map $\varepsilon\mapsto
\,c_{\Omega^{\varepsilon},\, A_{\varepsilon},\,
q_{\varepsilon},f_{\varepsilon}}^{*}(e)\;$ in $\varepsilon>0.$
\begin{thm}\label{hom with epsilon}
Let $e\in\mathbb{R}^{d}$ be a unit vector, and let
$\;\Omega\subseteq\mathbb{R}^{N}\;$ be a domain which is
$L-$periodic and satisfying (\ref{comega}). Assume that
$A=A(x,y),\;q=q(x,y),\;$ and $\;f=f(x,y,u)$ are $L-$periodic and
that they satisfy (\ref{cA}), (\ref{cq}), (\ref{cf1}) and
(\ref{cf2}) together with the assumptions
$\nabla.A\tilde{e}\equiv0\;\hbox{on}\;\overline{\Omega}\;\hbox{and}\;\nu.A\tilde{e}=0\;$
on $\partial\Omega.$ For each $\varepsilon
>0,$ consider the problem
\begin{eqnarray}\label{hom with u}
                                                            \left\{
                                                                \begin{array}{ll}
 \displaystyle{u^{\varepsilon}_t(t,x,y) =\displaystyle{\nabla\cdot(A_{\varepsilon}\nabla u^{\varepsilon})(t,x,y)\;+\,q_{\varepsilon}}\cdot\nabla u^{\varepsilon}\,+\,f_{\varepsilon}}(x,y,u^{\varepsilon}),\;t\in\mathbb{R},\;(x,y)\in\Omega^{\varepsilon},\\
       \displaystyle{
\displaystyle{\nu^{\varepsilon}}\cdot A_{\varepsilon}\;\nabla
u^{\varepsilon}(t,x,y) =0,\quad
        t\,\in\,\mathbb{R},\;(x,y)\,\in\,\partial\Omega^{\varepsilon},}
                                                                \end{array}
                                                              \right.
                                                            \end{eqnarray}
where $A_{\varepsilon},\,\displaystyle{f_{\varepsilon}}\,$ and
$\displaystyle{q_{\varepsilon}}\;$ are the coefficients defined in
the beginning of this section. Then, the minimal speed
$\;\displaystyle{c\,_{\displaystyle{\Omega^{\varepsilon},\,
A_{\varepsilon},\, q_{\varepsilon},f_{\varepsilon}}}^{*}(e)}$ of
pulsating travelling fronts propagating in the direction of $e$ and
solving (\ref{hom with u}) satisfies
\begin{equation}\label{homogenized speed}
\displaystyle{\lim_{\varepsilon\rightarrow0^{+}}\displaystyle{c\,_{\displaystyle{\Omega^{\varepsilon},\,
A_{\varepsilon},\,
q_{\varepsilon},f_{\varepsilon}}}^{*}(e)=}2\sqrt{\mi_{\!\!\!\!_C}\tilde{e}A\tilde{e}(x,y)dx\,dy}\,\sqrt{\mi_{\!\!\!\!_C}\zeta(x,y)dx\,dy}},
\end{equation}
 where $C$ is the periodicity cell of $\Omega$ and $\tilde{e}=(e,0,\cdots,0)\in\mathbb{R}^{N}.$
\end{thm}

\textbf{Proof.} As a first notice, we mention that the domain
$\Omega^{\varepsilon}\,$ is the image of $\Omega$ by the a dilation
whose center is the origin $O(0,\ldots,0)$ and whose scale factor is
equal to $\varepsilon.$ Consequently,
$$\hbox{for each}\;\varepsilon>0,\;\displaystyle{(\varepsilon x,\,\varepsilon y)\in\Omega^{\varepsilon}\quad\hbox{if and only if}\quad(x,y)\in\Omega,\;\hbox{and}}$$
$$\displaystyle{(\varepsilon x,\,\varepsilon y)\in\partial\Omega^{\varepsilon}\quad\hbox{if and only if}\quad(x,y)\in\partial\Omega.}$$

Moreover,
$$\displaystyle{\forall\varepsilon>0,\;\forall(x,y)\in\partial\Omega,\;\nu^{\varepsilon}(\varepsilon x,\,\varepsilon y)=\nu(x,y)}.$$

Consider now, for each $\varepsilon>0,$ the following change of
variables
$$v^{\varepsilon}(t,x,y)=u^{\varepsilon}(t,\varepsilon x,\varepsilon y)\,; \quad (t,x,y)\in\mathbb{R}\times\Omega.$$

One gets
$$\forall(t,x,y)\in\mathbb{R}\times\Omega,\;v^{\varepsilon}_t(t,x,y)=u^{\varepsilon}_t(t,\varepsilon x,\varepsilon y),$$
$$\displaystyle{\nabla_{x,y}\cdot(A(x,y)\nabla v^{\varepsilon})(t,x, y)= \nabla_{x,y}\cdot(A_{\varepsilon}\nabla u^{\varepsilon})(t,\varepsilon x,\varepsilon y)=\varepsilon^{2}\,\nabla\cdot(A_{\varepsilon}\nabla u^{\varepsilon})(t,\varepsilon x,\varepsilon y)},$$
and
\begin{eqnarray}\label{boundaryv}
    \begin{array}{ll}
\nu_{\varepsilon}(\varepsilon x,\varepsilon
y)\cdot\left[A_{\varepsilon}\nabla
u^{\varepsilon}\right](t,\varepsilon x,\varepsilon
y)&=\displaystyle{\nu(x,y)\cdot A\left(\frac{\varepsilon
x}{\varepsilon},\frac{\varepsilon y}{\varepsilon}\right)\nabla
u^{\varepsilon}(t,\varepsilon x,\varepsilon y)} \\
 &=\displaystyle{\frac{1}{\varepsilon}\;\nu(x,y)\cdot A(x,y)\nabla
v^{\varepsilon}(t,x,y)~\hbox{on}~\mathbb{R}\times\partial\Omega.}
\end{array}
\end{eqnarray}
The boundary condition in (\ref{hom with u}) yields that
$\nu_{\varepsilon}(\varepsilon x,\varepsilon y)\cdot
\left[A_{\varepsilon}\nabla u^{\varepsilon}\right](t,\varepsilon
x,\varepsilon y)=0,\,$ for all $\;(t,x,y)\in
\mathbb{R}\times\partial\Omega$ (which is equivalent to say: for all
$(t,\varepsilon x,\varepsilon
y)\in\mathbb{R}\times\partial\Omega^{\varepsilon}).$ It follows from
(\ref{boundaryv}) that
$$\forall(t,x,y)\in\mathbb{R}\times\partial\Omega,\quad \nu\cdot A\nabla v^{\varepsilon}(t,x,y)=0.$$

One can now conclude  that: for each
$\varepsilon>0,\;u^{\varepsilon}$ satisfies (\ref{hom with u}) if
and only if $v^{\varepsilon}\;$ satisfies
\begin{eqnarray}\label{hom with v epsilon}
                                                            \left\{
                                                                \begin{array}{ll}
 \displaystyle{v^{\varepsilon}_t(t,x,y) =\displaystyle{\frac{1}{\varepsilon^{2}}\nabla\cdot(A\nabla v^{\varepsilon})(t,x,y)\;+\,\frac{1}{\varepsilon}\,q\cdot\nabla v^{\varepsilon}\,+\,f(x,y,v^{\varepsilon})},\; t\in\mathbb{R},\;(x,y)\in\Omega},\vspace{3 pt}\\
\displaystyle{\nu\cdot A\;\nabla v^{\varepsilon}(t,x,y) =0,\;
        t\,\in\,\mathbb{R},\;(x,y)\,\in\,\partial\Omega.}
\end{array}
  \right.
            \end{eqnarray}

Having the assumptions (\ref{comega}), (\ref{cA}), (\ref{cq}),
(\ref{cf1}), and (\ref{cf2}) on $\Omega,\,A,\,q,\,$ and $f,$ one
gets that problem (\ref{hom with v epsilon}) admits, for each
$\varepsilon>0,$ a minimal speed of propagation denoted by
$\displaystyle{c_{\Omega,\,\left(\frac{1}{{\varepsilon}}\right)^{2}
A,\, \frac{1}{\varepsilon}\,q,\,f}^{*}(e)}.$ \vskip0.2cm

 Moreover,
due to the change of variables between $u^{\varepsilon}$ and
$v^{\varepsilon},$ it follows that for each
$\varepsilon>0,\;u^{\varepsilon}$ is a pulsating travelling front
propagating in the direction of $e$ within a speed $c$ and solving
(\ref{hom with u}) if and only if $v^{\varepsilon}$ is a pulsating
travelling front propagating in the direction of $e$ within a speed
$\;\displaystyle{\frac {c}{\varepsilon}}$ and solving (\ref{hom with
v epsilon}). This yields that
\begin{eqnarray}\label{relations between min speeds}
  \begin{array}{ll}
\forall \varepsilon>0, \quad
\displaystyle{\displaystyle{c\,_{\displaystyle{\Omega^{\varepsilon},\,
A_{\varepsilon},\,
q_{\varepsilon},f_{\varepsilon}}}^{*}(e)}}&=\displaystyle{\varepsilon
\,c_{\Omega,\,\left(\frac{1}{{\varepsilon}}\right)^{2} A,\,
\frac{1}{\varepsilon}\,q,\,f}^{*}(e)}
=\displaystyle{{c_{\Omega,\,M A,\,
\sqrt{M}\,q,\,f}^{*}(e)}/{\displaystyle{\sqrt{M}}}},
 \end{array}
\end{eqnarray}
where $M=\displaystyle{\left(\,{1}/{\displaystyle\varepsilon}\,\right)^{2}}.$

As $\varepsilon\rightarrow0^+,$ the variable $M\rightarrow+\infty.$
Applying Theorem \ref{lim as M}, with
$\displaystyle{\gamma=\frac{1}{2}},$ one gets that
$$\displaystyle{\lim_{M\rightarrow+\infty}\displaystyle{\frac{c_{\Omega,\,M A,\,
\sqrt{M}\,q,\,f}^{*}(e)}{\displaystyle{\sqrt{M}}}}=2\sqrt{\mi_{\!\!\!\!_C}\tilde{e}A\tilde{e}(x,y)dx\,dy}\,\sqrt{\mi_{\!\!\!\!_C}\zeta(x,y)dx\,dy}}.$$

Therefore,
$\displaystyle{\lim_{\varepsilon\rightarrow0^+}\displaystyle{c\,_{\displaystyle{\Omega^{\varepsilon},\,
A_{\varepsilon},\,
q_{\varepsilon},f_{\varepsilon}}}^{*}(e)}=2\sqrt{\mi_{\!\!\!\!_C}\tilde{e}A\tilde{e}(x,y)dx\,dy}\,\sqrt{\mi_{\!\!\!\!_C}\zeta(x,y)dx\,dy}},$
 and the proof of Theorem \ref{hom with
epsilon} is complete.\hfill$\Box$

\begin{remark}
It is worth noticing that, in formula \ref{homogenized speed}, the
homogenized speed depends on the averages of the diffusion and
reaction coefficients, but it does not depend on the advection.
\end{remark}

We move now to study the variation of the map
$\displaystyle{\varepsilon\mapsto\displaystyle{c\,_{\displaystyle{\Omega^{\varepsilon},\,
A_{\varepsilon},\, q_{\varepsilon},f_{\varepsilon}}}^{*}(e)}}$ with
respect to $\varepsilon>0.$ In other words, we want to check the
monotonicity behavior of the parametric minimal speed of
propagation, whose parameter $\varepsilon>0,$ as the periodicity
cell of the domain of propagation shrinks or enlarges within a ratio
$\varepsilon.$ In this study, we will consider the same situation of
Theorem \ref{variation of min speed with respect diffusion} and also
the same notations introduced in the beginning of section
\ref{homogenization}:

\begin{thm}\label{var with respect to the size of the cell}
Let $e=(1,0\ldots,0).$ Assume that $\Omega$ has the form
$\mathbb{R}\times\omega$ where $\omega$ may or may not be bounded
(precisely described in section \ref{as eps tends to zero}) and that
the diffusion matrix $A = A(y)$ satisfies (\ref{cA Y}) together with
the assumption that $e$ is an eigenvector of $A(y)$ for all
$y\in\overline{\omega},$ that is
\begin{equation}\label{Ae= alpha(y)e 1}
    A(x,y)e=A(y)e=\alpha(y)e,\;\hbox{for all}\; (x,y)\in\mathbb{R}\times\overline{\omega};
\end{equation}
where $y\mapsto\alpha(y)$ is a positive $(L_1,\ldots,L_d)-$ periodic
function defined over $\overline{\omega}.$ The nonlinearity $f$ is
assumed to satisfy (\ref{nonlinearity 1}) and (\ref{nonlinearity
2}). Assume further more that the advection field $q$ (when it
exists) is in the form $q(x,y)=(q_{_{1}}(y),0,\ldots,0)$ where
$q_{_{1}}$ has a zero average over $C,$ the periodicity cell of
$\omega.$ For $\varepsilon>0$ consider the
reaction-advection-diffusion problem
\begin{eqnarray}\label{hom with u 1}
                                                            \left\{
                                                                \begin{array}{ll}
                                                                \forall \,t\in\mathbb{R},\;\forall \,(x,y)\in\Omega^{\varepsilon}=\mathbb{R}\times\varepsilon\,\omega,\vspace{3 pt}\\
 \displaystyle{u^{\varepsilon}_t(t,x,y) =\displaystyle{\nabla\cdot(A_{\varepsilon}\nabla u^{\varepsilon})(t,x,y)\;+\,q_{\varepsilon}}\cdot\,\nabla u^{\varepsilon}\,+\,f_{\varepsilon}}(x,y,u^{\varepsilon});
 \vspace{3 pt}\\
       \displaystyle{
\displaystyle{\nu^{\varepsilon}}\cdot A_{\varepsilon}\;\nabla
u^{\varepsilon}(t,x,y) =0,\;
        t\in\mathbb{R},\;(x,y)\in\partial\Omega^{\varepsilon}.}
                                                                \end{array}
                                                              \right.
                                                            \end{eqnarray}
Then, the map
$\displaystyle{\varepsilon\mapsto\displaystyle{c\,_{\displaystyle{\Omega^{\varepsilon},\,
A_{\varepsilon},\, q_{\varepsilon},f_{\varepsilon}}}^{*}(e)}}$ is
increasing in $\varepsilon>0.$
\end{thm}

\textbf{Proof of Theorem \ref{var with respect to the size of the
cell}.} For each $\varepsilon>0,$ we consider the change of
variables
$$v^{\varepsilon}(t,x,y)=u^{\varepsilon}(t,\varepsilon x,\varepsilon y)\,; \quad (t,x,y)\in\mathbb{R}\times\Omega.$$

Owing to the justifications shown in the proof of Theorem \ref{hom
with epsilon}, one consequently obtains
\begin{eqnarray}\label{relations between min speeds1}
  \begin{array}{ll}
\forall \varepsilon>0, \quad
\displaystyle{\displaystyle{c\,_{\displaystyle{\Omega^{\varepsilon},\,
A_{\varepsilon},\,
q_{\varepsilon},f_{\varepsilon}}}^{*}(e)}}&=\displaystyle{\varepsilon
\,c_{\Omega,\,\left(\frac{1}{{\varepsilon}}\right)^{2} A,\,
\frac{1}{\varepsilon}\,q,\,f}^{*}(e)}=\displaystyle{{c_{\Omega,\,\beta A,\,
\sqrt{\beta}\,q,\,f}^{*}(e)}/{\displaystyle{\sqrt{\beta}}}},
 \end{array}
\end{eqnarray}
where $\beta(\varepsilon)=\displaystyle{\left(\,{1}/{\displaystyle\varepsilon}\,\right)^{2}}.$

Applying Theorem \ref{variation of min speed
with respect diffusion}, it follows that the map
$\displaystyle{\eta_{_1}:\;\beta\mapsto{c_{\Omega,\,\beta A,\,
\sqrt{\beta}\,q,\,f}^{*}(e)}/{\displaystyle{\sqrt{\beta}}}}$ is
decreasing in $\beta>0.$ On the other hand, the map
$\eta_{_2}:\;\displaystyle{\varepsilon\mapsto \beta(\varepsilon)}$
is also decreasing in $\varepsilon>0.$ Therefore,
$\displaystyle{\varepsilon\mapsto\displaystyle{c\,_{\displaystyle{\Omega^{\varepsilon},\,
A_{\varepsilon},\, q_{\varepsilon},f_{\varepsilon}}}^{*}(e)}},$
which is the composition $\displaystyle{\eta_{_1}\circ\eta_{_2}},$
is increasing in $\varepsilon>0$ and this completes our proof.\hfill
$\Box$\vskip0.3cm

Other homogenization results, concerning
reaction-advection-diffusion problems, were given in the case of a
combustion-type nonlinearity $f=f(u)$ satisfying
\begin{eqnarray}\label{ignition temp , homogen}
\left\{
  \begin{array}{ll}
   \exists \;\theta\,\in\,(0,1),\; f(s)=0\;\hbox{for all}\;s\in\left[0,\theta\right],\;
    f(s)>0\;\hbox{for all}\;s\in(\theta,1),\;f(1)=0,\\
    \exists\rho\in(0,1-\theta),\quad f \;\hbox{is non-increasing on}
    \;\left[1-\rho,1\right].
  \end{array}
\right.
\end{eqnarray}

Consider the equation
\begin{equation}\label{u eps in RN}
\displaystyle{u^{\varepsilon}_t(t,x)
=\displaystyle{\nabla\cdot(A(\varepsilon^{-1}\,x)\nabla
u^{\varepsilon})\;+\,\varepsilon^{-1}q(\varepsilon^{-1}\,x)\cdot\,\nabla
u^{\varepsilon}\,+\,f}(u^{\varepsilon})\quad\hbox{in}\;\mathbb{R}^{N},}
\end{equation}
where the nonlinearity $f$ satisfies (\ref{ignition temp ,
homogen}), and the drift and diffusion coefficients $q$ and $A$
satisfy the general assumptions (\ref{cA}) and (\ref{cq}), with
periodicity $1$ in all variables $x_{1},\ldots,x_{N}.$ Fix a unit
vector $e$ of $\mathbb{R}^{N}.$ From Berestycki and Hamel
\cite{BH1}, it follows that for each $\varepsilon>0,$ problem (\ref{u
eps in RN}) admits a unique pulsating front
$\left(c_{\varepsilon},u^{\varepsilon}\right)$ such that
$$u^{\varepsilon}(t,x)=\phi^{\varepsilon}(x\cdot e+c_\varepsilon t,\,x)$$
where $\phi^{\varepsilon}(s,x)$ is
$(\varepsilon,\ldots,\varepsilon)-$periodic in $x$ that satisfies
$\phi^{\varepsilon}(-\infty,.)=0$ and
$\phi^{\varepsilon}(+\infty,.)=1.$ The functions $u^{\varepsilon}$
are actually unique up to shifts in time, and one can assume that
$\displaystyle{\max_{\mathbb{R}^{N}}\phi^{\varepsilon}(0,.)=\theta}.$

Concerning problem (\ref{u eps in RN}), Heinze
\cite{Heinze Hom} proved that
$$\displaystyle{\hbox{as }\;\varepsilon\rightarrow0^{+},\;c_{\varepsilon}\rightarrow c_{0}>0,\;\hbox{and}\;u^{\varepsilon}(t,x)\rightarrow u_{0}(x\cdot e+c_{0}t)\;\hbox{weakly in}\;H^{1}_{loc},}$$
where $(c_{0},u_{0})$ is the unique solution of the one-dimensional
homogenized equation
\begin{eqnarray}
\left\{
  \begin{array}{ll}
    a^{*}\,u_{0}^{''}-c_{0}u_{0}^{'}+f(u_{0})=0\;\hbox{in}\;\mathbb{R},\\
u_{0}(-\infty)=0<u_{0}<u_{0}(+\infty)=1\;\hbox{in}\;\mathbb{R},\;u_{0}(0)=\theta
  \end{array}
\right.
\end{eqnarray}
and $a^{*}$ is a positive constant determined in \cite{Heinze Hom}.

In Theorem 1 of Caffarelli, Lee, Mellet \cite{CaffLeeMellet1}, the homogenization limit was combined with
the singular high activation limit for the reaction (one can also see \cite{CaffLeeMellet2} in this context) while the diffusion matrix was taken $\displaystyle{A=Id_{\mathbb{R}^N}}$. More precisely, the
nonlinearity had the form $\displaystyle{f_\varepsilon(u)=\frac{1}{\varepsilon}\beta(\frac{u}{\varepsilon})}$ with $\beta(s)$ a Lipschitz fucntion
satisfying $$\beta(s)>0\hbox{ in }(0,1)\hbox{ and }\beta(s)=0\hbox{ otherwise.}$$
  These nonlinearities approach a Dirac mass at $u=1.$
\section{Open problems} In all the results of this paper, we deal with nonlinearities of the ``KPP'' type. In the periodic framework of this paper,
 pulsating travelling fronts exist also with other types of nonlinearities (see Theorems 1.13 and 1.14 in \cite{BH1}). Namely, they exist when $f=f(x,y,u)$ is of the ``combustion'' type satisfying:

\begin{eqnarray}\label{f both}
    \left\{
      \begin{array}{ll}
f \;\hbox{is globally Lipschitz-continuous in}\,
\;\overline{\Omega}\times\mathbb{R},\\
\displaystyle{\forall\,(x,y)\in\,\overline{\Omega},\,\forall\,s\in(-\infty,0]\cup[1,+\infty),\,f(s,x,y)=0,}\\
 \exists
\,\rho\in\,(0,1),\;\forall(x,y)\,\in\overline{\Omega},\;\displaystyle{\forall\,
1-\rho\leq\,s\,\leq s'\,\leq\,1,}\;
\displaystyle{f(x,y,s)\;\geq\,f(x,y,s')},
  \end{array}
    \right.
\end{eqnarray}
and

  \begin{eqnarray}\label{combustion m}
\left\{
\begin{array}{ll}
f~ \hbox{ is $L-$periodic with respect to $x,$}\\
\exists\, \theta\,\in\,(0,1),\;
\forall(x,y)\in\overline{\Omega},\;\forall\, s\in[0,\theta],\;
f(x,y,s)=0,\\
\forall\,s\in\,(\theta,1),\; \exists
\,(x,y)\in\overline{\Omega}\;\;\hbox{such that}\;f(x,y,s)\,>\,0
\hbox{,}
\end{array}
\right.
\end{eqnarray}
or when $f=f(x,y,u)$ is of the ``ZFK'' (for Zeldovich-Frank- Kamenetskii) type satisfying (\ref{f both}) and

\begin{eqnarray}\label{ZFK m}
    \left\{
      \begin{array}{ll}
     f~ \hbox{ is $L-$periodic with respect to $x,$}\\
\exists\delta>0,\;\hbox{the restriction of $f$ to $\overline{\Omega}\,\times\,[0,1]$ is of class}\;C^{1,\,\delta},  \\
 \forall\,s\in\,(0,1),\; \exists \,(x,y)\in\overline{\Omega}\;\;\hbox{such that}\;f(x,y,s)\,>\,0  \hbox{.} \\
\end{array}
    \right.
\end{eqnarray}
In particular, the ``KPP'' nonlinearities are of the ``ZFK'' type.

Recently, El Smaily \cite{El Smaily min max} gave $\min-\max$ and $\max-\min$ formul{\ae} for the speeds of propagation of problem (\ref{front}) taken with a ``ZFK'' or a
``combustion'' nonlinearity. These formul{\ae}, together with the results of this paper, can give important estimates
 for the parametric minimal speeds of
the problem (\ref{front}) when $f$ is a ``ZFK'' nonlinearity which is not of the ``KPP'' type. Indeed, if $f$ is a ``ZFK'' nonlinearity, one can find a
 ``KPP'' function $h=h(x,y,u)$ such that $$\forall(x,y,u)\in\overline{\Omega}\times\mathbb{R},~ f(x,y,u)\leq h(x,y,u).$$

Referring to formula (1.17) in El Smaily \cite{El Smaily min max}, one can conclude that
$$\forall M>0,\forall B>0,~\forall \gamma\in\mathbb{R},~~\ds{c_{\Omega,MA,\displaystyle{M^{\gamma}\,q},Bf}^*(e)\leq c_{\Omega,MA,\displaystyle{M^{\gamma}\,q},Bh}^*(e).}$$
Moreover, if $f$ is a ``ZFK'' nonlinearity satisfying the additional assumption
\begin{equation}\label{add cond}
\forall (x,y)\in\o,~f'_u(x,y,0)>0,
\end{equation}
then one can find a ``KPP'' function $g=g(x,y,u)$ such that $g\leq f$ in $\o\times\R,$ and thus
\begin{equation}\label{estimating ZFK by KPP}
\begin{array}{c}
\forall M>0,\forall B>0,~\forall \gamma\in\mathbb{R},\vspace{5 pt}\\
 c_{\Omega,MA,\displaystyle{M^{\gamma}\,q},Bg}^*(e)\leq c_{\Omega,MA,\displaystyle{M^{\gamma}\,q},Bf}^*(e)\leq c_{\Omega,MA,\displaystyle{M^{\gamma}\,q},Bh}^*(e).
\end{array}
\end{equation}
As a consequence, under the assumptions that $0\leq\gamma\leq 1/2,$ $\nu\cdot A\tilde e=0$ on $\partial\Omega,$ and $\nabla\cdot A\tilde e\equiv0$ in $\Omega,$
Theorem \ref{lim as M} implies that
\begin{equation}\label{upper estimate}
\displaystyle{\limsup_{M\rightarrow+\infty}\displaystyle{\frac{c_{\Omega,\,M A,\,
M^{\gamma}q,\,f}^{*}(e)}{\displaystyle{\sqrt{M}}}}\leq2\sqrt{\mi_{\!\!\!\!_C}\tilde{e}A\tilde{e}(x,y)dx\,dy}\,\sqrt{\mi_{\!\!\!\!_C}g'_u(x,y,0)dx\,dy}},
\end{equation}
and
\begin{equation}\label{lower estimate}
\displaystyle{\liminf_{M\rightarrow+\infty}\displaystyle{\frac{c_{\Omega,\,M A,\,
M^{\gamma}q,\,f}^{*}(e)}{\displaystyle{\sqrt{M}}}}=2\sqrt{\mi_{\!\!\!\!_C}\tilde{e}A\tilde{e}(x,y)dx\,dy}\,\sqrt{\mi_{\!\!\!\!_C}h'_u(x,y,0)dx\,dy}}\;>0.
\end{equation}

If $f$ is a ``combustion'' nonlinearity, then problem (\ref{front}) admits a solution $\displaystyle{\left(c,u\right)}$
where $c=\ds{c_{_{\Omega, A,q,f}}(e)>0}$ is unique and $u=u(t,x,y)$ is increasing in $t$ and it is unique up to a translation in $t.$ Taking $g$ as a ``KPP'' nonlinearity such that
$g\geq f$ in $\overline{\Omega}\times\mathbb{R}$ and using Theorem \ref{lim as M}, it follows that
\begin{equation}\label{upper estimate 1}
\begin{array}{l}
\displaystyle{\limsup_{M\rightarrow+\infty}\displaystyle{\frac{c_{_{\Omega,\,M A,\,
M^{\gamma}\,q,\,f}}(e)}{\displaystyle{\sqrt{M}}}}\leq2\sqrt{\mi_{\!\!\!\!_C}\tilde{e}A\tilde{e}(x,y)dx\,dy}\,\sqrt{\mi_{\!\!\!\!_C}g'_u(x,y,0)dx\,dy}}\vspace{5pt}\\
\hbox{ together with }~~
\ds{\liminf_{M\rightarrow+\infty}\displaystyle{\frac{c_{_{\Omega,\,M A,\,
M^{\gamma}\,q,\,f}}(e)}{\displaystyle{\sqrt{M}}}}}\geq0.
\end{array}
\end{equation}

Similarly, one can get several estimates concerning the case of a small diffusion factors, small (resp. large) reaction factors, or small (resp. large)
periodicity parameters.

The above motivation gives several upper and lower estimates for the parametric speeds of propagation. However, the exact limits are not known.
This leads us to ask about the asymptotics of the minimal speeds of propagation with respect to diffusion, reaction and periodicity factors in
the ``ZFK'' case and about the asymptotics of the unique parametric speed of propagation in the ``combustion'' case.
These studies should help, as it was done in section \ref{homogenization}, in solving
some homogenization problems in the ``ZFK'' case.

Besides, Theorem \ref{hom with epsilon} gives the limit of $\displaystyle{c\,_{\displaystyle{\Omega^{\varepsilon},\,
A_{\varepsilon},\, q_{\varepsilon},f_{\varepsilon}}}^{*}(e)}$ as $\epsilon\rightarrow0^+.$ However, finding the homogenized equation of (\ref{hom with u}) in the
``KPP'' remains an open problem.

\section{Conclusions}
As we mentioned in the beginning of this paper, our first aim was to give a complete and rigorous analysis of the \emph{minimal speed of propagation} of pulsating travelling fronts solving \emph{parametric heterogeneous} reaction-advection-diffusion equations in a periodic framework. In the paper of Berestycki, Hamel and Nadirashvili \cite{BHN1}, several upper and lower estimates for the parametric minimal speed of propagation were given (see Theorems 1.6 and 1.10 in \cite{BHN1}). However, the exact asymptotic behaviors of the minimal speed with respect to diffusion and reaction factors and with respect to the periodicity parameter $L$ were not given there. In this paper, we determined the exact asymptotes of the minimal speed
in the ``KPP'' periodic framework. In sections \ref{as eps tends to zero}, \ref{as M goes to infty} and \ref{reaction}, we proved that
 (under some assumptions on $A,$ $q,$  $f$ and $\Omega$) the asymptotes of the parametric minimal speed are either $$ 2\sqrt{\max_{\overline{\omega}}\,\zeta}\sqrt{\max_{\overline{\omega}}eAe} \quad\hbox{ or }\quad \displaystyle{2\sqrt{\mi_{\!\!\!\!_C}\tilde{e}A\tilde{e}(x,y)dx\,dy}\,\sqrt{\mi_{\!\!\!\!_C}\zeta(x,y)dx\,dy}}.$$
 (see Theorems \ref{limit as eps}, \ref{lim as per L tend to infty}, \ref{lim as M}, \ref{lim as per L tend to 0}, \ref{lim as B tend to infty} and \ref{lim as B tend 0} above).
 Moreover, we found in section \ref{as eps tends to zero} that the presence of an advection field, in the general form or in the form of shear flows, changes the asymptotic behavior of the minimal speed within a small diffusion (see Theorem \ref{lim as eps but in presence of a shear flow} and Remark \ref{without shear}). Conversely, we proved in Section 4 that the presence of a general advection field $M^{\gamma}q$ (where $q$ satisfies (\ref{cq})) has no effect on $\displaystyle{\lim_{M\rightarrow+\infty}\frac{c_{\Omega,MA,\displaystyle{M^{\gamma}\,q},f}^*(e)}{\sqrt{M}}}$
 whenever $0\leq\gamma\leq 1/2$ (see Theorem \ref{lim as M}). Furthermore, we studied, in a particular periodic framework, the variations of the maps
 $\displaystyle{\,\beta\mapsto\frac{c_{\Omega,\beta A,\sqrt{\beta}\,q,f}^*(e)}{\sqrt{\beta}}}$ and $\displaystyle{L
\mapsto c_{\mathbb{R}^{N},A_{_L},\,q_{_{L}},f_{_L}}^{*}(e)}$ and $\displaystyle{B\mapsto \frac{c_{\Omega,
A,\sqrt{B}\,q,Bf}^*(e)}{\sqrt{B}}}$ with respect to the positive variables $\beta,$ $L$ and $B$ respectively. Roughly speaking, we found that the first and the third maps have opposite senses of variations (see Theorems \ref{variation of min speed with respect diffusion} and \ref{variation of c*(B)over rad B}). On the other hand, Theorem \ref{variation of the min speed with respect to period} and Theorem \ref{var with respect to the size of the cell} yield that the minimal speed increases when the medium undergoes a dilation whose scale factor is greater than $1.$

The second aim was to find the homogenized ``KPP'' minimal speed. We achieved this goal in section \ref{homogenization} (Theorem \ref{hom with epsilon}) under the assumptions of free divergence on $A(x,y)\te$ and invariance of the domain in the direction $A(x,y)\te.$ This was an application to the results obtained in section \ref{as M goes to infty}. The found homogenized speed should play an important role in finding the homogenized reaction-advection-diffusion equation in the ``KPP'' case. In a forthcoming paper \cite{EHR}, we find also the homogenized speed in the one dimensional case but in a more general setting (in fact, the assumption of divergence free is equivalent to the assumption that the diffusion term $x\mapsto a(x)$ is constant over $\R$ in the case $N=1$).

All the mathematical results obtained in this paper can be applied to study some \emph{spreading} phenomena. Referring to the results of Weinberger \cite{weinberger1}, one can conclude that the \emph{spreading speed} is equal to the ``KPP'' \emph{minimal speed of propagation} in the periodic framework under some assumptions on the initial data $u_0:=u_0(x,y)=u(0,x,y)$ which is defined on a periodic domain $\Omega$ of $\R^{N}.$ In such a setting, all our results can be applied to give rigorous answers on the asymptotic behavior of the \emph{parametric spreading speed} with respect to diffusion and reaction factors and with respect to the periodicity parameter.

\section*{Acknowledgments}

I am very grateful to Professor Fran\c{c}ois Hamel for his valuable
comments, directions and advices. I would like also to thank
Professor Mustapha Jazar for his support and his constant encouragement during
the preparation of this work.

\end{document}